\documentclass[10pt,leqno]{book}
\usepackage{fmkstyle}

\usepackage{makeidx}
\makeindex
\makeglossary
\begin{document}
\title{Sets and Their Sizes}
\author{Fred M. Katz}

\begin{titlepage}%
  \let\footnotesize\small
  \let\footnoterule\relax
  \let \footnote \thanks
    {\Large Fred M. Katz \par}
  \vspace{108pt}
    {\Huge
        \textit{Sets and Their Sizes}
     }
        \par%
  \vspace{18pt}
    {\Large}
  \vspace{60pt}
  \vspace{300pt}

  \par
  {\Large New York, New York}

\end{titlepage}%
  \setcounter{footnote}{0}%
  \global\let\thanks\relax
  \global\let\maketitle\relax

  \global\let\title\relax
  \global\let\author\relax
  \global\let\date\relax
  \global\let\and\relax
\renewcommand{\thepage}{}
\   \\
Fred M. Katz (fred@katz.com)   \\
Logic \& Light, Inc.
\\
 New York, NY 10028
\\
 USA
 \vspace{300pt}
\  \\
{\Large Sets and Their Sizes}    
\ \\
\ \\
\ \\
\copyright 1981, 2001 Fred M. Katz.
\\
All rights reserved.
This work may not be translated or copied in whole or in part
without the written permission of the publisher (Logic \& Light, Inc.
400 East 84th Street, New York, NY 10028, USA), except for
brief excerpts in connection with reviews or scholarly analysis.
\ \\
\\
    The author hereby grants to M.I.T. permission
    to reproduce and distribute copies of this
    thesis document in whole or in part.
\frontmatter
\Prefatory{Preface}
{
\renewcommand\thesubsection{}
Cantor's theory of cardinality violates common sense.
It says, for example, that all infinite sets of integers
are the same size.
This thesis criticizes the arguments for Cantor's theory
and presents an alternative.
\par
The alternative is based on a general theory,
\CS (for \readAs{Class Size}).
\CS consists of all sentences in the first order language
with a subset predicate and a less-than predicate which
are true in all interpretations of that language whose
domain is a finite power set.
Thus, \CS says that \notion{less than} is a linear
ordering with highest and lowest members and that every
set is larger than any of its proper subsets.
Because the language of \CS is so restricted, \CS will
have infinite interpretations.
In particular, the notion of \notion{one-one correspondence}
cannot be expressed in this language, so Cantor's definition
of similarity will not be in \CS, even though it is true for
all finite sets.
\par
We show that \CS is decidable but not finitely axiomatizable
by characterizing the complete extensions of \CS.
\CS has \notion{finite completions}, which are true only in finite
models and \notion{infinite completions},
which are true only in infinite models.
An infinite completion is determined by a set of \notion{remainder
principles}, which say, for each natural number, \n, how many
atoms remain when the universe is partitioned into \n disjoint
subsets of the same size.
\par
We show that any infinite completion of \CS has a model over
the power set of the natural numbers which satisfies an
additional axiom:
\namedSentence{\OP}
If initial segments of \A eventually become smaller than the
corresponding initial segments of \B, then \A is smaller than
\B.
\endNamedSentence
Models which satisfy \OP seem to accord with common intuitions
about set size.
In particular, they agree with the ordering suggested by the notion
of \notion{asymptotic density}.

 \Subsection{P1}{Acknowledgements}
This dissertation has been made possible by Cantor --- who
invented set theory, Tarski --- who discovered model theory,
Abraham Robinson --- who introduced model-completeness, and
George Boolos --- who made all of these subjects accessible to
me.
Prof. Boolos read several handwritten versions and as many typed
versions of this thesis, with care and patience each time around.
I owe him thanks for finding (and sometimes correcting) my
errors, for directing me to some useful literature, and
simply for being interested.
The reader should join me in thanking him for his help with
style and organization, and for his insistence on filling in
\scare{mere details}.
Nevertheless, whatever errors, gaps, or infelicities remain
are solely the responsibility of the author.
 \Subsection{P1}{Note}

Except for some minor wording changes, typographical corrections,
and  omissions of some elementary (but very long) proofs,
this is the text of my dissertation.

\par

This dissertation was
    submitted in partial fulfillment of the
    requirements for the degree of Doctor of
    Philosophy at the Massachusetts Institute of
    Technology in
    April, 1981

\par

    It was certified and accepted by Prof. George Boolos,
    Thesis Supervisor and Chairperson of the Department
    of Linguistics and Philosophy.

\ \\
\ \\
\ \\
Fred M. Katz
\\ June, 2001
 
}
\tableofcontents

\mainmatter
\Chapter{1}{Introduction}
\Section{1.1}{The Problem}

This paper proposes a theory of set size which is based
on intuitions, naive and otherwise.
The theory goes beyond intuitions, as theories will, so it needs both
justification and defense.
I spend very little time justifying the theory; it is so clearly true
that anyone who comes to the matter without prejudice will accept it.
I spend a lot of time defending the theory because no one who comes
to the matter comes without prejudice.
\par
The prejudice stems from Cantor's theory of set size, which is as old
as sets themselves and so widely held as to be worthy of the name
\emph{the standard theory}.
Cantor's theory consists of just two principles:
\namedSentence{\OO}
Two sets are the same size just in case there is a one-one
correspondence between them.
\endNamedSentence

\namedSentence{\CANTORLT}
A set, $x$, is smaller than a set, $y$, just in case $x$ is the
same size as some subset of $y$, but not the same size as $y$ itself.
\endNamedSentence
A \introSTerm{one-one correspondence} between two sets is a relation
which pairs each member of either set with exactly one member of the
other.
For example, the upper-case letters of the alphabet can be paired
with the lower-case letters:
\begin{center}
{\ttfamily
\par
ABCDEFGHIJKLMNOPQRSTUVWXYZ
\par
abcdefghijklmnopqrstuvwxyz
}
\end{center}
\par
So, the standard theory says, the set of upper-case letters
is the same size as the set of lower-case letters.
Fine and good.
\par
The standard theory also says that the set of even numbers is the
same size as the set of integers since these two sets can also be
paired off one-to-one:
\setcounter{MaxMatrixCols}{20}
\[
\begin{matrix}
\dots & -n & \dots & -3  & -2  & -1 & 0 & 1 & 2 & 3 & \dots & n & \dots\\
\dots & -2n & \dots & -6  & -4  & -2 & 0 & 2 & 4 & 6 & \dots & 2n & \dots
\end{matrix}
\]
\par
Similarly, the standard theory says that the set of positive even
integers is the same size as the set of prime numbers:
pair the n-th prime with the n-th positive even number.
In both of these cases, common sense chokes on the standard theory.
\par
In the first case, common sense holds that the set of integers is larger
than the set of even integers.
The integers contain all of the even integers and then some.
So it's just good common sense to believe there are more of the former
than the latter.
This is just to say that common sense seems to follow:
\namedSentence{\SUBSET}
If one set properly includes another, then the first is larger
than the second.
\endNamedSentence
even into the infinite, where it comes up against the standard theory.
\par
Common sense can make decisions without help from \SUBSET.
Though the set of primes is not contained in the set of even integers,
it is still clear to common sense that the former is smaller than
the latter.
One out of every two integers is even, while the prime numbers are
few and far between.
No doubt, to use this reasoning, you need a little number theory in
addition to common sense; but, given the number theory, it's the only
conclusion common sense allows.
\par
The theory proposed here accommodates these bits of common sense
reasoning.
It maintains \SUBSET and a \emph{few and far between} principle and much
else besides.
To state this theory, we use three two-place predicates:
\Symbol{<}, \Symbol{\eqsize}, and
\Symbol{>}.
If \A and \B name sets, then
\begin{itemize}
\item \schema{$A < B$} is read as \readAs{$A$ is smaller than $B$},
\item \schema{$A \eqsize B$} is read as \readAs{$A$ is the same size as $B$}, and
\item \schema{$A > B$} is read as \readAs{$A$ is larger than $B$}.
\end{itemize}
\par
Incidentally, we assume throughout this thesis that the following
schemata are equivalent, item by item, to the readings of the
three predicates given above, assuming that $A$ is the set of $\alpha$'s
and $B$ is the set of $\beta$'s:
\begin{enumerate}
\item
    \begin{itemize}
    \item There are fewer $\alpha$'s than $\beta$'s.
    \item There are just as many $\alpha$'s than $\beta$'s.
    \item There are more $\alpha$'s than $\beta$'s.
    \end{itemize}
\item
    \begin{itemize}
    \item The number of $\alpha$'s is less than the number of $\beta$'s.
    \item The number of $\alpha$'s is the same as the number of $\beta$'s.
    \item The number of $\alpha$'s is greater than the number of $\beta$'s.
    \end{itemize}
\item
    \begin{itemize}
    \item The size of $A$ is smaller than the size of $B$.
    \item $A$ and $B$ are (or, \emph{have}) the same size.
    \item The size of $A$ is larger than the size of $B$.
    \end{itemize}
\endEnum
\par
Regarding this last group, we emphasize  that we are not arguing
that there \emph{really are} such things as set sizes, nor that there
\emph{are not really} such things.
Statements about \emph{sizes} can be translated in familiar and
long-winded ways into statements about sets, though we will not bother
to do so.
\par
I have identified the standard theory, Cantor's, with two principles
about set size.
The term \emph{size}, however, is rarely used in connection with
Cantor's theory; so it might be wondered whether the standard theory is
really so standard.
In stating \OO and \CANTORLT, Cantor used the terms \notion{power}
and \notion{cardinal number} rather than \notion{size}.
In the literature, the term  \notion{cardinal number} (sometimes just
\notion{number}) is used most frequently.
If someone introduces \notion{cardinal number} as a defined predicate
or as part of a contextual definition (e.g. \qphrase{We say that two sets
have the \emph{same cardinal number} just in case \dots}), there
is no point in discussing whether that person is right about  \notion{size}.
\par
Though Cantor's theory is usually taken as a theory of set size, it
can also be taken as \emph{just} a theory of one-one correspondences.
More specifically, saying that two sets are \notion{similar} iff
they are in one-one correspondence can either be taken as a claim about
size or be regarded as a mere definition.
Whether or not \notion{similarity} is coextensive with \notion{being the same size}, the definition is worth making.
The relation picked out is well studied and well worth the study.
The technical brilliance of the theory attests to this:
it has given us the transfinite hierarchy, the continuum problem, and
much else.
In addition, the theory has consequences which do not \emph{prima facie}
seem to have anything to do with size or similarity: the existence
of transcendental numbers comes to mind.
All of this is to say that the interest in one-one correspondence
has not been sustained solely by its identification with the notion of
size.
Hence, denying that they are the same does not endanger the theory
of one-one correspondences, \emph{per se}.
\par
But most mathematicians and philosophers do not use \notion{cardinal number} as a mere abbreviation.
They use the term in just the way that we use \notion{size} and slide
freely among (A), (B), and (C).
This is true, in particular, of Cantor, who offered \OO as a theory;
indeed, he offered an argument for this theory.

\Section{1.2}{Cantor's argument}
Cantor bases his argument for \OO on the idea that the size of
a set, it's cardinal number, depends on neither the particular elements
it contains nor on how those elements are arranged:
\begin{quote}
The cardinal number, \Card{M}, of a set $M$ (is the general concept
which) arises from $M$ when we make abstraction of the nature of its
elements and of the order in which they are given.
\cite[p.86]{Cantor}

\end{quote}
But to say that the cardinal number of a set does not depend on certain
things is not to say what the cardinal is.
Neither does it insure that two sets have the same cardinal number
just in case they are in one-one correspondence.
To flesh out this notion of \notion{double abstraction}, Cantor reduces it
to a second abstraction operator, one which works on the elements
of sets rather than the sets themselves:
\begin{quote}
Every element, $m$, if we abstract from its nature becomes a \notion{unit}.
\cite[p.86]{Cantor}

\end{quote}
and so, concludes Cantor:
\begin{quote}
The cardinal number, \Card{M}, is a set composed of units (which has
existence in our minds as an intellectual image or projection of $M$.)
\cite[p.86]{Cantor}
\end{quote}
According to Jourdain, Cantor
\begin{quote}
distinguised very sharply between an aggregate and a cardinal number
that belongs to it: \qphrase{Is not an aggregate an object \emph{outside}
us, whereas its cardinal number is an abstract picture of it
\emph{in our} mind.}
\cite[p.80]{Cantor}
\end{quote}
I have parenthesized the expressions above where Cantor describes
cardinal numbers as mental entities.
Nevertheless, I can only make sense of his arguments insofar as
he treats cardinal numbers as sets:
he refers to them as \wq{definite aggregates}, supposes that
they have elements, and employs mappings between cardinal numbers
and other sets.
\par
The following three statements seem to express Cantor's intent:

\begin{gather}
\Card{M} = \theSetst{y}
    {\exists x \in M \Aand y = \tABST(x)} \\
\forall x \exists y (y = \tABST(x) \Aand \tUNIT(y))\\
\forall M \forall y (y \in \Card{M} \then \tUNIT(y))
\end{gather}

\emph{\tABST(x)} is to be read as \readAs{the result of abstracting from the
element $x$}, \emph{\tUNIT(x)} as \readAs{$x$ is a unit}, and $\Card{M}$
as \readAs{the cardinal number of $M$}.
\par
So (1.1) gives a definition of cardinal number, in terms of the
operation of abstraction, from which Cantor proves both \OOA
and \OOB.
\begin{align*}
M \eqsize N &\then \Card{M} = \Card{N} \tag{\OOA}\\
\Card{M} = \Card{N} &\then M \eqsize N   \tag{\OOB}
\end{align*}

\OOA is true, says Cantor, because
\begin{quote}
the cardinal number \Card{M} remains unaltered if in the place
of one or many or even all elements $m$ of $M$ other things are
substituted. \cite[p.80]{Cantor}
\end{quote}
and so, if $f$ is a one-one mapping from $M$ onto $N$, then in
replacing each element, $m$, of $M$ with $f(m)$
\begin{quote}
$M$ transforms into $N$ without change of cardinal number. (p.88)
\end{quote}
In its weakest form, the principle Cantor cites says that if a
single element of $M$ is replaced by an arbitrary element not in $M$
then the cardinal number of the set will remain the same.
That is,
    \begin{multline}
    (a \in M \Aand b \notin M)\\
    \Aand
    N = \theSetst{x}{(x \in M \Aand x \notin a) \Or x \in b}\\
    \then \Card{M} = \Card{N}
    \end{multline}
The reasoning is clear: so far as the cardinal number of a set
is concerned, one element is much the same as another.
It is not the elements of a set, but only their abstractions, that
enter into the cardinal number of a set.
But abstractions of elements are just units; so one is much the same as
another.
\par
\OOB is true, Cantor says, because
\begin{quote}
\dots \Card{M} grows, so to speak, out of $M$ in such a way that
from every element $m$ of $M$ a special unit of \Card{M} arises.
Thus we can say that $M \eqsize \Card{M}$.
\end{quote}
So, since a set is similar to its cardinal number, and similarity
is an equivalence relation, two sets with same cardinal number are
similar.
Unless each element of a set abstracts to a \wq{special}, i.e. distinct,
unit, the correspondence from $M$ to its cardinal number will be
many-one and not one-one.
A weak version of this principle is:

\begin{equation}
\begin{split}
 M = \theSet{a,b} &\Aand a \neq b\\
 &\then \Card{M} = \theSet{\tABST(a), \tABST(b)}   \\
    & \quad \Aand \tABST(a) \neq \tABST(b)
\end{split}
\end{equation}
These two arguments do one another in.
(1.4) says that replacing an element of a set with any element not
in the set does not affect the cardinality.
But, by the definition of \Card{M}, (1.1), this means that
\begin{equation}
\forall x \forall y (\tABST(x) = \tABST(y))
\end{equation}
For, consider an arbitrary pair of elements, $a$ and $b$.
Let $M = \theSet{a}$ and let $N = \theSet{b}$.
So, the conditions of (1.4) are met and $\Card{M} = \Card{N}$.
But $\Card{M} = \theSet{\tABST(a)}$ and $\Card{N} = \theSet{\tABST(b)}$, by (1.1).
So $\tABST(a) = \tABST(b)$.
Generalizing this argument yields (1.6).
\par
So Cantor's argument for \OOA only works by assigning all
non-empty sets the same, one-membered, cardinal number.
But, this contradicts \OOB.
\par
Conversely, the argument that a set is similar to its cardinal
number relies on (1.5), which entails
\begin{equation}
\forall x \forall y (x \neq y \then \tABST(x) \neq \tABST(y) )
\end{equation}
assuming only that any two objects can constitute a set.
But if the abstractions of any two elements are distinct, then no two
sets have the same cardinal number as defined by (1.1), \emph{contra}
\OOA.
\par
There is no way to repair Cantor's argument.
Rather than leading to a justification of \OO, Cantor's definition
of cardinal number is sufficient to refute the principle.
The negation of (1.6) is:
\begin{equation}
\exists x \exists y ( \tABST(x) \neq \tABST(y))
\end{equation}
So one of (1.6) and (1.8) must be true.
We have shown that (1.6) contradicts \OOB.
Similarly, (1.8) contradicts \OOA; if $a$ and $b$ have
distinct abstractions, then \theSet{a} and \theSet{b} have
distinct cardinal numbers, \theSet{\tABST(a)} and \theSet{\tABST(b)},
despite the fact that they are in one-one correspondence.
So \OO is false whether (1.6) or its negation, (1.8), is true.
\Section{1.3}{Cantor and the logicists}
Though both Frege and Russell accepted Cantor's theory of cardinality,
neither accepted Cantor's argument.
Frege spends an entire chapter of the \emph{Grundlagen} mocking
mathematicians from Euclid to Schroder for defining numbers as sets of
\scare{units}.
He neatly  summarizes the difficulty with such views:
\begin{quote}
If we try to produce the number by putting together different distinct
objects, the result is an agglomeration in which the objects remain still
in possession of precisely those properties which serve to distinguish
them from one another, and that is not the number.
But if we try to do it in the other way, by putting together identicals,
the result runs perpetually together into one and we never reach a
plurality \dots
\par
The word \wq{unit} is admirably adapted to conceal the difficulty \dots
We start by calling the things to be numbered \wq{units} without detracting
from their diversity; then subsequently the concept of putting
together (or collecting, or uniting, or annexing, or whatever we
choose to call it) transforms itself into arithmetical addition, while
the concept word \wq{unit} changes unperceived into the proper name \wq{one}.
 \cite[pp 50-51]{Frege}
\end{quote}
These misgivings about units do not prevent Frege from basing his
definition of \wq{number} and his entire reduction of arithmetic on
Cantor's notion of one-one correspondence.
\qphrase{This opinion}, says Frege, \qphrase{that numerical equality or identiy must
be defined in terms of one-one correlation, seems in recent years to
have gained widespread acceptance among mathematicians}
 \cite[pp. 73-74]{Frege}.
Frege cites Schroder, Kossak, and Cantor.
\par
Russell displays similar caution about Cantor's argument
 \cite[p.305]{Russell} and similar enthusiasm
for his theory (see the quote at the beginning of Chapter 2, for
example.)
\par
Of course, Frege and Russell \scare{cleaned up} Cantor's presentation of the
theory.
Russell, for example, notes that Cantor's statement (1) is not
a \wq{true definition} and
\begin{quote}
merely presupposes that every collection has some such property as
that indicated --- a property, that is to say, independent of its
terms and their order;
depending, we might feel tempted to add, only upon their number.
\cite[pp. 304-305]{Russell}
\end{quote}
So Russell, and similarly Frege, relied upon the principle of abstraction
to obtain a \wq{formal definition} of cardinal numbers, in contrast to
Cantor, who had \qphrase{taken} number \qphrase{to be a primitive idea} and had to
rely on \qphrase{the primitive proposition that every collection has
a number.}
\cite[p.305]{Russell}
\par
So, while some people regard Cantor's \OO as \emph{just} a definition and
others embrace it as a theory, the logicists have it both ways:
adding \OO as a formal definition to set theory (or, as they would
call it, logic) they have no obligation to defend it and can steer
clear of peculiar arguments about \notion{units};
at the same time, they can advance it as a great lesson for simple common
sense.
\par
The logicists' adoption of Cantor's theory of cardinality needs no
great explanation:
it \scare{came with} set theory and, to a large extent, motivated set theory
and determined its research problems.
But there are two specific reasons that they should have seized upon
\OO and \CANTORLT.
First, they both have the form of definitions, no matter how they are
intended.
So the notion of cardinality is \scare{born reduced}.
\par
Second, Cantor's theory clears the way for other reductions.
Suppose, for example, you wish to reduce ordered pairs to sets.
Well, you have to identify each ordered pair with a set and define
the relevant properties and relations among ordered pairs in
terms of properties and relations among sets.
One of the relations that has to be maintained is \notion{identity};
so each ordered pair must be identified with a distinct set.
In addition, the relative sizes of sets of ordered pairs should be
preserved under translation.
But, if \OO is the correct theory of size, then this second condition
follows from the first, since the existence of one-one correspondences
will be preserved under a one-one mapping.

\Section{1.4}{Aims and outline}
It would be naive to suppose that people's faith in Cantor's theory
would be shaken either by refuting specific arguments for \OO
or by associating the acceptance of that theory with a discredited
philosophy of mathematics.
Such points may be interesting, but in the absence of an alternative
theory of size they are less than convincing.
\par
This dissertation presents such an alternative.
Chapter 2 canvasses common sense intuitions for some basic principles
about set size.
Chapter 3 reorganizes those principles into a tidy set of axioms, offers
an account of where the intuitions come from (viz. known facts about
finite sets), and mines this source for additional principles.
Chapters 4 and 5 prove that the theory so obtained is \emph{complete},
in the sense that it embraces all facts about finite sets of
a certain kind (i.e. expressible in a particular language).
Finally, Chapter 6 elaborates additional principles that concern only
sets of natural numbers and demonstrates that these additional
principles, together with the theory in Chapter 3, are satisfiable
in the domain of sets of natural numbers.

\Chapter{2}{The General Theory}

\begin{quote}
The possibility that whole and part may have the same
number of terms is, it must be confessed, shocking to common sense
\dots Common sense, therefore, is here in a very sorry plight;
it must choose between the paradox of Zeno and the paradox of Cantor.
I do not propose to help it, since I consider that, in the face of
the proofs, it ought to commit suicide in despair.
\cite[p. 358]{Russell}
\end{quote}
Is common sense confused about set size, as Russell says, or is there
a way of elaborating on common sense to get a plausible and reasonably
adequate theory of cardinality?
To be \squote{plausible}, a theory should at least avoid principles and
consequences which violate common sense.
To be \squote{reasonably adequate}, a theory has to go beyond bare intuitions:
it should not rest with trivialities and it should answer as many
questions about set size as possible, though it need not be complete.
Plausibility and adequacy are conflicting demands:
the first says that there should not be too many principles (no false
ones, consistency), the second that there should not be too few principles.
\par
In the Introduction, I argued that a coherent theory of cardinality
has to contain some principles that refer to the kinds of objects
in sets, \emph{pace} Cantor.
In this chapter, however, I want to see how far we can go without
such principles; i.e. how much can you say about \notion{smaller than} without
using predicates (other than \notion{identity}) which relate the
members of the sets being compared?
I shall begin by stating a number of principles and explaining why
they are included in the general theory.

\Section{2.1}{The Theory CORE}
First, there is the \SUBSET principle:

\namedSentence{\SUBSET}
If $x$ is a proper subset of $y$, then $x$ is smaller than $y$.
\endNamedSentence
The reason for including \SUBSET should be obvious.
What has prompted the search for an alternative to the standard
theory of cardinality is the conflict between \OO and \SUBSET.
Now, it is often said that common sense supports both of these
principles and that it is in doing so that common sense is confused.
From this it is supposed to follow that common sense cannot
be relied upon, so we should opt for \OO,with the technical attractions
that it provides.
\par
But, there's a difference in the way that common sense supports these
two principles.
There is no doubt that you can lead an unsuspecting person to agree
to \OO by focusing their attention on forks and knives, husbands and
wives, and so forth: i.e. finite sets.
With carefully chosen examples, say the odds and the evens,
you might even convince someone that
\OO is true for infinite sets, too.
Now, I do not think that such guile is needed to lead someone to
agree to \SUBSET, but that's not what my argument depends on.
The argument hinges on a suggestion about how to resolve cases
where mathematical intuitions seem to conflict.
The suggestion is to see what happens with particular cases on which
the principles conflict \emph{before} you've lead someone to agree
to either of the general statements.\par
So, if you want to find out what common sense really thinks about
\SUBSET and \OO, you would present people with pairs of infinite sets,
where one was a proper subset of the other.
I've actually tried this, in an unscientific way, and what I've
gotten, by and large, is what I expected: support for \SUBSET.
(\squote{By and large} because many people think that all infinite sets
have the same size: \emph{Infinity}.)
\par
Naturally, I would not venture that this sort of technique, asking
people, is any way to find out which of \SUBSET and \OO is \emph{true}.
People's intuitions about mathematics are notoriously unreliable, not
to mention inconstant.
Of course, harping on this fact might engender some unwarranted
skepticism about mathematics.
What I am suggesting is that there might be a rational way of
studying mathematical intuitions and that we should at least explore
this possibility before proclaiming common sense to be hopelessly
confused on mathematical matters.
\par
So, \SUBSET, all by itself seems to be a plausible alternative to
Cantor's theory. though it surely is not enough.
Given just this principle, it's possible that one set is smaller
than another just in case it is a proper subset of the other.
It's clear that we need additional principles.
All of the following seem worthy (where $x < y$ is to be read as
\readAs{$x$ is smaller than $y$},
$x \eqsize y$ is to be read as
\readAs{$x$ is the same size as $y$},
$x > y$  as
\readAs{$x$ is larger than $y$}.
The Appendix gives a full account of the notations used throughout
this thesis.)
\aTheory{2.1.1}{QUASI-LOGICAL}
\begin{gather*}
x < y \then \neg y < x\tag{\ASYMLT}\\
x > y \then \neg y > x\tag{\ASYMGT}\\
(x < y \Aand y < z) \then x < z\tag{\TRANSLT}\\
(x > y \Aand y > z) \then x > z\tag{\TRANSGT}\\
 x \eqsize y \then \tINDISC(x,y) \tag{\INDISCEQ}\\
x \eqsize x \tag{\REFEQ}\\
x \eqsize y \then y \eqsize x\tag{\SYMEQ}\\
(x \eqsize y \Aand y \eqsize z) \then x \eqsize z\tag{\TRANSEQ}\\
x >y \bicond y < x\tag{\DEFGT}
\end{gather*}
where $\tINDISC(x,y)$ abbreviates
    $\forall z ( (z < x \bicond z < y) \Aand (z > x \bicond z > y))$
\endTheory
\par
We  call the principles listed above \introSTerm{quasi-logical principles}
because it is tempting to defend them as logical truths.
Consider the first principle, for example, in unregimented English:
\namedSentence{\ASYMLT}
If $x$ is smaller than $y$, then $y$ is not smaller than $x$.
\endNamedSentence
This sentence can be regarded as an instance of the schema:
\namedSentence{\ASYMF}
If $x$ is $F$-er than $y$, then $y$ is not $F$-er than $x$.
\endNamedSentence
where $F$ is to be replaced by an adjective from which comparatives
can be formed, e.g. \qword{tall}, \qword{short}, \qword{happy}, but
not \qword{unique} or \qword{brick}.
It appears that every instance of this schema is true, so it could
be maintained that each is true in virtue of its form, that each is
a logical truth.
\par
The other principles might be defended in the same way, though
the schema for \INDISCEQ would have to be restricted to
triples of corresponding comparatives, for example:
\emph{is smaller than},
\emph{is larger than},
\emph{is the same size as}.
\par
But using such observations to support these principles would be
problematic for two reasons.
First, it would require taking positions on many questions about
logical form and grammatical form which would take us far afield
and, possibly, antagonize first-order logicians.
Second, there are some instances of the schemata that make for
embarassing counterexamples: \qword{further east than} (in a round
world) and \qword{earlier than} (in, I'm told, a possible world).
\par
So, it might be that casting the principles above as instances
of the appropriate schema would only explain why they are part
of common sense.
What remains clear is that a theory of cardinality which openly
denied any of these principles would be implausible:
it would be ridiculed by common sense and mathematical sophisticates
alike.
I can just barely imagine presenting a theory which, for fear of
inconsistency, withheld judgment on one or more of these
principles.
But to do so without good reason would be counterproductive.
It seems that if a case could be made that these statements, taken
together, are inconsistent with \SUBSET, that would be good reason
to say that there is no reasonably adequate alternative to Cantor's
theory.
Since my goal is to counter such a conclusion, it seems that the
proper strategy is to include such seemingly obvious truths and
to show that the resulting theory is consistent.
So the strategy here is not to adduce principles and argue for the
truth of each. This would be impossible, given that the principles
are logically contingent.
Instead, our approach is to canonize what common sense holds to be
true about cardinality and show that the result is consistent and
reasonably adequate.
\par
As long as we restrict ourselves to \SUBSET and the quasi-logical
principles, consistency is no problem.
After all, what do the quasi-logical principles say?
Only that \notion{smaller than} is a partial ordering, the \notion{larger than} is
the converse partial ordering, that \notion{the same size as}
is an equivalence relation, and that sets of the same size are
indiscernible under the partial orderings.
So, if we are given any domain of sets, finite or infinite, we get
a model for our theory by assigning to $<$ the relation of
\notion{being a proper subset of}, assigning to $>$ the relation
of \notion{properly including}, and assigning to $\eqsize$
the identity relation.
Since common sense knows that different sets can be the same size, there
must be some additional principles to be extracted from common sense.
\par
We shall now consider some principles which cannot be regarded as
quasi-logical.
\par
First, there is the principle of \introSTerm{trichotomy}.
\namedEq{TRICH}{x < y \Or x \eqsize y \Or x > y}
which says that any two sets are comparable in size.
While a theory of set size which excluded \TRICH might escape
ridicule, it would surely be regarded with suspicion.
Indeed, if the principles of common sense were incompatible with \TRICH, this would undoubtedly be used to discredit them.
\par
Second, there is the \introSTerm{representation principle}.
\namedEq{\REPLT}{
    x < y \then \exists \xP ( \xP \eqsize x \Aand \xP \subset y)
}
which says that if a set, $x$, is smaller than another, $y$,
then $x$ is the same size as some proper subset of $y$.
Now, this is a principle which common sense has no particular
feelings about.
Analogous statements about physical objects are neither intuitive
nor very clearly true.
For example, (1) does not stand a chance of being regarded as true:
\numbSentence{1}
If one table is smaller than another, then the first is the same
size as some proper part of the second.
\endNumbSentence
if \qword{part} is taken to mean \qword{leg or top or rim or ...}.
Even if common sense can be persuaded to take particles and arbitrary
fusions of such as parts of tables, no one should condemn its residual
caution about (1).
If \REPLT is true, then it seems to be an interesting and special
fact about sets.
\par
\REPLT was originally included in this theory for technical
reasons;
it makes it easier to reduce the set of axioms already presented
and it provides a basis for several principles not yet presented.
\REPLT may be open to doubt, but it is not a principle that
Cantorians could complain about, for it is entailed by Cantor's
definition of $<$:
\namedEq{\CANTORLT}
{x < y \bicond \neg(x \eqsize y) \Aand
    \exists \xP ( \xP \eqsize x \Aand \xP \subset y)
}
If \CANTORLT is regarded as a principle instead of
a definition, then it is entailed by the principles we have
already mentioned:
\begin{quote}
If $x < y$, then $\neg(x \eqsize y)$, by \INDISCEQ and \ASYMLT.
By \REPLT, some proper subset of \y, say \xP, must be the
same size as \x.
But $\xP < y$, by \SUBSET:
so $x < \yP$, by \INDISCEQ.
\par
Conversely, if $\xP \eqsize x$ and $\xP \subset y$, then $\xP < y$, by \SUBSET.
So $x < y$, by \INDISCEQ.
\end{quote}
There are more principles to come, but before proceeding, I'd
like to take stock of what we already have.
First, I want to reduce the principles mentioned above to
a tidy set of axioms.
Second, I want to estimate how far we've gone.
\par
The entire set of principles already adopted are equivalent
to the following, which will be referred to as \emph{the core theory}.
\aTheory{2.1.2}\introTheory{CORE}{core theory}
\begin{align*}
x \subset y &\then x < y \tag{\SUBSET}\\
x \eqsize y &\bicond \tINDISC(x,y) \tag{\DEFEQ}\\
(x \eqsize y \Aand y \eqsize z) &\then x \eqsize z  \tag{\TRANSEQ}\\
 x >y &\bicond y < x \tag{\DEFGT}\\
  \neg( x &< x ) \tag{\IRREFLT} \\
  x < y \Or x &\eqsize y \Or x > y \tag{\TRICH}
\end{align*}
\endTheory
\par
The only axiom in \CORE that has not already been introduced is \DEFEQ,
which is logically equivalent to the conjunction of \INDISCEQ and
its converse \EQINDISC:
\begin{align*}
    x \eqsize y &\then \tINDISC(x,y) \tag{\INDISCEQ}\\
    \tINDISC(x,y) &\then x \eqsize y \tag{\EQINDISC}
\end{align*}
\EQINDISC says that if two sets fail to be the same size, then their
being different in size is attributable to the existence of some
set which is either smaller than one but not smaller than the other,
or larger than one but not larger than the other.
\Theorem{2.1.3}
Let $T = \QUASI \adjoin \SUBSET \adjoin \REPLT$.
Then $T \equiv \CORE$.
\endTheorem
\Proof
\BothWays
\ProvesWay  $T \Proves  \CORE$.
        We only need to show that \EQINDISC is entailed by $T$.
        Suppose $\neg(x \eqsize y)$.
        So $x < y$ or $y < x$, by \TRICH.
        But $\neg(x < x)$ and $\neg(y < y)$, by \IRREF.
\IsProvedWay   $\CORE \Proves T$
    \Enum
    \item \TRANSLT. If $y < z$, there is a \yP such that
        $y' \eqsize y$ and $\yP \subset  z$ by \REPLT.
        If $x < y$, then $x < \yP$ because $\yP \eqsize y$,
        by \DEFEQ.
        So, there is an $\xP$ such that $\xP \eqsize x$ and
        $\xP \subset \yP$.
        So $\xP \subset z$ and, by \SUBSET, $\xP < z$.
        But then $x < z$ by \DEFEQ.
    \item \ASYMLT. If $x < y$ and $y < x$, then $x < x$, by
        \TRANSLT, {\it contra} \IRREFLT.
    \item \TRANSGT, \ASYMGT, and \IRREFGT follow from the corresponding
        principles for $<$ and \DEFGT.
    \item \INDISCEQ, \SYMEQ, \TRANSEQ, and \REFEQ are logical
        consequences of \DEFEQ.
    \endEnum
\enBoth
\endProof
\CORE is consistent.
In fact, two kinds of models satisfy \CORE.
\Def{2.1.4}
\Enum
\item \Model{A} is a \introSTerm{finite class model} with basis $x$
    \iff
    \Enum
    \item $\Domain{A} = \PowerSet{x}$, where $x$ is finite.
   \item \Model{A} \Satisfies $a \subset b$ \iff $a \subset b$
    \item \Model{A} \Satisfies $a < b$ \iff $\Card{a} < \Card{b}$
    \endEnum
\item \Model{A} is a \introSTerm{finite set model} \iff
    \Enum
    \item $\Domain{A} = \theSetst{x}{x \emath{\ is a finite subset of\ } Y}$,
        for some infinite $Y$.
    \item \Model{A} \Satisfies $a \subset b$ \iff $a \subset b$
    \item \Model{A} \Satisfies $a < b$ \iff $\Card{a} < \Card{b}$
    \endEnum
\endEnum
\endDef
Models can be specified by stipulation the \notion{smaller than}
relation since \notion{larger than} and \notion{same size as} are defined in
terms of \notion{smaller than}.
\Fact{2.1.5}
\Enum
\item If \Model{A} is a finite class model,
        then \Model{A} \Satisfies \CORE.
\item If \Model{A} is a finite set model,
        then \Model{A} \Satisfies \CORE.
\endEnum
\endFact
\Proof
In both cases, the finite cardinalities determine a quasi-linear
ordering of the sets in which any set is higher than any of its
proper subsets.
\endProof
The normal ordering of finite sets is, in fact, the only one that satisfies
\CORE.
By adding \TRICH we have ruled out all non-standard interpretations
of $<$.
\Theorem{2.1.6}
Suppose that \Model{A} is a model such that
\Enum
\item If $x \in \Domain{A}$, then $x$ is finite,
\item If $x \in \Domain{A}$ and $y \subset x$, then $y \in \Domain{A}$, and
\item $\Model{A} \Satisfies \CORE$.
\endEnum
Then \Model{A} \Satisfies  $a < b$ \iff $\Card{a} < \Card{b}$
\endTheorem
\Proof
We shall prove (*) by induction on \n:
\begin{gather*}
\tag{*}
    \emath{If} \Card{a} = n \emath{, then}
    \mA \Satisfies(a\eqsize b)
    \tiff \Card{b} = n
\end{gather*}
\begin{flalign*}
    \Suppose    n = 0\\
    \Then       a = \NullSet\\
    \Butif      \Card{b} = 0\\
    \Then       \mA \Satisfies (a \eqsize b)
                        \by{\REFEQ}\\
    \Andif      \Card{b} \neq 0\\
    \Then       a \subset b\\
    \so         \mA \Satisfies (a < b) \by{\SUBSET}\\
    \so         \mA \nSatisfies (a \eqsize b)
                        \by{\INDISCEQ}
\end{flalign*}
    Now, suppose that (*) is true for all $i \leq n$:
\begin{flalign*}
    \If         \Card{a} = \Card{b} = n+1\\
    \Then       \mA \nSatisfies (a \eqsize b)
                   \emath{, as follows:}\\
    \Suppose    \mA \Satisfies (a < b) \\
    \So       \mA \Satisfies (\aP \eqsize a)
                    \ForSome{\aP \subset b}
                    \by{\REPLT}\\
    \Butif      \mA \Satisfies \aP \subset b \\
    \Then       \Card{\aP} \leq n \\
    \so         \Card{a} \leq n
                    \by{(*), contra our hypothesis.}\\
    \Butif      \Card{a} = n+1\\
    \AND        \mA \Satisfies (a \eqsize b)\\
    \Then       \Card{b} \geq n+1 \by{induction.}\\
    \Butif      \Card{b} > n+1\\
    \pick       \bP \subset b, \bP \in A,\\
    \with       \Card{\bP} = n+1 \by{condition (b)}\\
    \So         \mA \Satisfies (\bP < b) \by{\SUBSET}\\
    \AND        \mA \Satisfies (\bP \eqsize a)\\
    \So         \mA \Satisfies (a < b)
                    &&\emath{contra our supposition.}\\
    \so         \Card{b} = n+1
\end{flalign*}
\endProof
\Section{2.2}{Addition of Set Sizes}
We shall now extend \CORE to an account of
addition of set sizes.
Since the domains of our intended models contain only
sets and not sizes of sets, we have to formulate our
principles in terms of a three-place predicate true of
triples of sets:
\tSUM(\x, \y, \z) is to be read as
\readAs{the size of \z is the sum of the sizes of \x and \y}.
\par
The following principles are sufficient for a theory of
addition:
\aTheory{2.2.1}{Addition}
\begin{align*}
\intertext{Functionality of addition}
\tSUM(x,y,z) &\then (x \eqsize \xP \bicond \tSUM(\xP,y,z)\tag{a}\\
\tSUM(x,y,z) &\then (y \eqsize \yP \bicond \tSUM(x,\yP,z) \tag{b}\\
\tSUM(x,y,z) &\then (z \eqsize \zP \bicond \tSUM(x,y,\zP)\tag{c}\\
\intertext{Addition for disjoint sets}
    x \intersect y = \emptyset &\then
    \tSUM(x,y,x \union y) \tag{\DISJPL}\\
\intertext{Monotonicity of Addition}
    \tSUM(x,y,z) &\then x < z \Or x \eqsize z \tag{\MONOT}
\end{align*}
\endTheory
\FUNCPL says that sets bear \tSUM relations to one
another by virtue of their sizes alone.
This condition must clearly be met if \tSUM is to be read as
specified above.
\par
\DISJPL tries to say what function on sizes the \tSUM
relation captures by fixing the function on paradigm cases:
disjoint sets.
But \FUNCPL and \DISJPL leave open the possibility that
addition is cyclic:
suppose we begin with a finite class model whose basis has
\n elements and assign to \tSUM those triples \tuple{x,y,z} where
\[
    \Card{z} = (\Card{x} + \Card{y}) \bmod (n+1)
\]
Both \FUNCPL and \DISJPL will be satisfied, though the
interpretation of \tSUM does not agree with the intended
reading.
\MONOT rules out such interpretations.
\par
Given an interpretation of '<' over a power set there
is at most one way of interpreting \tSUM which satisfies
\ADDITION.
We shall show this by proving that \ADDITION and \CORE
entail \DEFPL:
\begin{multline}
    \tSUM(x,y,z) \bicond        \tag{\DEFPL}\\
    \exists \xP \exists \yP
    (x \eqsize \xP \Aand y \eqsize \yP \Aand
    \xP \intersect \yP = \emptyset \Aand
    \xP \union \yP = z)
\end{multline}
\DEFPL says that the extension of \tSUM is determined by the extension of <.
\par
A model of \CORE must satisfy an additional principle,
\DISJU, if \tSUM is to be interpreted in a way compatible
with \ADDITION,
\begin{multline}
( x \eqsize \xP \Aand y \eqsize \yP \Aand
x \intersect y = \emptyset \Aand
\xP \intersect \yP = \emptyset)    \tag{\DISJU}   \\
\then
(x \union y) \eqsize (\xP \union \yP)
\end{multline}
(Note: The proofs in this chapter will use boolean
principles freely, despite the fact that we have not
yet introduced them.)
\Theorem{2.2.2}
\CORE + \ADDITION \Proves \DISJU
\endTheorem
\Proof
\begin{flalign*}
    \Suppose    (x \eqsize \xP \Aand y \eqsize \yP
                \Aand x \intersect y = \emptyset
                \Aand  \xP \intersect \yP = \emptyset)\\
    \Then       \tSUM(x,y,x\union y)\\
    \AND        \tSUM(\xP,\yP,\xP\union\yP)
                    \by{\DISJPL}\\
    \so         \tSUM(\xP,y,x \union y) \by{\FUNCPL(a)}\\
    \so         \tSUM(\xP,\yP,x \union y) \by{\FUNCPL(b)}\\
    \so         x \union y \eqsize \xP\union\yP
                        \by{\FUNCPL(c)}
\end{flalign*}
\endProof
If the minimal
conditions on addition are to be satisfied in a model
of \CORE, then \tSUM has to be definable by \DEFPL.
\Theorem{2.2.3}
\CORE + \ADDITION \Proves \DEFPL
\endTheorem
\Proof

\IfWay
\begin{flalign*}
    \Suppose    \tSUM(x,y,z)\\
    \so         x < z \Or x \eqsize z \by{\MONOT}\\
    \Butif      x \eqsize z, \\
    \Let        \xP = z\\
    \AND        \yP = \NullSet;\\
    \Then       \tSUM(\xP,\yP,z) \by{\DISJPL}\\
    \so         \tSUM(x,\yP,z) \by{\FUNCPL(a)}\\
    \so         \yP \eqsize y \by{\FUNCPL(b)}\\
    \Andif      x < z\\
    \pick       \xP \subset z     \\
    \with       \xP \eqsize x \by{\REPLT.}\\
    \Let        \yP = z - \xP\\
    \but        \yP \eqsize y\emath{, as before.}
 \end{flalign*}
\OnlyWay
 \begin{flalign*}
    \Suppose    (x \eqsize \xP \Aand y \eqsize \yP \Aand
    \xP \intersect \yP = \emptyset \Aand
    \xP \union \yP = z)\\
    \Then       \tSUM(\xP,\yP,z) \by{\DISJPL}\\
    \so         \tSUM(x,\yP,z) \by{\FUNCPL(a)}\\
    \so         \tSUM(x,y,z) \by{\FUNCPL(b)}
 \end{flalign*}
\endProof
\aTheory{2.2.4} \introTheory{\EXCORE}{extended core}
is $\CORE \union \ADDITION$.
\endTheory
\Fact{2.2.5}
The following are consequences of \EXCORE:
\begin{align*}
    x \intersect y_1 = x \intersect y_2 = \emptyset
    \Aand y_1 \eqsize y_2
        &\then x \union y_1 \eqsize x \union y_2\\
    x \intersect y_1 = x \intersect y_2 = \emptyset
    \Aand y_1 < y_2
        &\then x \union y_1 < x \union y_2\\
     x_1 \intersect y_1 = x_1 \intersect y_2 = \emptyset
    \Aand x_1 \eqsize x_2
    \Aand y_1 < y_2
        &\then x_1 \union y_1 < x_2 \union y_2\\
    x \subset z \Aand y \subset z \Aand x < y
    &\then
    (z - y) < (z - x)\\
    x \subset z \Aand y \subset z \Aand x \eqsize y
    &\then
    (z - y) \eqsize (z - x)\\
    x < y &\then
    \exists \yP (\yP \eqsize y \Aand x \subset \yP)\\
    x \eqsize y &\then
    x - (x \intersect y) \eqsize
    y - (x \intersect y)\\
    x < y &\then
    x - (x \intersect y) <
    y - (x \intersect y)
\end{align*}
\endFact
\Proof
The proofs are elementary.
\endProof

\Chapter{3}{A Formal Theory of Class Size}

In the preceding chapter, we searched for principles that accord with pre-Cantorian ideas about sizes of sets.
We produced several such principles, constituting \EXCORE,
and found two kinds of interpretations which satisfied these principles.
In one case, the domains of the interpretations were finite power
sets.
In the other case, the domains consisted of all finite subsets of a given infinite set.
Insofar as both kinds of interpretation have quantifiers ranging over finite sets, we may say that they
demonstrate that \EXCORE is true when it is construed as being about finite sets.
\par
Our goal is to show that this general theory of
set size can be maintained for infinite sets as well as for finite sets.
We shall show this by constructing a model for the general
theory whose domain is the power set of the natural numbers.
But the ability to construct such a model is interesting
only to the extent that the general theory it satisfies is
reasonably adequate.
Suppose, for example, that we offered as a general theory
of size the axioms of \CORE other than \REPLT.
Call this theory \textit{T}.
So $T$ just says that \notion{smaller than} is a quasi-linear ordering which extends the partial ordering given by the proper subset relation.
Since any partial ordering can be extended to a quasi-linear ordering, $T$ has a model over \PowerN.
But unless we have a guarantee that the model constructed will satisfy, say, \DISJU,
the existence of the model does not rule out the possibility that $T$ is incompatible with
\DISJU.
For a particular principle, $\phi$, in this case \DISJU, we
may take one of three tacks:
(1) add $\phi$ to T, obtaining \Tp, and show that \Tp has a model over \PowerN;
(2) show that $\phi$ \emph{is} inconsistent with \T and argue that \T is somehow more fundamental or more intuitive than \gphi;
(3) acquiesce in ignorance of whether \T and \gphi are compatible and argue that \emph{if} they are incompatible, then T
should be maintained anyway.
\par
Below, we deal with \DISJU as in (1), since \DISJU is in \EXCORE.
Cantor's principle \OO is dealt with as in case (2).
It seems futile to try to rule out the need to resort to
the third approach for any cases at all, but we can reduce
this need to the extent that we include in our general theory, T, as many plausible statements as possible.
\par
We cannot construct \T by taking \emph{all} statements which are true for
finite sets.
Not only is \OO such a statement, but using the notion of \emph{all statements true for finite sets} presupposes that we have some idea of
the range of \squote{all statements}.
To avoid the problems involved in speaking of \squote{all statements},
we might instead settle for \squote{all statements in $L$}, where $L$ is some judiciously chosen language.
To avoid \OO, $T$ must fall short of the full expressive
power of the language of set theory.
\par
Consider, now, the axioms in \EXCORE.
Other than size relations, these axioms involve only
boolean operations and inclusion relations among sets.
They do not use the notions \notion{ordered pair},
\notion{relation}, or \notion{function}.
In short, the only set theory implicit in these axioms is boolean algebra, or a sort of \scare{Venn diagram set theory}.
This is not to say the axioms do not apply to relations,
functions, or other sets of ordered pairs, but only that they do not refer to these sorts of objects as such.
\par
In the next section, we define a language just strong enough to express \EXCORE.
We then construct a theory by taking all statements \emph{in that language} which are true over all finite power sets.
By drawing statements only from this relatively weak language, we arrive at a theory which can be satisfied over infinite power sets.
But, since we include in the language \emph{all} statements of the language which are true over any finite
power set, we know that no statement in that language can arise as something which ought to be true over infinite power sets but might be incompatible with our theory.
\par
There remains the possibility that we could follow the same strategy with a more expressive language, though it would have to remain less expressive that the full language of set theory.
In fact, such a language can be obtained by including a notion of the product of set sizes.
This in turn opens the possibility of a succession of richer languages and a corresponding succession of stronger theories of size.
At this point, the possible existence of any such hierarchy is sheer speculation; we mention it only to
emphasize that no claim is made here that we have the
strongest possible general theory of set size.
\Section{3.1}{CS - The Theory of Class Size}
The theories discussed in this paper will be formulated within first order predicate logic with identity.
To specify the language in which a theory is expressed, then, we need only list the individual constants, predicates, and operation symbols of the language and stipulate the \emph{rank}, or number of argument places, for each predicate and each operation symbol.

\Def{3.1.1}
\Enum
\item \LC, \introLanguage{\LC}{classes} is the first order language
 with individual constants $\NullSymbol$ and
\Identity,
the one-place predicate, \tATOM, the two-place predicate
\Symbol{\subset}, and the two place operation symbols, \Symbol{-},
\Symbol{\cap}, and \Symbol{\cup}.

\item \introLanguage{\LS}{size}, is the first order
language with the one-place predicate \tUNIT, the two-place
predicates $<$ and $\subset$, and the three place
predicate \tSUM.
\item \introLanguage{\LCS}{class size} is the first order language
containing just the non-logical constants in \LC and \LS.
\endEnum
\endDef
Following the strategy outlined above, we define the theory of class size in terms of interpretations of \LCS over finite power sets.

\Def{3.1.2}
\Enum
\item If $\LC \subset \L$, then \Model{A} is a \introSTerm{standard interpretation} of \L \iff
    \Enum
    \item \Domain{A} = $\PowerSet{x}$, for some $x$, and
    \item \Model{A} assigns the usual interpretations to all constants of \LC:
    \begin{align*}
    \Model{A}(\Identity) &= x,\\
    \Model{A}(\NullSymbol) &=  \NullSet,\\
    \Model{A} \Satisfies a \subset b &\tiff a \subset b
    \end{align*}
    \endEnum

\item If \Model{A} is a standard interpretation and
    \Domain{A} = \PowerSet{x}, then x is the \introSTerm{basis} of \Model{A}, \Basis{\Model{A}}.
\endEnum
\endDef

\Def{3.1.3} \Model{A} is a \introSTerm{standard finite interpretation} of \LCS \iff
    \Enum
    \item \Model{A} is a standard interpretation of \LCS,
    \item \Model{A} has a finite basis, and
    \item
    \begin{align*}
     \Model{A} \Satisfies a < b &\tiff
                \Card{a} < \Card{b},\\
     \Model{A} \Satisfies a \eqsize b &\tiff
                \Card{a} = \Card{b}\emath{, and}\\
     \Model{A} \Satisfies \tUNIT(a) &\tiff \Card{a} = 1
    \end{align*}
   \endEnum
\endDef

\Def{3.1.4}\introTheory{\CS}{theory of class size} is the set of all sentences of \LCS which are true in all standard finite interpretations of \LCS.
\endDef

\par
By drawing only on principles which can be stated in \LCS
we have at least ruled out the most obvious danger of paradox.
Since the notion of one-to-one correspondence cannot be expressed in this language,
Cantor's principle \OO will not be included in the
theory \CS, even though it is true over any finite power set.
\par
Since CS has arbitarily large  finite models, it has infinite models.
It is not obvious that \CS has \emph{standard} infinite models, in which the universe is an infinite power set.
In Chapter 6 we show that such models do exist.
\par
The present chapter is devoted to getting a clearer picture of the theory \CS.
Section 2 develops a set of axioms, \CA, for \CS.
Section 3 outlines the proof that \CA does indeed axiomatize \CS.
This proof is presented in Chaper~5, after a slight detour in Chapter~4.

\Section{3.2}{CA - Axioms for CS}

Here we shall develop a set of axioms, \CA, for the theory \CS.
This will be done in several stages.
\par
\Subsection{3.2.1}{\BA - Axioms for atomic boolean algebra}
\par
We'll begin with the obvious.
Since all of the universes of the interpretations mentioned in the definition of \CS are power sets, they must be atomic boolean algebras and, so, must satisfy
\BA:
\Def{3.2.1} \introTheory{\BA}{theory of atomic boolean algebras} is the theory consisting of the following axioms:
\begin{gather*}
x \cup y = y \cup x\\
x \cap y = y \cap x\\
x \cup (y \cup z) = (x \cup y) \cup z)\\
x \cap (y \cap z) = (x \cap y) \cap z)\\
x \cap (y \cup z) = (x \cap y) \cup (y \cap z)\\
x \cup (y \cap z) = (x \cup y) \cap (y \cup z)\\
x \cap (\Identity - x) = \NullSet\\
x \cup (\Identity - x) = \Identity\\
x \subset y \bicond (x \cup y) = y \Aand x \ne y\\
\tATOM(x) \bicond (\forall y)(y \subset x \then y = \NullSet)\\
x \ne \NullSet \then (\exists y)(\tATOM(y) \Aand
        (y \subset x \Or y = x)
\end{gather*}
\endDef
\par
These axioms are adapted from \cite[Def 9.3, p. 141 and Def. 9.28, p. 151]{Monk}.
\par
\BA is clearly not a complete axiomatization of \CS, since
\BA does not involve any size notions.
But \BA does entail all the sentences in \CS which do not
themselves involve size notions.
To show this we need to draw on some established facts
about the complete extensions of \BA (in the language \LC).
The key idea here is that complete extensions of \BA can
be obtained either by stipulating the finite number of atoms in a model or by saying that there are infinitely many atoms.

\Def{3.2.2} For $n \ge 1$,
\Enum

\item $\ATLEASTn$ is a sentence which says that there are
    at least $n$ atoms:
    \[
    \exists x_1 ... \exists x_n(\tATOM(x_1) \Aand ... \Aand
    \tATOM(x_n)
    \Aand \theSetst{x_i \neq x_j}{0 < i < j \leq n})
    \]
\item $\EXACTLYn$ is a sentence which says that there are exactly $n$ atoms:
    \[
    \ATLEASTn \Aand \neg \ATLEAST{n+1}
    \]
\item $\INF$ is a set of sentences which is satisfied in all
    and only infinite models of \BA:
    \[
    \INF \eqdef \theSetst{\ATLEASTn}{n > 1}
    \]
\item $\BAn \eqdef \BA \adjoin \EXACTLYn$
\item $\BAI \eqdef \BA \cup \INF$
\endEnum
\endDef
\Fact{3.2.3}
\Enum
\item For $n \geq 1$, \BAn is categorical. \cite[Cor 9.32, p.152]{Monk}
\item For $n \geq 1$, \BAn is complete. (Immediate from above)
\item For $n \geq 1$, an \n-atom atomic boolean algebra is
isomorphic to any standard finite interpretation of \LC with
an \n-element basis. \cite[Prop. 9.30, p.151]{Monk}
\item \BA is complete. \cite[Theorem 21.24, p.360]{Monk}
\item \BAI and the theories \BAn, $n \geq 1$, are the only complete, consistent extensions of \BA.
\endEnum
\endFact
\Fact{3.2.4}
If $\BAI \Proves \gphi$, then \gphi is true in some finite
model of \BA.
\endFact
\Proof
$\BA \union \INF \Proves \gphi$.
By compactness, then, there is a \k such that
\[
    \BA \union \theSetst{\ATLEASTn}{1 \leq n \leq k}
    \Proves \gphi
\]
So \gphi is true in any atomic boolean algebra with
more than \k atoms.
\endProof
\Theorem{3.2.5}
If $\CS \Proves \gphi$, and $\gphi \in \LC$, then $\BA \Proves \gphi$.
\endTheorem
\Proof
If $\BA \nProves \gphi$, then \notphi is true is some
atomic boolean algebra, \ModelA.
\par
If \ModelA is finite, then \ModelA is isomorphic to some
standard finite interpration \ModelAp of \LC.
But, then \notphi is true in \ModelAp, so \gphi is not true
in \ModelAp and $\gphi \notin \CS$.
\par
If \ModelA is infinite, then \notphi is consistent with \BAI.
But \BAI is complete, so $\BAI \Proves \notphi$.
By 3.2.4, \notphi is true in some finite model \ModelA of \BA.
Hence, \notphi is true in some standard finite interpretation
of \LC and, again, $\gphi \notin \CS$.
\endProof

\Subsection{3.2.2}{Size principles}
Here, we just gather the principles presented above as \EXCORE:
\aTheory{3.2.6}{CORE}
\SIZE consists of the following axioms:
\begin{align*}
    x \subset y &\then x < y     \tag{\SUBSET}  \\
    x >y &\bicond y < x          \tag{\DEFGT}\\
    x \eqsize y &\bicond \tINDISC(x,y) \tag{\DEFEQ}\\
    \neg( x &< x )               \tag{\IRREFLT} \\
    x < y \Or x &\eqsize y \Or x > y \tag{\TRICH} \\
        \tUNIT(x) &\bicond \tATOM(x)        \tag{\DEFID}\\
    \tSUM(x,y,z) &\bicond
        \exists \xP \exists \yP (  \tag{\DEFPL} \\
        &\quad \quad x \eqsize \xP \Aand y \eqsize \yP  \\
       &\quad  \Aand \xP \cap \yP = \NullSet  \\
       & \quad \Aand \xP \cup \yP = z)\\
    ( x \eqsize \xP \Aand y \eqsize \yP   \tag{\DISJU}     \\
        \Aand \xP \cap \yP = \NullSet     \\
         \Aand \xP \cup \yP = z)
         &\then x \cup y \eqsize \xP \cup \yP
\end{align*}
\endTheory
Combining the principles of boolean algebra and the size principles,
we obtain our first serious attempt at a general theory of size:
\Def{3.2.7} \introTheory{\BASIC}{basic theory} is defined as:
\[
\BA \union \SIZE
\]
\endDef

\Subsection{3.2.3}{Division principles}
\BASIC is not a complete axiomatization of \CS.
In this section we shall exhibit an infinite number of principles which need to be added to \BASIC in order to
axiomatize \CS.
When we are done, we will have an effective set of sentences,
\introTheory{\CA}{Class Size Axioms} though
we will not prove that $\CA \equiv \CS$ until Chapter 5.
\par
To show that the new principles really do need to be
added, we'll need some non-standard models of \BASIC.
These models will be similar in that (1)their universes will
be subsets of \PowerN, (2) their atoms will be singletons in \PowerN, and (3) all boolean symbols will receive their usual interpretations.
The models will, however, include different subsets of $\Naturals$ and  assign different size orderings to
these sets.
\par
In Chapter 6, these models of \BASIC will reappear as
submodels of various standard models of \CS over \PowerN.
So, in addition to the immediate purpose of establishing
independence results, these models provide a glimpse of
how sets of natural numbers are ordered by size.
\par
Every standard finite interpretation, \ModelA, of \LCS
satisfies exactly one of the following:
\begin{gather*}
    \exists x \exists y (x \eqsize y
    \Aand x \cap y =
    \NullSet \Aand x \cup y = \Identity)
    \tag{\Theory{EVEN}}\\
\begin{split}
\exists x \exists y \exists y
    (
    x \eqsize y \Aand
     x \cap y = x \cap z = y \cap z = \NullSymbol \\
     \Aand \tATOM(z) \Aand x \cup y \cup z = \Identity
    )
\end{split} \tag{\Theory{ODD}}
\end{gather*}
$\ModelA \Satisfies \EVEN$ if \Card{\Basis{A}} is even and
$\ModelA \Satisfies \ODD$ if \Card{\Basis{A}} is odd.
But $\BASIC \nProves (\EVEN \Or \ODD)$.
Consider the model \Model{F} whose universe consists of
all only the finite and cofinite subsets of \Naturals, where, for $a, b \in \Domain{F}$:
\begin{align*}
\Model{F} \Satisfies (\Constant{a} < \Constant{b})
\tiff& a \emath{and} b \emath{are both finite and}
\Card{a} < \Card{b}\\
&\emath{or} a \emath{and} b \emath{are both cofinite and}
\Card{\fN -a} >  \Card{\fN - b}\\
&\emath{or}a \emath{is finite and} b \emath{is infinite.}
\end{align*}
\par
\Model{F} is a model of \BASIC.
But neither \EVEN nor \ODD is true in \Model{F},
for any two sets that are the same size are either both finite,
in which case their union is also finite,
or both cofinite, in which case they cannot be disjoint.
\par
So $\EVEN \Or \ODD$ is in \CS, but not entailed by \BASIC.
As you might suspect, this is just the tip of the iceberg of
principles missing from an axiomatization of \CS.
Informally, we can extend \Model{F} to a model that
satisfies \EVEN by including the set of even numbers and the set of odd numbers and making them the same size.
To round out the result to a model of \BASIC, we also
need to include all sets which are \notion{near} the set
of evens or the set of odds, i.e. those that differ from
the evens or odds by a finite set.
With these additions made, the new model will be closed
under boolean operations and will satisfy \BA.
There is a (unique) way of ordering these added sets by
\notion{smaller than} that will satisfy \BASIC:
rank them according to the size and direction of their
finite difference from the odds or the evens.
So, we can construct an infinite model of \BASIC \adjoin
(\EVEN \Or \ODD).
\par
But this model will still not satisfy \CS, as we can
see by generalizing the argument above.
$(\EVEN \Or \ODD)$ says that the universe is \introSTerm{roughly divisible} by two:
\EVEN says that the universe is divisible by two without
remainder.
\ODD says that there is a remainder of a single atom.
We can construct a similar statement that says the
universe is roughly divisible by three --- with remainder 0, 1, or 2.
As with $(\EVEN \Or \ODD)$, this statement is in \CS
but is satisfied by neither our original model
nor the model as amended.
Again, we can extend the model and again we can produce
a statement of \CS which is false in the resulting model.
\par
We  now formalize this line of reasoning.

\Def{3.2.8}
If $0 \leq m < n$, then
\Enum
\item \MODn{m} is the sentence
\begin{align*}
    \exists x_1 ... \exists x_n
    \exists y_1 \dots \exists y_m(\\
& x_1 \eqsize x_2 \eqsize \dots \eqsize x_n      \\
& \Aand \tATOM(y_1) \Aand \dots \Aand \tATOM(y_m)\\
& \Aand \Andst{(x_i \intersect x_j = \emptyset)}{1 \leq i < j \leq n}   \\
&    \Aand \Andst{(y_i \intersect y_j = \emptyset)}{1 \leq i < j \leq m}\\
& \Aand (x_1 \union \dots x_n)
            \intersect (y_1 \union \dots y_m) = \emptyset   \\
& \Aand (x_1 \union \dots x_n)
            \union (y_1 \union \dots y_m) = \Identity)
\end{align*}
\MODn{m} says that the universe is divisible into $n$ sets of the same size with
$m$ atoms remaining.
\item \DIVn is the sentence
\[
    \MODn{0} \Or \dots \Or \MODn{n-1}
\]
\endEnum
\endDef
\Fact{3.2.9}
If $0 \leq m < n$ and \ModelA is a standard finite interpretation of
\LCS, then
\begin{align*}
\ModelA &\Satisfies \MODn{m} \tiff \Card{B(\ModelA)} \Cmod{m}{n}\tag{a}\\
\ModelA &\Satisfies \DIVn\tag{b}\\
\CS &\Proves \DIVn\tag{c}
\end{align*}
\endFact
\Theorem{3.2.10}
    $\BASIC \nProves \CS$
\endTheorem
\Proof
    If $n>1$, $\BASIC \nProves \DIVn$.
    The model, \Model{F} defined above satisfies \BASIC but not
    \DIVn, for any \n finite sets have a finite union and any
    two cofinite sets overlap.
    So, $\BASIC \nProves \CS$, by 3.2.9c.
\endProof

We could consider adding all \DIVn sentences to \BASIC in the
hope that this would yield a complete set of axioms for \CS.
We did consider this, but it does not work.
To demonstrate this, we need some independence results for
sets of \DIVn sentences.
\Def{3.2.11} If $J$ is a set of natural numbers, then
\Enum
\item   $\DIVJ = \theSetst{\DIVn}{n \in J}$
\item   $\BDIVJ = \BASIC \union \DIVJ$
\item   $\BDIVj = \BDIV{\theSet{j}}$
\endEnum
\endDef
Our independence results will be obtained by constructing
models of \BASIC which satisfy specific sets of DIV sentences.
To build such models from subsets of \Naturals, we shall
include sets which can be regarded as fractional portions of
\fN.
\Def{3.2.12}
    For $n \geq 0$,
\Enum
\item   \x is an \introMTerm{$n$-congruence class}{n-congruence class} \iff
        \x = \MSk{n}{m} for some \m, where
        $0 \leq m < n$.
\item   \x is an \introMTerm{$n$-quasi-congruence class}{n-quasi-congruence class} \iff
        \x is the union of finitely many \n-congruence classes.
\item   \x is a \introSTerm{congruence class} \iff
        \x is an n-congruence class for some \n.
\item   \x is a \introSTerm{quasi-congruence class} \iff
        \x is an n-quasi-congruence class for some \n.
\item
        \QCn = \theSetst{x}{x \emath{is an n-quasi-congruence class}}
\item   \QC = \Unionst{\QCn}{n > 0}
\endEnum
\endDef
\Examples
\Enum
\item
    The set of evens, \evens, and the set of odds, \odds, are
    2-congruence classes.
\item
    \MSk{3}{2} is a 3-congruence class.
\item
    $\fN$ is a 1-congruence class.
\item
    $\fN$ is an \n-quasi-congruence class for every $n > 0$:
    \[
        \fN = \MSk{n}{0} \union \dots \union \MSk{n}{n-1}
    \]
\endEnum

\endExamples
\Fact{3.2.13}
\Enum
\item   If $x \in \QCn$ and $y \in \QCn$, then
        $x \union y \in \QCn$.
\item   If $x \in \QCn$, then
        $\Naturals - x \in \QCn$.
\endEnum
\endFact
\Proof
\Enum
\item
    Suppose
    $x = a_1 \union \dots a_k$ and
    $y = b_1 \union \dots b_j$.
    Then
    $x \union y = a_1 \union \dots a_k \union
    b_1 \union \dots b_j$.
\item
    $\fN$ itself is the union of \n-congruence
    classes.
    If \x is the union of \m of these classes, then
    $\fN - x$ is the union of the remaining $(n-m)$
    \n-congruence classes.
\endEnum
\endProof
\Def{3.2.14}
    $x$ is \introSTerm{near} $y$, NEAR(x,y), \iff $x - y$
    and $y - x$ are both finite.
\endDef
\Fact{3.2.15}
    $x$ is near $y$ \iff there are finite sets, $w_1$ and $w_2$
    such that $x = (y \union w_1) - w_2$.
\endFact
\Proof
\BothWays
\IfWayB
    Let $w_1 = x -y$ and let $w_2 = y-x$.
\OnlyWayB
    Suppose $x = (y \union w_1) - w_2$.
    Then $x -y \subseteq w_1$ and $y-x \subseteq w_2$,
    so $x-y$ and $y-x$ are finite.
\enBoth
\endProof

\Fact{3.2.16}
    If $x_1 \subseteq x \subseteq x_2$, $x_1$ is near \y, and
    $x_2$ is near \y, then \x is near \y.
\endFact
\Proof
    Since $x \subseteq x_2$, $(x-y)\subseteq(x_2-y)$.
    But $(x_2-y)$ is finite, so $(x-y)$ is finite.
    \par
    Since $x_1 \subseteq x$, $(y-x)\subseteq(y-x_1)$.
    But $(y-x_1)$ is finite, so $(y-x)$ is finite.
\endProof

\Fact{3.2.17}
    NEAR is an equivalence relation.
\endFact
\Proof
\Enum
\item
    \x is near \x, since $x-x =\NullSet$, which is finite.
\item
    If \x is near \y, then \y is near \x. Immediate.
\item
    Suppose \x is near \y and \y is near \z.
    Note that
    \[
        (z-x)=((z \intersect y)-x)\union((z-y)-x)
    \]
    But $(z-y)-x$ is finite because $(z-y)$ is finite
    and $((z\intersect y) - x)$ is finite because
    $((z\intersect y) - x) \subseteq (y-x)$, which is
    finite.
    So the union, $(z-x)$, is finite.
    \par
    Similarly,
    \[
        (x-z)=((x\intersect y)-z)\union((x-y)-z)
    \]
    where $((x\intersect y) - z) \subseteq (y-z)$
    and $((x-y)-z \subseteq (x-y)$.
    So $(x-z)$ is also finite.
    \par
    Hence, \x is near \z.
\endEnum
\endProof

\Fact{3.2.18}
\Enum
\item
    If $x_1$ is near $x_2$, then
    $x_1 \union y$ is near $ x_2 \union y$.
\item
    If $x_1$ is near $x_2$ and $y_1$ is near $y_2$, then
    $x_1 \union y_1 $ is near $ x_2 \union y_2$.
\item
    If \x is near \y, then $\fN - x$ is near $\fN - y$
\endEnum
\endFact
\Proof
\Enum
\item
    Since $x_1$ is near $x_2$,
    $x_1-x_2$ and $x_2-x_1$ are finite.  But
    \begin{align*}
    (x_1 \union y) - (x_2 \union y) &\subseteq (x_1-x_2)
      \\
      \emath{and} (x_2 \union y) - (x_1 \union y) &\subseteq (x_2-x_1)
    \end{align*}
\item
    By (a), $x_1 \union y_1$ is near $x_2 \union y_1$ , which
    is near $x_2 \union y_2$.
    So $x_1 \union y_1$  is near $x_2 \union y_2$, by
    transitivity.
\item
    $(\fN - x) - (\fN -y) = y-x$ and
    $(\fN - y) - (\fN -x) = x-y$.
    So if \x and \y are near each other, so are their
    complements.
\endEnum
\endProof

We can now define the domains of the models we will use to establish the
independence results.

\Def{3.2.19}
\Enum
\item
    \Qn = \theSetst{y}{y \emath{is near an \n-quasi-congruence
                        class}}
\item
    \Q = \Unionst{\Qn}{n > 0}
\endEnum
\endDef

\Examples
\Enum
\item
    $Q_1 = \theSetst{y \subseteq \fN }{\emath{\y is finite or cofinite}}$.
\item
    $Q_2 = \theSetst{y\subseteq \fN}{
        \emath{\y is  finite, cofinite, near \evens or near \odds}}$.
    \par
    $Q_2$ is the domain of the model constructed above to
    satisfy \DIV{2}.
\endEnum

\endExamples

\Fact{3.2.20} If \A is a class of sets such that
\Enum
\item
    $\Union A \in A$,
\item
    If $x \in \A$, then $\Union A - x \in A$, and
\item
    If $x \in A$ and $y \in A$, then $x \union y \in A$,
\endEnum
then \A forms a boolean algebra under the usual set-theoretic
operations, where $\Identity$ is interpreted as \A.
\cite[Def. 9.1, p.141 and Corr 9.4, p.142]{Monk}
\endFact
\Theorem{3.2.21}
If $n > 0$, then \Qn forms an atomic boolean
algebra under the usual set-theoretic operations.
\endTheorem
\Proof
Using Fact 3.2.20:
\Enum
\item
    $\Union{\theSet{\Qn}} = \fN$,
    since $\fN \in \theSet{\Qn}$ and
    if $x \in \Qn$, then $x \subseteq \fN$.
\item
    If $x \in \Qn$, then
    \begin{flalign*}
    \Suppose   x \isNear y    \\
    \AND    y \in \QCn\\
    \So     \fN - y  \in \QCn \by{3.2.13b}\\
    \AND    \fN - x \isNear \fN - y
                        \by{3.2.18b}\\
    \So     \fN - x \in \Qn.
    \end{flalign*}
\item
    \begin{flalign*}
    \Suppose x \in \Qn      \\
    \AND    y \in \Qn       \\
    \Then   x \isNear \xP   \\
    \AND    \xP \in \QCn \ForSome{\xP}\\
    \AND    y \isNear \yP   \\
    \AND    \yP \in \QCn \ForSome{\yP}\\
    \Butn   \xP \union \yP \in \QCn \by{3.2.13a}\\
    \AND    x \union y \isNear \xP \union \yP
                        \by{3.2.18b}\\
    \So     x \union y \in \QCn
    \end{flalign*}
\endEnum
\endProof
Thus, \Qn is a boolean algebra.
Moreover, every singleton is in $\Qn$ since all singletons are near \NullSet.
So \Qn is an atomic boolean algebra.
\par
We now define a size function on all sets which are near
quasi-congruence classes.
The sizes assigned to sets are ordered pairs.
The first member is a rational between 0 and 1 which represents
the density of the set.
The second member is an integer which represents the finite
(possibly negative) deviation of a set from \scare{average} sets of
the same density.
First, we define the ordering and arithmetic for these sizes
with the intention of inducing the size ordering and \tSUM relation
for sets from the assignment of sizes to sets.

\Def{3.2.22}
\Enum
\item
    A \introSTerm{size} is an ordered pair \Pair{\rho}{\delta}, where
    $\rho$ is rational and $\delta$ is an integer.
\item
    If $\theta_1 = \rd{1}$ and $\theta_2 = \rd{2}$ are
    sizes, then
    \Enum
    \item   $\theta_1 < \theta_2$ \iff
                $\rho_1 < \rho_2$
                or $(\rho_1 = \rho_2 \Aand \delta_1 < \delta_2)$.
    \item   $\theta_1 + \theta_2 =
                \Pair{\rho_1 + \rho_2}{\delta_1 + \delta_2}$
    \endEnum
\endEnum
\endDef
\par
Only some of these sizes will actually be assigned to sets.
Specifically, a size will be assigned to a set only if
$0 \leq \rho \leq 1$.
Moreover, if $\rho = 0$, then $\delta \geq 0$ and
if $\rho = 1$, then $\delta \leq 0$.
\par
Our intention in assigning sizes to sets is as follows:
Suppose \x is near an \n-quasi-congruence class, \xP.  So \xP
is the union of $\k \leq \n$ \n-quasi-congruence classes.
The set \xP has density $k/n$ and this is the value, $\rho$,
assigned to \x.
The $\delta$ value assigned to \x is the finite number of elements
added to or removed from \xP to obtain \x.

The definitions and facts below formalize this intention and
demonstrate that the assignment of sizes to sets is well-defined.

\Fact{3.2.23}
\Enum
\item
    If $x \in \QC$, $y \in \QC$, and $x \neq y$, then \x is
    not near \y. (That is, no two quasi-congruence classes are
    near each other.)
\item
    Any set is near at most one quasi-congruence class.
\endEnum
\endFact
\Proof
\Enum
\item
    Let \n be the least number such that
    \[
        x \in \QCn \emath{and} y \in \QCn
    \]
    So each is a disjoint union of \n-congruence classes:
\begin{align*}
    x &= x_1 \union \dots \union x_j \\
\emath{and} y &=  y_1 \union \dots \union y_k
\end{align*}
\begin{flalign*}
\Suppose    a \in x-y       \padleft\\
\Then       a \in x_i \emath{for some} i.\\
\But        a \notin y_j\emath{, for any \y.}\\
\So         x_i \ne y_j\emath{, for any} j\\
\So         x_i \subseteq x - y\\
\So         x-y \emath{is infinite, since} x_i \emath{is infinite} \\
\end{flalign*}
Similarly, if $a \in y-x$, then $y-x$ is infinite.
\item
    If \x were near two quasi-congruence classes, the two
    would have to be near each other, since NEAR is transitive.
    But this is impossible by (a).
\endEnum
\endProof

\Def{3.2.24} If \x is near a quasi-congruence class, then
\Enum
\item
    $C(x)$ is the quasi-congruence class near \x.
\item
    $\Dox = x - C(x)$.
\item
    $\Dtx = C(x) - x$.
\endEnum
\Dox and \Dtx are finite and
\[
    x = (C(x) \union \Dox) - \Dtx
\]
\endDef

\Def{3.2.25}
    If $x \in \QC$, then
\Enum
\item
    \ax = the least \n such that $x \in \QCn$.
\item
    \bx =  the unique \k
        such that \x is the disjoint union of \k
        \ax-congruence classes.

\endEnum
\endDef

\Examples
\Enum
\item
    If $x = \odds$, $\ax = 2$ and $\bx = 1$.
\item
    If $x = \MSn{4}{1} \union \MSn{4}{2}$, $\ax = 4$ and $\bx = 2$.

\item
    If $x = \MSn{4}{1} \union \MSn{4}{3}$, $\ax = 2$ and $\bx = 1$,
    since $x = \odds$.

\endEnum
\endExamples

\Def{3.2.26} If $x \in \Q$, then
\begin{align*}
    \rho(x) &= \frac{\beta(C(x))}{\alpha(C(x))}\\
    \delta(x) &= \Card{\Dox} - \Card{\Dtx}\\
    \theta(x) &= \Pair{\rho(x)}{\delta(x)}
\end{align*}
\endDef

We can, at last, define the models to be used in our independence
proof.
\Def{3.2.27}
    For $n > 0$, \mQn is the interpretation \mQ of \LCS
    in which boolean symbols receive their usual interpretation
    and
\begin{align*}
    \DomainQ &= \Qn\\
    \mQ \Satisfies x < y &\tiff \theta(x) < \theta(y)\\
    \mQ \Satisfies x \eqsize y &\tiff \theta(x) = \theta(y)\\
    \mQ \Satisfies \tUNIT(x) &\tiff \theta(x) = \Pair{0}{1}\\
    \mQ \Satisfies \tSUM(x,y,z) &\tiff \theta(z) = \theta(x) + \theta(y)
\end{align*}
\endDef
To show that the models \mQn satisfy \BASIC, we will need the
following facts about congruence classes.

\Fact{3.2.28}
\Enum
\item
    If $x = \MSn{a}{b}$ and $a_2 = a c$, then
    \[
        x = \Unionst{\MSn{a_2}{(i a + b)}}{0 \leq i < c}
    \]
 
\item
    If $x \in QC_n$ and $m = kn$, then $x \in QC_m$.
\item
    If $x \in QC$ and $y \in QC$, there is an \n such that
    $x \in QC_n$ and $y \in QC_n$.

\endEnum
\endFact

\Proof
\Enum
\item
    If $k \in \MSn{a_2}{(i a + b)}$ for some $i$, $0 \leq i < c$, then, for some
    $n_1$,
    \begin{align*}
    k &= a_2 n_1 + i a + b    \\
    &= a (c n_1) + i a + b \\
    & = a (c n_1 + i) + b
    \end{align*}
    So $k \in x$.
    \par
    If $k \in x$, there is an $n_1$ such that $k = a n_1 +b$.
    Let $n_2$ be the greatest $n$ such that $a_2 n \leq k$ and let
    $k^\prime = k - a_2 n_2$.
    \begin{flalign*}
    \Since k \Cmod{b}{a}\\
    \AND   a_2 n_2 \Cmod{0}{a}    \padleft \\
    \Then   k^\prime \Cmod{b}{a}\\
    \So k^\prime = i a + b \emath{, where} 0 \leq i < c\\
    \But    k = a_2 n_2 + k^\prime = a_2 n_2 + i a + b\\
    \So     k \in  \MSn{a_2}{(i a + b)}
    \end{flalign*}
\item
    By (a), each \n-congruence class is a disjoint union of
     \m-congruence classes.
\item
    Suppose $x \in \QC_{n_1}$ and $y \in \QC_{n_2}$.
    Then, by (b), both \x and \y are in $\QC_{n_1 n_2}$.
\endEnum
\endProof

\Theorem{3.2.29}
    For any $n > 0$, $\mQn \Satisfies \BASIC$
\endTheorem
\Proof
    By 3.2.21, \mQn is an atomic boolean algebra; so
    $\mQn \Satisfies \BA$.
    The $<$-relation of \mQn is induced from the linear ordering
    of sizes; so it is a quasi-linear ordering and \IRREF,
    \TRICH, and \DEFEQ are satisfied.
    As for the remaining axioms:
    \Enum
    \item \SUBSET: Suppose $x \subset y$.
    If $C(x) = C(y)$
    then $\rho(x) = \rho(y)$
    and $\Dox \subseteq \Doy$ and $\Dtx \subseteq \Dty$, where at least one
    of these inclusions is proper.
    So $\delta(x) < \delta(y)$.
    \par
    But if $C(x) \neq C(y)$, then $C(x) \subset C(y)$, so $\rho(x) < \rho(y)$.
    \par
    In either case, $\theta(x) < \theta(y)$, so $\mQn \Satisfies x < y$.
    \item \REPLT:
    Suppose $\mQn \Satisfies x < y$.
    So $\theta(x) = \Pair{k_1/n}{\delta_1}$,
      $\theta(x) = \Pair{k_2/n}{\delta_2}$, and either $k_1 < k_2$ or
      $\delta_1 < \delta_2$.
      \par
      We want to find some $\xP \subset y$
      such that $\mQn \Satisfies x \eqsize \xP$.
      \par
      If $k_1 = k_2 > 0$, then \y must be infinite; so, \xP can be obtained
      by removing $\delta_2 - \delta_1$ atoms from \y.
      \par
      If $k_1 = k_2 = 0$, then $0 \leq \delta_1 < \delta_2$; so, again,
      \xP can be obtained
      by removing $\delta_2 - \delta_1$ atoms from \y.
      \par
      If $k_2 > k_1 > 0$, then let $y_1$ be the union of $k_1$ $n$-congruence
      classes contained in $C(y)$.
      So $y - y_1$  is infinite and $y_1 - y$ is finite.
      Let $y_2 = y_1 - (y_1 - y) = y_1 \cap y$.
      So $\theta(y_2) = \Pair{k_1/n}{-\delta_3}$
       where $\delta_3 = \Card{y_1 - y}$.
       Finally, let $\delta_4 = \delta_3 + \delta_1$ and
       $\xP = y_2 \cup y_3$, where
       $y_3 \subseteq y_1 - y$ with $\Card{y_3} = \delta_4$.
       \par
       If $k_1 = 0$, then $x$ is finite.
       If $k_2 > 0$, then $y$ is infinite, so there is no problem.
       If $k_2 = 0$, then $y$ is finite, but has more members than $x$, since
       $\delta_1 < \delta_2$. So, let \xP be any proper subset of $y$ with
       $\delta_1$ members.
    \item \DISJU:
    It is enough to show that if \x and \y are disjoint, then
    $\theta(x \cup y) = \theta(x) + \theta(y)$.
    We need the following three facts:
    \Enum
    \item
    $\C(x \cup y) = \C(x) \cup \C(y)$. (See Fact 3.2.18b.)
    \item
    $\Delta_1(x \cup y) = (\Delta_1(x) \cup \Delta_1(y)) - (C(x) \cup C(y))$.
    \par
    If $z \in x \cup y$ but $a \notin C(x \cup y)$, then
    $a \in \Delta_1(x)$ or $a \in \Delta_1(y)$;
    any element of $\Delta_1(x)$ is also in $\Delta_1(x \cup y)$
    unless it is in $C(y)$;
     any element of $\Delta_1(y)$ is also in $\Delta_1(x \cup y)$
    unless it is in $C(x)$.
    \item
    $\Delta_2(x \cup y) = (\Delta_2(x) \cup \Delta_2(y)) -
        (\Delta_1(x) \cup \Delta_1(y))$.
    \par
    Note that if $x \in \Q$, $y \in \Q$, and $x \cap y = \NullSet$,
    then $C(x) \cap C(y) = \NullSet$; otherwise, $C(x)$ and $C(y)$ have
    an infinite intersection.
    \par
    \begin{flalign*}
    \Henceif    a \in  C(x) \cup C(y) - (x \cup y)  \padleft\\
    \Then       a \in C(x) - x = \Dtx \\
    \OR         a \in C(y) - y = \Dty \\
    \Andif      a \in \Dtx, \etext{ then } a \in \Delta_2(x \cup y) \\
    \unless     a \in \Doy \\
    \Andif      a \in \Dty, \etext{ then } a \in \Delta_2(x \cup y) \\
    \unless     a \in \Dox \\
    \end{flalign*}.
    \endEnum
    From (ii) we obtain (iia):
    \begin{multline*}
    \Card{\Delta_1(x \cup y)} =
        \Card{\Delta_1(x) \cup \Delta_1(y)}   \\
            -  \Card{
                    (\Delta_1(x) \cup \Delta_1(y)
                    \cap (C(x) \cup C(y))
                }     \tag{iia}
   \end{multline*}
     and from (iii) we obtain (iiia):
    \begin{multline*}
   \Card{\Delta_2(x \cup y)} =
        \Card{\Delta_2(x) \cup \Delta_2(y)}       \\
            -  \Card{\Delta_1(x) \cup \Delta_1(y)
                    \cap (\Delta_2(x) \cup \Delta_2(y) ) }     \tag{iiia}
    \end{multline*}
    But $\Delta_1(x)$ and $\Delta_1(y)$ are disjoint, since \x and \y are
    disjoint. So:
    \begin{gather*}
        \Card{\Delta_1(x) \cup \Delta_1(y)}
            = \Card{\Delta_1(x)} + \Card{\Delta_1(y)}
                = \delta_1(x) + \delta_1(y)
            \tag{iv}
    \end{gather*}
    Since $\Delta_2(x)$ and $\Delta_2(y)$ are contained, respectively,
    in $C(x)$ and $C(y)$, which are disjoint, they are also disjoint.
    So:
    \begin{gather*}
      \Card{ \Delta_2(x) \cup \Delta_2(y)}
        = \delta_2(x) + \delta_2(y)      \tag{v}
    \end{gather*}
    \par
    So,
    \begin{align*}
    \delta_1(x \cup y) - &\delta_2(x \cup y)     \\
        & = (\delta_1(x) +  \delta_1(y)) - (\delta_2(x)  + \delta_2(y))  \\
        & = (\delta_1(x) - \delta_2(x))  + (\delta_1(y) - \delta_2(y))  \\
        & = \delta(x) + \delta(y)
    \end{align*}
    Since $\rho(x \cup y)  = \rho(x) + \rho(y)$, by (i), we know that
    $\theta(x \cup y)  = \theta(x) + \theta(y)$.

    \item \DEFPL:
    \begin{align*}
    \Suppose       \mQn \Satisfies \tSUM(x,y,z)     \\
    \So            \theta(z) = \theta(x) + \theta(y) \\
    \Clearly       \theta(x) \leq \theta(z)  \\
    \Assume        \theta(x) = \theta(z)  \\
    \Then          \theta(y) = \Pair{0}{0}  \\
    \So            y = \NullSet \\
    \Let           \xP = z \emath{, } \yP = \NullSet \emath{to satisfy \DEFPL}\\
    \Assume        \theta(x) < \theta(z)    \\
    \Then          \mQn \Satisfies \xP \eqsize x \\
    \AND           \xP \subset z \emath{for some} \xP
                    \since{ \mQn \Satisfies \REPLT} \\
    \Let           \yP = z - \xP \\
    \So             z = \xP \Cup \yP \\
    \CLAIM         \theta(y) = \theta(\yP)
                         \emath{so }\mQn \Satisfies y \eqsize \yP \\
    \For           \theta(\xP) + \theta(\yP) = \theta(z)
                        \by{\DISJU} \\
    \But           \theta(\xP) = \theta(x)\\
    \So            \theta(\yP) = \theta(z) - \theta(x) = \theta(y)
    \end{align*}
    (Cancellation is valid for sizes because it is valid for rationals and integers.)
    \par
    Conversely,
    \begin{flalign*}
    \If       \theta(z) = \theta(\xP) + \theta(\yP)     \\
    \AND      \theta(\xP) = \theta(x)    \\
    \AND      \theta(\yP) = \theta(y)     \\
    \Then     \theta(z) = \theta(x) + \theta(y)
     \end{flalign*}

    \endEnum

\endProof
\Theorem{3.2.30}
    For any $n > 0$, $\mQn \Satisfies \DIV{m}$ \iff $\m \divides \n$.
\endTheorem
\Proof
\BothWays
\IfWayB
    $\theta(\fN) = \Pair{1}{0}$.
    Hence, if \m disjoint sets of the same size exhaust
    \fN, they must each have size \Pair{1/m}{0}.
    But if $x \in \Qn$, then
    $\theta(x) = \Pair{a/n}{b}$, for integral $a$ and $b$.
    So $b = 0$ and $a = n/m$.
\OnlyWayB
For each $i$, $0 \leq i < n$, let $A_i = \MSk{n}{i}$.
So $\fN = \Unionst{A_i}{0 \leq i < n}$.
\par
If $i \neq j$, then $A_i \cap A_j = \NullSet$
and $\mQn \Satisfies (A_i \eqsize A_j)$
since $\theta(A_i) = \theta(A_j) = \Pair{1/n}{0}$.
\par
Letting $p = n/m$, group the $n$ sets $A_i$ into $m$ collections
with $p$ members in each:
\[
    B_1, \dots , \B_m
\]
Letting $b_j = \Union B_j$ for $1 \leq j \leq m$, we have $b_j \in \Qn$
and $\theta(b_j) = \Pair{p/n}{0} = \Pair{1/m}{0}$.
\par
Furthermore, $b_1 \cup \dots b_m = \fN$.
\enBoth
\endProof
\Def{3.2.31}
    If $J \neq \emptyset$ and $J$ is finite, then the
    \introSTerm{least common multiple} of $J$,
    $\mu(J)$, is the least \k which is divisible by
    every member of $J$.
\endDef
\Remark
$\mu(J)$ always exists since the product of all members of $J$
is divisible by each member of $J$.
Usually, the product is greater than $\mu(J)$.
\endRemark
\Corollary{3.2.32}
If $J$ is finite, then
\Enum
\item
    If $\BDIVn \Proves \DIVm$, then $m \divides n$.
\item
    If $\BDIVJ \Proves \DIVm$, then $m \divides \mu(J)$.
\item
    There are only finitely many \m for which
    $\BDIVJ \Proves \DIVm$.
\endEnum
\endCorollary
\Proof
Based on Theorem 3.3.20:
\Enum
\item
    If $m \ndivides n$, then $\mQn \Satisfies \BDIVn\adjoin\neg\DIVm$
\item
    $\Model{Q}_{\muJ} \Satisfies \BDIVJ$ since it satisfies
    \DIV{j} for each $j \in J$.
    But, if $m \ndivides \muJ$, then
    $\mQ_{\muJ} \nSatisfies \DIVm$.
\item
    Immediate from (b), since only finitely many \m divide \muJ.
\endEnum
\endProof
We are now ready to show that $\BDIVn \nProves \CS$ by
finding a sentence in \CS which entails infinitely many \DIVn
sentences.
Such sentences can be produced by generalizing the notion
of divisibility to all sets instead of applying it only to
the universe.
\Def{3.2.33}
\Enum
\item
    If $0 \le n$, then $\Timesn(x,y)$ is the formula:
    \[
        \exists x_0 \dots \exists x_n
        [x_0 = \NullSet \Aand x_n = y \Aand
            \Andst{\tSUM(x_{i-1},x,x_i)}{1 \leq i \leq n}
        ]
    \]
    $\Timesn(x,y)$ says that \y is the same size as
    the \tSUM of \n sets, each the same size
    as \x.
\item
    If $0 \le m < n$, then $\Modnm(z)$ is the formula:
    \[
        \exists x \exists y \exists v \exists w
        [
            Times_n(x,v)
            \Aand \tUNIT(y)
            \Aand Times_m(y,w)
            \Aand \tSUM(v,w,z)
        ]
    \]
    $\Modnm(z)$ says that \z can be partitioned in \n
    sets of the same size and \m atoms.
\item
    $\Divn(z)$ is the formula
    \[
        \Modn{0}(z) \Or \dots \Or \Modn{n-1}(z)
    \]
\item
    \ADIVn is the sentence
    \[
        \forall x \Divn(x)
    \]
\endEnum
\endDef

\Remark
We have taken this opportunity to formulate the divisibility
predicates purely in terms of size predicates.
Notice that in the presence of \BASIC,
\begin{align*}
    \MODn{m} &\equiv \Modn{m}(\Identity)   \\
  \emath{and }       \DIVn  & \equiv \Divn(\Identity)
\end{align*}
\endRemark
\Fact{3.2.34}
    $\CS \Proves \ADIVn$, for every \n.
\endFact
\Proof
Every set in every standard finite interpretation is a finite set, and all finite
sets are roughly divisible by every \n.
\endProof

\Fact{3.2.35}
    $\BASIC;\ADIVn \Proves$:
\Enum
\item
    \ADIV{{n^m}}, for all \m
\item
    \DIVn, and
\item
    \DIV{{n^m}}, for all \m.
\endEnum
\endFact

\Proof
\Enum
\item
    By induction on \m:
    if $m = 1$, then $n^m = n$, so
    $\ADIVn \Proves \ADIV{n^m}$.
    \par
    If $T \Proves \ADIV{n^k}$, $\mA \Satisfies T$,
    and $x \in \Domain{A}$, then \x can be partitioned into
    $n^k$ sets of the same size and \i atoms,
    where $i < n^k$.
    Each non-atomic set in the partition can be further
    partitioned into \n sets of the same size and
    \j atoms, where $j < n$.
    Thus, we have partitioned \x into $n^k n$ sets of
    the same size and $n^k j + i$ atoms.
    But $n^k n = n^{k+1}$ and, since $i < n^k$ and
    $j < n$, $n^k j + i < n^{k+1}$.
    Hence, $\mA \Satisfies \ADIV{n^{k+1}}$.
\item
    Obvious.
\item
    Immediate from (a) and (b).
\endEnum
\endProof

\Theorem{3.2.36}
    $\BDIV{\fN} \nProves \ADIVn$ for any \n > 1.
\endTheorem
\Proof
    If $\BDIV{\fN} \Proves \ADIVn$, then by compactness
    there is a finite set \J such that
    $\BDIVJ \Proves \ADIVn$.
    But then $\BDIVJ \Proves \DIV{n^k}$ for every
    \k, by Fact 3.2.35c.
    But this contradicts Fact 3.2.32, which says that \BDIVJ entails
    only finitely many \DIVn sentences.
\endProof

So, even if we add all of the DIV sentences to \BASIC,
we are left with a theory weaker than \CS.
Since this weakness has arisen in the case of ADIV sentences,
it is reasonable to attempt an axiomatization of \CS as
follows:

\Def{3.2.37}
    $\CA \equiv \BASIC \union \theSetst{\ADIVn}{n > 0}$
\endDef

The remainder of this chapter and the next two are devoted to showing
that \CA is, indeed, a complete set of axioms for \CS.
\Section{3.3}{Remarks on showing that CA axiomatizes CS}

We know that $\CS \Proves \CA$ and we want to show that $\CA \equiv \CS$,
i.e. that $\CA \Proves \CS$.
To do so, it will be sufficient to show that every consistent
extension of \CA is consistent with \CS.

\Fact{3.3.1}
\Enum
\item
    (Lindenbaum's lemma) Every consistent theory has a consistent,
    complete extension.
\item
    If every consistent, complete extension of $T_2$ is consistent
    with $T_1$, then $T_2 \Proves T_1$.
\endEnum
\endFact

\Proof
\Enum
\item
    \cite[Theorem 11.13, p.200]{Monk}
\item
    Suppose that $T_1 \Proves \gphi$,
    $T_2 \nProves \gphi$.
    Then $T = T_2 \adjoin \neg\gphi$ is consistent.
    \T has a consistent, complete extension, \Tp, by
    Lindenbaum's lemma.
    Since $T_2 \subseteq T$, \Tp is also a consistent,
    complete extension of $T_2$.
    But \Tp is not consistent with $T_1$.
\endEnum
\endProof

\Def{3.3.2}
    \Tp is a \introSTerm{completion} of \T \iff \Tp is a
    complete, consistent extension of \T.
\endDef

To prove that every completion of \CA is consistent with \CS,
we define two kinds of completions of a theory.
\Def{3.3.3}
\Enum
\item
    \Tp is a \introSTerm{finite completion} of \T \iff
    \Tp is true in some finite model of \T.
\item
    \Tp is an \introSTerm{infinite completion} of \T \iff
    \Tp is true in some infinite model of \T.
\endEnum
\endDef

\Fact{3.3.4}
    If \Tp is a completion of \T and $\BA \subseteq \T$, then
    (\Tp is a finite completion of \T \iff \Tp is not an
    infinite completion of \T.)
\endFact
\Proof
\BothWays
\IfWayB
    Suppose \mA is a finite model of \Tp with $n$ atoms.
    Since \Tp is complete,
    $\Tp \Proves \EXn$.
    So $\Tp \Proves \neg \ATLEAST{n+1}$ and has no
    infinite models.
\OnlyWayB
    $\Tp \Proves \ATLEASTn$ for every \n, so \Tp
    has no finite models.
\enBoth
\endProof
\Fact{3.3.5}
\Enum
\item
    Every finite completion of \CA is equivalent to
    \CA;\EXn, for some \n.
\item
    Every finite completion of \CA is consistent with \CS.
\endEnum
\endFact
\Proof
\Enum
\item
    $\CA \Proves \BASIC$ and, by Theorem 2.1.6,
    \BASIC is categorical in every finite power.
\item
    The model \mA of $\CA \adjoin \EXn$ is a standard
    finite interpretation.
    So $\mA \Satisfies \CS$.
\endEnum
\endProof

\Def{3.3.6}
\Enum
\item
    $\CAI = \CA \union \INF$
\item
    $\CSI = \CS \union \INF$
\endEnum
\endDef
So, to show that every completion of \CA is consistent with
\CS, we may now concentrate on showing that every completion of \CAI
is consistent with \CSI.
\par
What, then, are the completions of \CAI?
Recall that \CA entails \DIVn for every $n > 0$,
where \DIVn is
\[
    \MODn{0} \Or \dots \Or \MODn{n-1}
\]
Any completion, \T, of \CAI has to solve the disjunction
\DIVn for each $\n$; \T has to entail one of the
disjuncts.   \emph{Remainder theories} specify, for each \n,
the number of atoms remaining
when the universe is divided into \n disjoint sets of the same
size.
\Def{3.3.7}  Remainder functions and remainder theories
\Enum
\item
    $f\colon \fNp \to \fN$ is a \introSTerm{remainder function}
    \iff $0 \leq f(n) < n$ for all $n \in \Domf$.
    (Henceforth $f$ ranges over remainder functions.)
\item
    \f is \introSTerm{total} \iff $\Domf = \fNp$; otherwise \f
    is \introSTerm{partial}.
\item
    \f is \introSTerm{finite} \iff \Domf is finite.
\item
    \n is a \introSTerm{solution} for \f \iff for any $i \in \Domf$,
    $n \Cmod{f(i)}{i}$.
\item
    \f is \introSTerm{congruous} \iff for any $i, j \in \Domf$, then
    $\gcd{i}{j} \divides f(i) - f(j)$;
    otherwise, \f is \introSTerm{incongruous}.
\item
    The \introSTerm{remainder theory} specified by \f, \RTf, is the set of
    sentences
    \theSetst{\MODn{m}}{f(n) = m}.
\item
    If $T$ is a theory, $\Tf = T \union \RTf$.
\endEnum
\endDef
Chapter 5 will show that if \f is total, \CAIf is
complete and that these are the only complete extensions of \CAI.
In this section, we will show that \CAIf is consistent
just in case \CSIf is consistent.

\Fact{3.3.8}
\Enum
\item
    If \f is finite, then \f has a solution \iff \f is congruous
    \iff \f has infinitely many solutions.
    \cite[Theorem 5-11, p.80]{Griffin}
\item
    \f is congruous \iff every finite restriction of \f is congruous.
\item
    There are congruous \f without any solutions. (Let
    $f(p) = p - 1$, for all primes $p$.
    Any solution would have to be larger than every prime.)

\endEnum
\endFact

\Theorem{3.3.9}
\Enum
\item
    If \f is finite, then \CSf is consistent \iff \f is congruous.
\item
    \CSf is consistent \iff \f is congruous.
\item
    \CSIf is consistent \iff \f is congruous.

\endEnum
\endTheorem
\Proof
\Enum
\item
\IfWay
    Let $\gphi = \bigwedge \RTf$,
    Since $\CS\adjoin\gphi$ is consistent,
    there is some \n such that $\mAn \Satisfies \gphi$.
    So \n is a solution of \f and, hence, \f is
    congruous by 3.3.8a.
\OnlyWay
    If \f is congruous, \f has a solution, \n.
    So $\mAn \Satisfies \CS\adjoin\gphi$
\item
\IfWay
    For every finite restriction, \g, of \f,
    \CSg is consistent.
    By (a), each such \g is congruous. Hence, \f
    is congruous by 3.3.8b.
\OnlyWay
    Every finite restriction, \g, of \f is congruous.
    So \CSg is consistent, by (a).
    By compactness, then, \CSf is consistent.
\item
\IfWay
    If \CSIf is consistent, so is \CSf. So, by (b),
    \f is congruous.
\OnlyWay
    By compactness, it is sufficient to show that
    every finite subtheory, \T, of \CSIf is
    consistent.
    But, if \T is such a theory, then
    \[
        T \subseteq \CS_{g} \union
            \theSetst{ATLEAST{i}}{i < n}
    \]
    for some \n and some finite restriction, \g, of
    \f.
    Since \f is congruous, \g is as well, by 3.3.8b.
    So $g$ has arbitrarily large solutions and
    \CSg has finite models large enough to satisfy \T.
    Hence, \T is consistent.
\endEnum
\endProof

We now want to prove a similar theorem for \CA, our proposed
axiomatization of \CS.
To do this, we must first establish that certain sentences
are theorems of \CA.
\Lemma{3.3.10}
    If $n \divides m$,
    $0 \leq q < n$, and
    $p \Cmod{q}{n}$, then
    \[
    \CA \Proves \MODm{p} \then \MODm{q}.
    \]
\endLemma
\Proof
    Suppose $\mA \Satisfies \MODm{p}$,
    $k_1 = m/n$, and $p = k_2 n + q$.
    So, \bA can be partitioned into \m sets of the
    same size
    \begin{align*}
        b_{1,1} &\dots b_{1,n}\\
        b_{2,1} &\dots b_{2,n}\\
        b_{k_1,1} &\dots b_{k_1,n}
    \end{align*}
    and \p atoms
     \begin{align*}
        a_{1,1} &\dots a_{1,n}\\
        a_{2,1} &\dots a_{2,n}\\
        a_{k_2,1} &\dots a_{k_2,n}\\
        c_1 &\dots \dots c_q
    \end{align*}
    For $1 \leq i \leq n$, let
    \[
        B_i = \Unionst{b_{j,n}}{1 \leq j \leq k_1} \union
                \Unionst{a_{j,n}}{1 \leq j \leq k_2}
    \]
    Since $\mA \Satisfies \DISJU$,
        $\mA \Satisfies B_i \eqsize B_j$, for
        $1 \leq i,j \leq n$.
    Furthermore,
    \[
        \bA = \Unionst{B_i}{1 \leq i \leq n}
                \union (c_1 \dots c_q)
    \]
    So, $\mA \Satisfies \MODn{q}$.
\endProof

\Lemma{3.3.11}
    If $0 \leq p < q < m$, then
    \[
        \CA \Proves \MODm{p} \then \neg\MODm{q}
    \]
\endLemma
\Proof
Suppose $\mA \Satisfies \MODm{p} \Aand \neg\MODm{q} $.  Then:
\[
    \mA \Satisfies x_1 \cup \dots x_m \cup a_1 \cup \dots a_p = I
\]
\[
    \mA \Satisfies y_1 \cup \dots y_m \cup b_1 \cup \dots b_q = I
\]
where the $a$'s and $b$'s are atoms and the $x$'s ($y$'s) are disjoint
sets of the same size in \mA.
Let
\begin{align*}
    X &= x_1 \cup \dots \cup x_m     \\
    Y &= y_1  \cup \dots \cup y_m    \\
    A &=  a_1  \cup \dots \cup a_p          \\
    B &=  b_1  \cup \dots \cup b_p          \\
    \BP &=  b_{p+1}  \cup \dots \cup b_q
\end{align*}
We claim that $y_1 < x_1$. For, if $y_1 \eqsize x_1$, then
$X \cup A \eqsize Y \cup B$
and if $x_1 < y_1$, then
$X \cup A < Y \cup B$;
neither is possible since
$Y \cup B \subset I = X \cup a$.
So $y_i < x_i$ for
$1 \leq i \leq m$. Since $\mA \Satisfies \REPLT$, there is a proper
subset $y^\prime_i$ of $x_i$ which is the same size at $y_i$.
\begin{flalign*}
\Let    z_i = x_i - y^\prime_i \\
\AND     Y^\prime = y^\prime_i \cup \dots \cup y^\prime_m    \\
\So     Y^\prime \cup Z = X = I - A\\
\AND    Y^\prime \cup \B^\prime  = I -B \emath{.}
\end{flalign*}
\par
But $A \eqsize B$, since each is the disjoint union of $p$ atoms.
Thus $(I-A) \eqsize (I-B)$, by \RCEQ, so $\YP \cup Z \eqsize Y \cup \BP$.
But $\YP \eqsize Y$, since for each component of $Y$, there is a component of
\YP of the same size.
So $Z \eqsize \BP$, by \RCEQ.
\par
But $Z$ must be larger than \BP, for \BP is the union of fewer
than $m$ atoms while $Z$ is the union of $m$ non-empty sets.
So, the original supposition entails a contradiction.

\endProof

\Theorem{3.3.12}
\Enum
\item   If \f is finite, then \CAf is consistent \iff \f is congruous.
\item   \CAf is consistent \iff \f is congruous.
\item   \CAIf is consistent \iff \f is congruous.

\endEnum
\endTheorem
\Proof
\Enum
\item

\IfWay
    Supposing that \f is incongruous, there exist
    \i, $j$, and \k such that $k = \gcd{i}{j}$ and
    $k \ndivides (f(i) - f(j))$.
    We will show that
    \begin{gather*}
    \tag{*}
        \CA \Proves \neg (\MOD{i}{f(i)} \Aand \MOD{j}{f(j)})
    \end{gather*}
    from which it follows that \CAf is inconsistent.
    \par
    Let \p and \q be such that
    \begin{gather*}
    0 \leq p, q < k \emath{,}\\
    f(i) \Cmod{p}{k} \emath{, and}\\
    f(j) \Cmod{q}{k}
    \end{gather*}
    By lemma 3.3.10, we have
    \begin{gather*}
    \tag{1}
        \CA \Proves \MOD{i}{f(i)} \then \MOD{k}{p}
        \emath{, and}\\
    \tag{2}
        \CA \Proves \MOD{j}{f(j)} \then \MOD{k}{q}
    \end{gather*}
    since $k \divides i$ and $k \divides j$.
    Since $k \ndivides f(i) - f(j)$, $p \neq q$.
    So lemma 3.3.11 yields
    \begin{gather*}
    \tag{3}
        \CA \Proves \MOD{k}{p} \then \neg \MOD{k}{q}
    \end{gather*}
    From (1), (2), and (3) we may conclude (*).
\OnlyWay
    Follows from 3.3.9a since $\CAf \subseteq \CSf$.
\item
    See proof of 3.3.9b.
\item
    See proof of 3.3.9c.
\endEnum
\endProof

\Corollary{3.3.13}
    \CAIf is consistent \iff \CSIf is consistent.
\endCorollary
\Proof
Immediate from 3.3.9c and 3.3.12c.
\endProof
It might help to review our strategy before presenting the difficult
parts of the proof that $\CA \equiv \CS$.
The main objective is (1), which follows from (2) by 3.3.1b.

\numbSentence{1}
$\CA \Proves \CS$
\endNumbSentence

\numbSentence{2}
Every completion of \CA is consistent with \CS.
\endNumbSentence

We already know that the finite completions of \CA are consistent
with \CS (see 3.3.5c) and that if \CAIf is consistent, then
\CSIf is also consistent (see 3.3.12).
So (2) is a consequence of (3).
\numbSentence{3}
If $T$ is a completion of \CAI, then $T \equiv \CAIf$
    for some total, congruous \f.
\endNumbSentence

To establish (3), it is sufficient to prove (4) because every
completion of \CAI entails \CAIf for some total \f.
\numbSentence{4}
If \f is total and congruous, then \CAIf is complete.
\endNumbSentence

To prove (4), we invoke the prime model test:
If $T$ is model complete and $T$ has a prime model, then $T$ is
complete (see Appendix, Fact C.3).
So (4) follows from (5) and (6).
\numbSentence{5}
If \f is total and congruous, then \CAIf has a prime model.
\endNumbSentence

\numbSentence{6}
For any \f, \CAIf is model complete.
\endNumbSentence

Finally, since any extension of a model complete theory is also
model complete (see Appendix, Fact C7a), we can infer (6) from (7).
\numbSentence{7}
\CAI is model complete.
\endNumbSentence
So (1) follows from (5) and (7).
\par
The proof outlined here will be carried out in Chapter 5.
But first, we consider a simpler theory, \PSIZE, which deals
only with sizes of sets and ignores boolean relations.
Chapter 4 formulates \PSIZE and establishes that it is model
complete, a result we need for showing that \CAI is model complete.

\Chapter{4}{The Pure Theory of Class Sizes}
\CS is about sets; it makes claims about sets in terms
of their boolean relations and their size relations.
In this chapter, we identify a theory, \PCS (the pure
theory of class sizes), which is not about sets, but only
about sizes of sets. \emph{Sizes}, here, are equivalence
classes of sets which have the same size.
\PCS is formulated in the language of size relations, \LS.
\par
\PCS is worth examining in its own right, for if
number theory is  the theory of cardinal
numbers, then \PCS is our version of number theory.
But our main reason for introducing \PCS is to aid the
proof that \CA is model complete.
For this reason,we only give  a sketchy treatment of
\PCS itself.
\par
Section 1 defines \PCS and develops a set of axioms, \PCA,
for it, as follows:
for each model, \ModelA, of \BASIC, the \introSTerm{size model}, \ModelSA,
consists of
equivalence classes drawn from \DomainA under the \notion{same size} relation;
\PCS is the set of statements true in \ModelSA for any
standard finite model, \ModelA.
\PCA consists of a theory, \PSIZE, which holds in \ModelSA
whenever \ModelA \Satisfies BASIC and a set of divisibility
principles.
\par
Using some results about model theory in section 2, Section
3 establishes that \PCA is model complete, the main result
of this chapter and the only result needed for subsequent proofs.
This is done by reducing \PCA to the theory \Zgm, whose models
are \fZ-groups taken modulo some specific element.
\par
Finally, section 4 indicates how \PCA could be shown to axiomatize
\PCS.
This method is the same outlined in chapter 3 to show that
$\CA \equiv \CS$.

\Section{4.1}{Size models and PCS}
\Def{4.1.1} Suppose \ModelA \Satisfies \BASIC.
\Enum
\item   If \x is a member of \DomainA, then $\SIGMA{x}{\mA}$
        is the \introSTerm{size} of \x in \ModelA:
        \[
            \SIGMA{x}{\mA} =
                \theSetst{y}{\ModelA \Satisfies (x \eqsize y)}
        \]
\item   \ModelSA, the \introSTerm{size model} for \ModelA, is the interpretation
        of \LS whose domain is
        $\theSetst{\SIGMA{x}{\mA}}{x \in \DomainA}$ where
        \begin{align*}
        \ModelSA \Satisfies
            \SIGMA{x}{\mA} \eqsize  \SIGMA{y}{\mA}
            &\tiff
            \ModelA \Satisfies x \eqsize y\\
       \ModelSA \Satisfies
            \SIGMA{x}{\mA} < \SIGMA{y}{\mA}
            &\tiff
            \ModelA \Satisfies x < y\\
       \ModelSA \Satisfies
            \tSUM(\SIGMA{x}{\mA},  \SIGMA{y}{\mA},\SIGMA{z}{\mA})
            &\tiff
            \ModelA \Satisfies \tSUM(x, y, z)\text{, and}\\
       \ModelSA \Satisfies
            \tUNIT(\SIGMA{x}{\mA})
            &\tiff
            \ModelA \Satisfies \tATOM(x)
        \end{align*}
        \item  The one-place operator  $\Comp{x}$ is to be read
        as \readAs{the complementary size of $x$}:
        \[
            \ModelSA \Satisfies (y = \Comp{x}) \tiff \ModelSA
            \Satisfies \tSUM(x,y, \Identity)
        \]
\endEnum
\endDef
These interpretations are well defined because the
predicates are satisfied by elements of \DomainA in virtue of
their sizes.
For example, if  $\ModelA \Satisfies \tSUM(x, y, z)$ and
$\ModelA \Satisfies z \eqsize \zP$ , then $\ModelA \Satisfies \tSUM(x, y, \zP)$.

\Def{4.1.2}
    \introTheory{\PCS}{pure theory of class sizes} consists
    of all sentences of \LS which are true in the size model
    of every  standard finite interpretation of \LCS.
\endDef

\Def{4.1.3}
    \introTheory{\PSIZE}{theory of sizes} consists of the following axioms:
\begin{gather*}
\intertext{Order axioms}
    \neg(x < x)         \tag{\IRREF}    \\
    x < y \Aand y < z \then x < z    \tag{\TRANSLT}\\
    x \eqsize y \bicond x = y       \tag{\UNIQEQ}\\
    x \leq \Identity                \tag{\MAX} \\
    0 \leq x                        \tag{\MIN} \\
    x < y \Or x \eqsize y \Or y < x \tag{\TRICH} \\
\intertext{Unit axioms}
    \tUNIT(x) \bicond (y < x \bicond y = 0)\\
    \exists x \tUNIT(x)\\
\intertext{Sum axioms}
    \tSUM(x, 0, x)                    \tag{\IDENT}  \\
    \tSUM(x,y,z) \bicond \tSUM(y,x,z)   \tag{\COMM}  \\
    \tSUM(x_1,y,z_1) \Aand \tSUM(x_2,y.z_2) \then   \quad\quad\quad\quad  \quad\quad\quad\quad  \tag{\MONOT}  \\
  \quad\quad\quad\quad    \quad\quad\quad\quad    (x_1 < x_2 \bicond z_1 < z_2)\\
   \tSUM(x,y,w_1) \Aand \tSUM(w_1,z,w)
            \Aand \tSUM(y,z,w_2) \quad\quad\quad\quad  \quad\tag{\ASSOC}\\
    \quad\quad\quad\quad \quad\quad\quad    \quad\quad\quad\quad \then \tSUM(x, w_2,w)  \\
    \exists z\tSUM(x,y_1,z) \Aand
                    y_2 \leq y_1 \then
                    \exists z \tSUM(x, y_2,z)
                    \tag{\EXISTPL}  \\
    x \leq z \then \exists y \tSUM(x,y,z)   \tag{\EXISTM}   \\
    \tSUM(x, \Comp{x}, \Identity)  \tag{\COMP}
\end{gather*}
\endDef
\Fact{4.1.4}
    If \ModelA \Satisfies \BASIC, then \ModelSA \Satisfies \PSIZE
\endFact
\Fact{4.1.5}
\PSIZE \nProves \PCS.
\endFact
\Proof
\PSIZE fails to axiomatize \PCS for the same reason that \BASIC
fails to axiomatize \CS: the lack of divisibility principles.
\endProof
We offer \PCA as an axiomatic version of \PCS:
\Def{4.1.6} \introTheory{\PCA}{pure class size axioms} is defined as:
    \[
     \PSIZE \union \theSetst{\ADIVn}{n > 0}
    \]
\endDef

\Fact{4.1.7}
    If $\ModelA \Satisfies \CA$, then $\ModelSA \Satisfies \PCA$
\endFact
If \ModelA is a standard finite interpretation of \BASIC with
\n atoms, then the elements of \ModelSA can be regarded as
the sequence 0, \dots, \n
with the usual ordering, where
\[
\ModelSA \Satisfies \tUNIT(x) \tiff  x = 1
\]
and
\[
\ModelSA \Satisfies \tSUM(i,j,k)  \tiff (i + j) = k \leq n
\]
Once \n is fixed, this is the only interpretation allowed by
the axioms \PSIZE.
In particular, \MONOT rules out
the interpretation in which $\tSUM(i,j,k)$ is satisfied just in
case $(i + j) \Cmod{k}{n+1}$.

\Section{4.2}{Some models and theories}
We shall use Theorem 4.2.1a to show that \PCA is model complete and
Theorem 4.2.1b to show, in chapter 5, that
\CAI is model complete.
\Theorem{4.2.1}
\Enum
\item
    If $\T$ satisfies \introSTerm{Monk's Condition}, then \T is model complete.
    \begin{quote}
    \textbf{Monk's Condition.}
    If $\ModelA \Satisfies \T$, $\ModelB \Satisfies \T$,
    $\ModelA \submodel \ModelB$, and \ModelC is a finitely
    generated submodel of \ModelB,
    then there is an isomorphism, $\f \colon \ModelC \to \ModelA$
    such that if $x \in \DomainC \intersect \DomainA$, then
    $f(x) = x$.   $f$ is a \introSTerm{Monk mapping}.
    \end{quote}
\item
    If $\T$ is model complete and \LT has no function symbols, then
    \T satisfies Monk's Condition.
\endEnum
\endTheorem
\Proof
\Enum
\item
    \cite[p.359]{Monk}
\item
    If \LT has no function symbols, then any finitely
    generated structure over \LT is finite.
    So, suppose \dC contains \listxx{a}{n} (from
    \dA) and \listxx{b}{m} (from \dB-\dA).
    Let $\phi_1$ be the diagram of \mC and obtain
    $\phi_2$ from $\phi_1$  by substituting the variable
    $x_i$ for each constant $a_i$ and the variable
    $y_i$ for each constant $b_i$.
    Finally, let $\phi_3$ be
    \[
        \exists y_1 \dots \exists y_m \phi_2
    \]
    So, $\phi_3$ is a primitive formula.
    $\mB\Satisfies \phi_3(\listAn)$, so \mA does as
    well, by Fact C.5d in the Appendix.
    So, to obtain the desired isomorphism, map the
    $a_i$'s into themselves and map the $b_i$'s into
    a sequence of elements of \dA which can stand
    in for the existentially quantified variables
    of $\phi_3$.
\endEnum
\endProof
In chapter 5, we use Monk's Theorem to infer the model completeness
of the theory \CA from that of \PCA.
To establish the model completeness of \PCA , we use Fact 4.2.2.
\Fact{4.2.2}
\Enum
\item
    If \To is model complete and \To\Proves\Tt, then \Tt
    is also model complete.
\item
    If \T is model complete in \L, and \Lpr is an expansion
    of \L by adjoining new individual constants, then \T
    is model complete in \Lpr. \cite[p.355]{Monk}
\endEnum
\endFact

\Def{4.2.3}
    Suppose \Lo and \Lt are first order languages and
    $\Lot = \Lo - \Lt$.
    A \introSTerm{translation} (\introSTerm{simple translation}) of \Lo into \Lt is
    a function, \gtau, which:
\Enum
\item
    assigns to the universal quantifier a (quantifier free)
    formula, \tauA, of \Lt with exactly one free variable.
\item
    assigns to each \n-place predicate, \P, in \Lot a
    (quantifier free) formula, \tauP, of \Lt with exactly
    \n free variables, and
\item
    assigns to each \n-place function symbol, $O$, in \Lot
    a (quantifier free) formula, \tauO, of \Lt with
    exactly $(n+1)$ free variables.
\endEnum
\endDef

\Def{4.2.4}
    If \gtau is a translation of \Lo into \Lt, then \gtau
    extends to all formulae of \Lo as follows:
\Enum
\item
    Predicates and function symbols of \Lt are translated
    into themselves.
\item
\begin{align*}
    \tau(\phi \Or \psi)& = \tau(\phi) \Or \tau(\psi)\\
    \tau(\phi \Aand \psi)& = \tau(\phi) \Aand \tau(\psi)\\
    \tau(\neg \phi) & =  \neg \tau(\phi)\\
    \tau(\forall x \phi) & =  \forall x (\tauA(x) \then \phi)\\
    \tau(\exists x \phi) & =  \exists x (\tauA(x) \Aand \phi)\\
\end{align*}
\endEnum
\endDef

\Def{4.2.5}
    If \gtau is a translation of \Lo into \Lt, then
\Enum
\item
    The \introSTerm{functional assumptions} of \gtau are the
    sentences:
    \begin{multline*}
        \forall x_1 \dots \forall x_n(
        \tauA(x_1) \Aand \dots \Aand \tauA(x_n)) \then   \\
        \exists y_1 (
            \tauA(y_1) \Aand \forall y_2(\tauO(x_1 ,\dots, x_n, y_2)
                \bicond y_1 = y_2)
        )
    \end{multline*}
    where $O$ is a function symbol in \Lo but not in \Lt.
\item
    The \introSTerm{existential assumption} of \gtau is
    \[
        \exists x \tauA(x)
    \]
\endEnum
\endDef
The functional assumptions of a translation say that the formulas
which translate function symbols yield unique values within
the relevant part of the domain when given values in the relevant
part of the domain.
The \squote{relevant part of the domain} is the set of elements
which satisfy the interpretation of universal quantifier.
The existential assumption of a translation says that that
subdomain is non-empty.
Notice that the existential and functional assumptions of a
translation are sentences of \Lt.
\par
A translation from \Lt into \Lo
induces a mapping from interpretations of
\Lt into interpretations of \Lo.

\Def{4.2.6}
If \gtau is a translation from \Lo into \Lt and \mB is an
interpretation of \Lt which satisfies the existential and
functional assumptions of \gtau, then $\tau(\mB)$ is the
interpretation \mA of \Lo such that:
\Enum
\item
    \DomainA is the set of elements of \DomainB which satisfy
    \tauA.
\item
    \mA interprets all predicates and function symbols common
    to \Lo and \Lt in the same way that \mB does.
\item
    \mA interprets all predicates and function symbols in
    \Lot in accordance with the translations assigned by \gtau.
    \begin{align*}
    \mA \Satisfies P(\xbar) & \tiff  \mB \Satisfies \tauP(\xbar)\\
    \mA \Satisfies y = O(\xbar) & \tiff \mB \Satisfies \tau_O(\xbar,y)
    \end{align*}

\endEnum
\endDef

The following condition on theories
allows us to infer the model completeness of one from
the model completeness of the other.
\Def{4.2.7}

\Enum
\item
    If \gtau is a translation from \Lo to \Lt, then
    \To is \introMTerm{\gtau-reducible}{tau-reducible} to \Tt \iff for every
    model \mA of \To there is a model \mB of \Tt such that
    $\mA = \tau(\mB)$.
\item
    \To is \introSTerm{reducible} (\introSTerm{simply reducible})   to \Tt \iff there is
    a (simple) translation, \gtau, for which \To is \gtau-reducible
    to \Tt.
\item
    \To is  \introSTerm{uniformly \gtau-reducible} to \Tt \iff
    for any models, \mAo and \mBo, such that $\mAo \Satisfies \To$
    and $\mBo \Satisfies \To$, and  $\mAo \submodel \mBo$,
    there exist models \mAt and \mBt such that
    $\mAt \Satisfies \Tt$
    and $\mBt \Satisfies \Tt$, and  $\mAt \submodel \mBt$
    $\mAo = \tau(\mAt)$ and $\mBo = \tau(\mBt)$,
\endEnum
\endDef

\Lemma{4.2.8}
    Suppose that \To is \gtau-reducible to \Tt and that
    $\mAo = \tau(\mAt)$.
    Then, for any primitive formula, \gphi, of \Lo and any
    sequence, $\xbar \in \A_1$,
    \[
        \mAo \Satisfies \phi(\xbar) \tiff
        \mAt \Satisfies \tauph(\xbar)
    \]

\endLemma
\Proof
Suppose
    \[
        \phi(\xbar) = \exists y_1\dots \exists y_n
                \phiP(\xbar)
    \]
    where $\phiP(\xbar)$ is a conjunction of atomic
    formulae and negations of atomic formulae.  Then
    \begin{align*}
    \mAo \Satisfies \phi(\xbar)
    &\tiff    \mAo \Satisfies
                \Exists{y}{n}\phiP(\xbar,\listxx{y}{n})\\
    &\tiff     \mAo \Satisfies
                \phiP(\xbar,\listBn) \etext{, for } b_i \in \dA  \\
     &\tiff     \mAt \Satisfies
                \tauphp(\xbar,\listBn)   \\
     &\tiff     \mAt \Satisfies
                \Exists{y}{n} \\
       &\quad\quad\quad(
                \tauA(y_1) \Aand \dots \Aand \tauA(y_n)
                \Aand \phiP(\xbar,\listxx{y}{n})
                )\\
     &\tiff     \mAt \Satisfies
                \tauph(\xbar)
    \end{align*}
\endProof
\Theorem{4.2.9}
    If \Tt is model complete, \gtau is a simple translation
    from \theLang{\To} to \theLang{\Tt}, and \To is uniformly \gtau-reducible
    to \Tt, then \To is also model complete.
\endTheorem
\Proof
By Fact C.5d in the Appendix, it is enough to show that given models \mAo and \mBo of \To, where
$\mAo \submodel \mBo$ and a primitive formula, $\phi$:
\begin{flalign*}
\If     \mBo \Satisfies \phi(\xbar) \emath{for} \xbar \in \Domain{A} \\
\Then   \mAo \Satisfies \phi(\xbar)
\end{flalign*}
Since   \To is uniformly \gtau-reducible to \Tt, there are models of \Tt,
$\mAt \submodel \mBt$, where $\mAo = \tau(\mAt)$ and $\mBo = \tau(\mBt)$
\begin{flalign*}
\Since  \mBo \Satisfies  \phi(\xbar)        \by{assumption}\\
\Then   \mBt \Satisfies \tauph(\xbar)    \by{lemma 4.2.8}\\
\So     \mAt \Satisfies \tauph(\xbar)    \since{\Tt is model complete}\\
\AND    \mAo \Satisfies \phi(\xbar)            \by{lemma 4.2.8}
\end{flalign*}
\endProof
We shall now define several theories, all more or less familiar, which will
serve in showing that our theory of size is model complete.
\Def{4.2.10}
\Enum
\item
    The \introTheoryof{abelian groups with identity} has the following
    axioms
    \begin{gather*}
    x + (y + z) =(x+y) + z \tag{1} \\
    x + y  = y + x          \tag{2} \\
    x + 0 = x               \tag{3} \\
    \exists y (x + y = 0)   \tag{4}
    \end{gather*}
\item
    The \introTheoryof{cancellable abelian semigroups with identity}
    consists of (1), (2), and (3) above and:
    \begin{gather*}
    x+y=x+z \then y =z      \tag{$4^\prime$}
    \end{gather*}
\item
    The axioms of \introSTerm{simple order} are:
    \begin{gather*}
    x \leq y \Aand \y \leq z \then x \leq z \tag{5}  \\
    x \leq y \Aand \y \leq x \then x = y    \tag{6}  \\
    x \leq x                                \tag{7}  \\
    x \leq y \Or \y \leq x                  \tag{8}
    \end{gather*}
\item
    The \introMTheoryof{\protect \fZ-groups, \Zg}{Z-groups}, has the following axioms
    \Enum
    \item
        The axioms for abelian groups with identity,
    \item
        The axioms for simple order,
    \item
        The following additional axioms:
        \begin{gather*}
        y \leq z  \then x+y \leq x + z  \tag{9} \\
        1 \emath{is the least element greater than} 0      \tag{10} \\
        \forall x \exists y (ny = x \Or \dots \Or ny = x+(n-1))  \tag{11}
        \end{gather*}
        for each positive \n,  where $ny$ stands for:
        \[
        \underbrace{y + \dots + y}_{ n \emath{times}}
        \]
    \endEnum
\item
     The \introMTheoryof{\protect \fN-semigroups}{N-semigroups} has the following axioms:
     \Enum
     \item
     The axioms for cancellable abelian semigroups with identity,
     \item
     The axioms of simple order,
     \item
     Axioms (9) and (10) of \Zg, and
     \item
     The additional axiom:
     \begin{gather*}
     0 \leq x \tag{12} \\
     \end{gather*}
     \endEnum
\item
    The \introMTheoryof{\protect \fZ-groups modulo $I$}{Z-groups modulo $I$}, \Zgm, consists of the following axioms:
     \Enum
     \item
     The axioms for abelian groups with identity,
     \item
     The axioms for simple order
     \item
     Axioms (10) and (11) of \Zg, axiom (12) from the theory of \fN-semigroups, and
     \item
     The following additional axioms
     \begin{gather*}
     (y  \leq z \Aand x \leq x + z) \then x + y \leq x + z)   \tag{13} \\
     x \leq I \tag{14}
     \end{gather*}
     \endEnum
\endEnum
\endDef

The theory of \fZ-groups is taken from \cite[p. 291]{Chang}

\Fact{4.2.11}
\Enum
\item
    \Zg is the complete theory of \tuple{\fZ, +, 0, 1, \leq}.
    \cite[p.291]{Chang}
\item
    \Zg is model complete.\cite{Robinson}
\endEnum
\endFact

\Theorem{4.2.12}
\Enum
\item
    The theory of \fN-semigroups is model complete.
\item
    The theory of \fZ-groups modulo $I$ is model complete.
\endEnum
\endTheorem
\Proof
\Enum
    \item
    Every abelian semigroup with cancellation can be isomorphically embedded
    in an abelian group \cite[pp. 44-48]{Kurosh}.
    It is clear from the construction in \cite{Kurosh} that if the semigroup is
    ordered, the abelian group in which it is embedded may also be ordered
    and that
    the elements of the semigroup will be the positive elements of the group.
    Moreover, the (rough) divisibility of the elements in the semigroup will
    be carried over to the group.
    \par
    Consequently, the theory of \fN-semigroups is uniformly reducible to
    the theory of \fZ groups by the translation:
    \[
        \tauA = \schema{0 \leq x}
    \]
    Since the latter is model complete, so is the former, by 4.2.9.
    \item
    First, consider the theory of \fN-semigroups in the language which
    contains, besides the constant symbols in the original theory, an individual
    constant, $I$.
    The theory of \fN-semigroups is model-complete in this language, by 4.2.2b.
    \par
    We claim that the theory of \fZ-groups modulo $I$ is uniformly reducible
    to this new theory by the following translation:
    \begin{gather*}
        \tauA = \schema{x \leq I}        \\
        \tau_+ = \schema{x + y = z \Or x + y = I + z}
    \end{gather*}
    (The construction: given a model of \Zgm, stack up $\omega$ many copies
    of the model, assigning interpretations in the obvious way.
    The result is an \fN-semigroup and the original
    model is isomorphic to the first copy of itself.)
\endEnum
\endProof
In the next section, we use Theorem 4.2.9 to show that
\PSIZE is model complete, by reducing it to the theory of
\fZ-groups modulo $I$.
Chapter 5 uses Monk's Theorem to show that \CA is model
complete.
\Section{4.3}{PCA is model complete}

To show that \PCA is model complete, we shall reduce it to \Zgm,
the theory of \fZ-groups with addition taken modulo some constant
(see 4.2.11).
The model completeness of \PCA then follows from the model
completeness of \Zgm (Fact 4.2.12) and Theorem 4.2.9.
Specifically, we shall show that every model of \PCA
is the \gtau-image of a model of \Zgm, where \gtau is the following
translation.

\Def{4.3.1}
Let \gtau be the translation from \theLang{\PSIZE} to \theLang{\Zgm} where:
\begin{align*}
\tauA               &= \schema{x = x}       \\
\tau_{\tSUM(x,y,z)}   &= \schema{x+y = z \Aand x \leq z}    \\
\tau_{\tUNIT(x)}      &=  \schema{x = 1}    \\
\tau_{\Comp{x}}     &=   \schema{x + y = \Identity}   \\
\tau_{x < y}        &=   \schema{x \leq y \Aand x \neq y}   \\
\tau_{\NullSet}     &= \schema{x = 0}
\end{align*}
\endDef
Given a model, \mA, of \PCA, we can construct a model, \ZgmA, of
\Zgm directly:

\Def{4.3.2}
If \mA \Satisfies \PSIZE, then $\ZgmA$ is the interpretation
of \theLang{\Zgm} in which:
\begin{align*}
\Domain{\ZgmA} &= \dA   \tag{a}   \\
\ZgmA \Satisfies (x \leq y) &\tiff \mA \Satisfies (x <y \Or x = y) \tag{b}   \\
\ZgmA(I) &= \mA(I)  \tag{c}\\
\ZgmA(0) &= \mA(0)  \tag{d}\\
\ZgmA \Satisfies (x + y = z) &\tiff \mA \Satisfies \tSUM(x,y,z) \tag{e} \\
    & \quad \emath{or}
    \mA \Satisfies \exists a \exists w \tUNIT(a)\\
     & \quad \quad \quad \Aand \tSUM(\Comp{x},\Comp{y},w) \Aand \tSUM(w,a,\Comp{z}) \notag\\
\ZgmA \Satisfies (x =1) &\tiff \mA \Satisfies \tUNIT(x) \tag{f}
\end{align*}
\endDef
Fact 4.3.3j establishes that 4.3.2 gives \tSUM a functional interpretation.
Theorem 4.3.8 establishes that $\ZgmA\Satisfies\Zgm$, on the basis
of the intervening facts:
4.3.3 deals with the model \mA of \PSIZE; 4.3.4 deals with
the corresponding model $\ZgmA$;
4.3.6 verifies some connections between \mA and $\ZgmA$.
\Fact{4.3.3}
The following are theorems of \PSIZE. ($\tSm(x,y)$ abbreviates
$\exists z \tSUM(x,y,z)$. )

\begin{gather*}
\tSUM(x,y,z_1) \Aand \tSUM(x,y,z_2) \then z_1 = z_2 \quad\quad\etext{(\UNIQUEPL)} \tag{a}  \\
\tSUM(x,y_1,z) \Aand \tSUM(x,y_2,z) \then y_1 = y_2 \quad\quad \etext{(\UNIQUEMN)} \tag{b}  \\
 \tSUM(x,y,w_1) \And \tSUM(w_1, z,w) \then
    \quad\quad\quad \quad\quad\quad \quad\quad \etext{(\ASSOCT)} \tag{c}  \\
    \quad\quad\quad
    \exists w_2 (\tSUM(y,z,w_2) \Aand \tSUM(x,w_2,w) )  \\
 \Comp{\Comp{x}} = x  \tag{d}  \\
x < y \bicond \Comp{y} < \Comp{x}    \tag{e}   \\
x < y \Aand y < \Comp{y} \then x < \Comp{x}   \tag{f}   \\
\tSm(x,y) \Or \tSm(\Comp{x}, \Comp{y})  \tag{g} \\
\tSUM(x,y,z) \then \tSUM(\Comp{z},y,\Comp{x})  \tag{h1} \\
\tSUM(x,y,z) \then \tSUM(x, \Comp{z},\Comp{y})  \tag{h2} \\
\tSUM(x,y,z_1) \Aand \tSUM(\Comp{x}, \Comp{y}, \Comp{z_2}) \then
    (z_1 = z_2 = I) \tag{i}  \\
 \exists z_1 \tSUM(x,y,z_1) \bicond \quad\quad\quad \quad\quad\quad
 \quad\quad\quad \quad\quad\quad \quad\quad\quad \quad\quad\quad \tag{j} \\
 \quad\quad\quad
\neg \exists a \exists w \exists z_2
    (\tUNIT(a)
        \Aand \tSUM(\Comp{x}, \Comp{y}, w)
        \Aand  \tSUM(w,a, \Comp{z_2})
    )
\end{gather*}

\endFact
\Proof
The proofs are elementary.
\endProof

In addition to the theorems of \PSIZE listed in 4.3.3, we require a battery of
tedious facts about the model \ZgmA.
We shall state these in terms of the model \mZgastar, an expansion of both \mA
and \ZgmA, which interprets two additional operators, as follows:

\Def{4.3.4}
    Given a model \mA of \PSIZE,
    $\mZgastar$  is the expansion of \ZgmA induced by the definitions
    in 4.3.2 together with (a) and (b):
    \begin{align*}
        \mZgastar \Satisfies y = -x  &\tiff \ZgmA\Satisfies x+y=0    \tag{a} \\
        \mZgastar \Satisfies z = x -y &\tiff \ZgmA\Satisfies z = x+-y \tag{b}
    \end{align*}
\endDef

\Remark
    Given a model \mA of \PSIZE: for each $x$ in
    \Domain{\ZgmA}, there is
    a unique $y$ such that:
    \[
        \ZgmA \Satisfies x + y = 0
    \]
\endRemark
\Proof
{\allowdisplaybreaks
\Enum
\item   If $x = 0$, then $\ZgmA \Satisfies x + y = 0$ iff $y=0$.
\BothWays
\IfWayB
    \begin{flalign*}
    \If             \ZgmA \Satisfies 0 + y = 0  \\
   \Then           \mA \Satisfies \tSUM(\NullSymbol, y, \NullSymbol)
                            \by{4.3.2e} \\
   \But             \mA \Satisfies \tSUM(\NullSymbol, \NullSymbol, \NullSymbol)\by{IDENT}  \\
   \So             \ZgmA \Satisfies  y = 0 \by{\UNIQUEMN}
   \end{flalign*}
\OnlyWayB    Obvious
\enBoth
\item    If $x > 0$, there is a unique $y$ which satisfies
\begin{align*}
        \tSUM(\Comp{x},a,y) \emath{, where \tUNIT($a$)} \tag{*}
\end{align*}
since there is a unique unit and $\mA \Satisfies \UNIQUEPL$.
\par
But $\ZgmA \Satisfies x + y = 0$ iff $\mA \Satisfies $ (*).
\BothWays
\IfWayB
   \begin{flalign*}
   \If      \ZgmA \Satisfies x + y = 0 \\
   \Then    \mA \Satisfies \tSUM(\Comp{x}, \Comp{y}, w)   \\
   \AND     \mA \Satisfies \tSUM(w, a, \Comp{\NullSymbol})
                \since{$\mA \nSatisfies \tSUM(x,y, 0)$}\\
   \But     \mA  \Satisfies \Comp{\NullSymbol} = \Identity \\
   \So      \mA  \Satisfies w = \Comp{a} \by{\COMP and \UNIQUEMN}\\
   \So      \mA  \Satisfies \tSUM(\Comp{x}, \Comp{y}, \Comp{a})   \\
   \So      \mA  \Satisfies  \tSUM(\Comp{x}, a, y)  \by{4.3.3d and h2}
   \end{flalign*}

\OnlyWayB
   \begin{flalign*}
   \If      \mA \Satisfies \emath{(*)}  \\
   \Then    \mA \Satisfies  \tSUM(\Comp{x},\Comp{y},\Comp{a})   \by{4.3.3h2}   \\
   \AND     \mA \Satisfies  \tSUM(\Comp{a},a, \Comp{\NullSymbol})   \by{COMP}   \\
   \So      \ZgmA \Satisfies  x + y = 0   \by{4.3.2e}
   \end{flalign*}

\enBoth
\endEnum
}   
\endProof
\Fact{4.3.5}
    \mZgastar satisfies the following:
\begin{gather}
    -(x+y) = -x + -y   \tag{a} \\
    -(x -y) = y -x     \tag{b} \\
    (x+y)-y =x         \tag{c} \\
    -(-x) = x          \tag{d} \\
     x\neq0 \Aand x < y \then -y < -x   \tag{e}
\end{gather}
\endFact
\Proof
Omitted.
\endProof

\Fact{4.3.6}
    \mZgastar satisfies the following:
\Enum
\item
    $\tSm(a,b) \bicond \tSUM(a,b,a+b)$
\item
    $c \leq b \bicond \tSUM(c, b-c, b)$
\item
    $\tSUM(a,b,c) \Aand 0 < b \then \tSUM(a,-c,-b)$
    \par
    $\tSUM(a,b,c) \Aand 0 < a \then \tSUM(-c,b,-a)$

\item
    $\neg \tSm(a,b) \then ( \tSm(-a,-b) \Or b = -a )$
\endEnum
\endFact
{\allowdisplaybreaks
\Proof
\Enum
\item
\IfWay
 \begin{flalign*}
    \Suppose \tSm(a,b)\\
    \so \tSUM(a,b,w) \emath{for some $w$}\\
    \so (a+b) = w \by{4.3.2e}\\
    \so \tSUM(a,b,a+b)
 \end{flalign*}
\OnlyWay
    Obvious.

\item
\IfWay
 \begin{flalign*}
    \Suppose    c \leq b\\
    \Then       \tSUM(c,w,b) \ForSome{w} \by{\EXISTM}\\
    \so         b = c+w \by{(a) and \UNIQUEPL}\\
    \so         b-c = (c+w)-c\\
    \so         b-c = w \by{4.3.5c}\\
    \so         \tSUM(c,b-c,b)
 \end{flalign*}

\OnlyWay
\begin{flalign*}
    \Suppose    \tSUM(c,b-c,b)\\
    \but        \tSUM(c,0,c) \by{\IDENT}\\
    \so         0 \leq (b-c) \bicond b \leq c
                    \by{\MONOT and \UNIQUEPL}\\
    \so         b \leq c \by{MIN}
 \end{flalign*}
\item
\begin{flalign*}
    \Suppose    \tSUM(a,b,c)\\
    \AND        0 < b\\
    \Then       c = (a+b) \by{(a) and \UNIQUEPL}\\
    \so         c-b = (a+b) -b\\
    \so         c-b = a \by{4.3.5c}\\
    \so         -b = a -c \by{4.3.5c}\\
    \so         \tSUM(a,-c,-b) \emath{if} \tSm(a,-c)
                                        \by{(a)}\\
    \Butif      a =0\\
    \Then       \tSm(a,-c) \by{\IDENT}\\
    \Andif      a \neq 0\\
    \Then       a < c\emath{, since}0<b
                    \by{\MONOT and \UNIQUEPL}\\
    \so         -c < -a \by{4.3.5e}\\
    \so         -c \leq \Comp{a} \by{4.3.4}\\
    \but        \tSm(a,\Comp{a}) \by{\IDENT}\\
    \so         \tSm(a,-c) \by{\EXISTPL}
 \end{flalign*}

\item
\begin{flalign*}
\Suppose    \neg \tSm(a,b)\\
\but        b = (a+b) -a,               \by{4.3.5c}\\
\so         \neg \tSUM(a,(a+b) -a, a+b),  \by{(a)}\\
\so         \neg(a \le a+b),            \by{(b)}\\
\so         (a+b) < a,                  \by{\TRICH}\\
\Butif      (a+b) \ne 0\\
\Then       -a < -(a+b),                \by{4.3.5e}\\
\so         -a < -a + -b,               \by{4.3.5a}\\
\so         \tSUM(-a,(-a + -b) - (-a), -a + -b), \by{(b)}\\
\but        (-a + -b) - (-a) = -b,      \by{4.3.5c}\\
\so         \tSUM(-a,-b,-a+-b)\\
\so         \tSm(-a,-b)\\
\Andif      (a+b) = 0\\
\Then       b = -a\\
\end{flalign*}
\endEnum
\endProof }

We can now show that $\ZgmA\Satisfies\Zgm$.
The only real difficulty arises in verifying that
addition is associative.
We need the following lemma.
\Lemma{4.3.7}$\ZgmA\Satisfies (a+c) + (b-c) = a+b$
\endLemma

 \Proof
{\allowdisplaybreaks
We shall work through successively more general cases:
\Enum
\item  When $\tSm(a,b)$ and $c \leq b$
\begin{flalign*}
\Weknow      \tSUM(b-c,c,b) \by{4.3.6b}     \\
\AND         \tSUM(b,a,a+b) \by{4.3.6a}         \\
\So          \exists w (\tSUM(c,a,w) \\
            & \quad\quad\Aand \tSUM(b-c,w,a+b) \by{\ASSOCT}       \\
\So          w= c+a \by{4.3.6a and \UNIQUEPL}       \\
\AND         a+b = (b-c) + (c + a) \by{4.3.6a and \UNIQUEPL}       \\
\So          a+b = (a+c) + (b-c)
\end{flalign*}
\item When $\tSm(a,b)$ and $\tSm(a,c)$
\begin{flalign*}
\If         c \leq b \emath{, Case (a) applies directly}      \\
\SoAssume   b < c      \\
\Then       a+c = (a+b) + (c-b) \by{Case (a)}     \\
\So         (a+c) - (c-b) = a+b \by{4.3.5c}    \\
\So         (a+c) + -(c-b)  = a+b \by{4.3.4b}  \\
\So         (a+b) + (b-c) = a+b \by{4.3.5b}
\end{flalign*}
 
\item $\tSm(a,b)$
\begin{flalign*}
\If        \tSm(a,c) \emath{Case (b) applies}      \\
\SoAssume  \neg \tSm(a,c) \emath{and thus} b<c \by{\EXISTPL}       \\
\So        \tSUM(b,c-b,c) \by{4.3.6b}        \\
\AND       0 < c - b \by{\MONOT}     \\
\So        \tSUM(b,-c,-(c-b)) \by{4.3.6c}       \\
\so        \tSUM(b,-c,b-c) \by{4.3.5b}           \\
\so        \tSUM(-c,b,b-c)        & \quad\quad\quad\quad\quad\quad\etext{(c1)}      \\
\\
\Either    \tSm(-a,-c) \emath{or} c = -a \emath{since} \neg \tSm(a,c) \by{4.3.6d}    \\
\Assume    \tSm(-a,-c)           \\
\Then      \tSUM(-a,-c,-a + -c) \by{4.3.6a}           \\
\So        \tSUM(-a,-c,-(a+c)) \by{4.3.5a}    \\
\AND      0 < -a   \\
\So       \tSUM(-(-(a+c)), -c, -(-a)) \by{4.3.6c}   \\
\So       \tSUM(a+c, -c,a)    \by{4.3.5d} \quad\etext{(c2)}\\
\\
\Assume   c = -a   \\
\Then     a+c = 0  \\
\AND      a = -c   \\
\Hence       \tSUM(a+c, -c,a) \by{(c2)}  \\
\AND      \tSUM(-c,b,b-c) \by{(c1)}   \\
\AND    \tSUM(a,b,a+b) \since{$\tSm(a,b)$}      \\
\so     \tSUM(a+c, b-c,a+b) \by{\ASSOC}       \\
\So   (a+c) + (b-c) = a+b \by{4.3.6a}
\end{flalign*}

  \item Whenever
\begin{flalign*}
\Suppose    \neg \tSm(a,b)   \emath{, for otherwise case c applies.} \\
\Either   \tSm(-a,-b) \emath{or} b = -a \by{4.3.6d}   \\
\Assume   b = -a   \\
\Then     a+ b = 0   \\
\AND      b-c = (-a)-c \\
            &\quad\quad= -a + -c \\
            &\quad\quad= -(a+c)   \\
\So       (a+b) + (b-c) \\
            &\quad\quad= (a+c) + -(a+c)  \\
            &\quad\quad= 0 = a+b\\
\Assume    \tSm(-a,-b)   \\
\Then      -a + -b = (-a + -c) + (-b - -c)\by{case c}    \\
\So        -(a+b) = -(a+c) + -(b-c) \by{4.3.5a}   \\
\So         -(a+b) = -((a+c) + (b-c)) \by{4.3.5a}   \\
\So        a+b = (a+c) + (b-c) \by{4.3.5d}
\end{flalign*}

\endEnum
    } 
\endProof
 
\Theorem{4.3.8}
    $\ZgmA \Satisfies \Zgm$
\endTheorem
\Proof
{\allowdisplaybreaks
The only axiom for abelian groups that needs further verification is associativity.
We prove this using lemma 4.3.7:
\begin{align*}
    x + (y + z) &= (x + y) + ((y+z) -y)\emath{, by 4.3.7} \\
                &= (x + y) + z \emath{, by 4.3.5c}
\end{align*}

The ordering axioms of \Zgm are satisfied in \ZgmA because \ZgmA uses the same
ordering as \mA, \mA \Satisfies \PSIZE, and \PSIZE includes the same ordering
axioms.

 \par
\ZgmA satisfies axiom (10) of \Zgm because
\mA satisfies the \Theory{UNIT} axiom of \PSIZE.
\par
The divisibility of elements in \ZgmA required by axiom (11) of \Zgm
is guaranteed
by the fact that \mA satisfies the divisibility principles of \PCA.
\par
Similarly, \ZgmA satisfies axioms (12), and (14) of \Zgm because
\mA satisfies  \MIN and \MAX.
\par
Axiom (13) of \Zgm is:
\[
    y \le z \Aand x \le x + z \then x+y \le x+z
\]
\begin{flalign*}
\If     \ZgmA \Satisfies x \le x + z \\
\Then   \mA \Satisfies \tSUM(x,z,x+z)   \\
\Since    \ZgmA \Satisfies y < z\\
\Then   \mA \Satisfies \tSUM(x,y,x+y)  \by{ \EXISTPL}    \\
\AND     \mA \Satisfies x+y \le x+z    \by{ \MONOT  }  \\
\So      \ZgmA \Satisfies x+y \le x+z
\end{flalign*}

} 
\endProof

So, we \PCA is \gtau-reducible to \Zgm.
The reduction is  uniform since each model \mA
of \PCA has the same domain as its \Zgm-model.
So, we may conclude that \PCA is model complete.

\Theorem{4.3.9}
\Enum
\item
    \PCA is model complete.
\item
    \PCA satisfies Monk's condition.
\endEnum
\endTheorem
\Proof
\Enum
\item
    Apply Theorem 4.2.9.
\item
    Immediate from (a) and 4.2.1(b).
\endEnum
\endProof
\Section{4.4}{Remarks on showing that PCA axiomatizes PCS}

To prove that $\PCA \equiv \PCS$, we could follow the method
outlined at the end of chapter 3 for showing that $\CA \equiv \CS$.

We already know that $\PCS \Proves \PCA$, so we need only
prove (1), which follows from (2) by 3.3.1b.
\numbSentence{1}
    $\PCA \Proves \PCS$
\endNumbSentence
\numbSentence{2}
    Every completion of \PCA is consistent with \PCS.
\endNumbSentence
But (2) is equivalent to the conjunction of (2a) and (2b).

\numbSentence{2a}
    Every finite completion of \PCA is consistent with \PCS.
\endNumbSentence
\numbSentence{2b}
    Every infinite completion of \PCA is consistent with \PCS.
\endNumbSentence

The finite completions of \PCA are the
theories $\PCA;\EXn$. But $\PCA;\EXn$ is true in \ModelSA, where
\mA is the standard finite interpretation of \LCS containing
\n atoms.
So, the finite completions of \PCA are consistent with \PCS.
(Formally, we would have to define $\EXn$ in terms of units
rather than atoms.)
\par
Letting $\PCAI = \PCA + \INF$ and
$\PCSI = \PCS + \INF$,
(2b) is a consequence
of (3a) and (3b).
\index{PCA}
\numbSentence{3a}
    If \PCAIf is consistent, then \PCSIf is consistent.
\endNumbSentence
\numbSentence{3b}
    If \T is a completion of \PCAI, then $T \equiv \PCAIf$, for
    some \f.
\endNumbSentence
To prove (3a), we would have to prove analogues of 3.3.9 through
3.3.13 for \PCA and \PCS.
This seems straightforward, but tedious.
The trick is to show that enough axioms about size have been
incorporated in \PCA to establish the entailments among MOD
statements.
\par
To prove (3b), it is sufficient to demonstrate (4), because every
completion of \PCAI entails \PCAIf for some total \f.
\numbSentence{4}
    If \f is total, then \PCAIf is complete.
\endNumbSentence
But PCA is model complete, so only (5) remains to be shown.
\numbSentence{5}
    If \f is total and congruous, then \PCAIf is has a prime model.

\endNumbSentence
We will not construct prime models for the extensions of \PCAI.
It's apparent that the size models of the prime models for \CAIf
would do nicely.
Alternatively, the construction could be duplicated in this simpler
case.

\Chapter{5}{Completeness of CA}

\Section{5.1}{Model completeness of CA}
We shall show that \CA is model complete by showing
that it satisfies Monk's criterion (see 4.3.1).
So, given Assumption 5.1.1, we want to prove 5.1.2.
\Assumption{5.1.1}
    $\ModelA \Satisfies \CA$, $\ModelB \Satisfies \CA$,
    $\ModelA \submodel \ModelB$
    , and
    \ModelC is a finitely
    generated substructure of \ModelB.
\endAssumption
\Theorem{5.1.2}
    There is an isomorphic embedding, $f\colon\ModelC \to \ModelA$, where
    \[
        f(x) = x \emath{if} x \in \CIA
    \]
\endTheorem
The Monk mappings for \PCA can serve as
a guide in constructing Monk mappings for \CA.
The existence of Monk mappings for \PCA tells us that
we can find elements with the right sizes.
\DISJU and \REPLT then allow us to find elements with
those  sizes that fit together in the right way.
\par
Strictly speaking, \ModelSA is not a submodel of \ModelSB,
so we cannot apply Monk's Theorem directly.
But, let
\[
    \ModelSR{X}{\ModelB} = \emath{the submodel of} \ModelSB
    \emath{whose domain is} \theSetst{\SIGMA{x}{\mB}}{x \in X}
\]
Then, clearly,
\[
    \ModelSA \iso \ModelSR{A}{\ModelB} \submodel \ModelSB
    \emath{, and}
\]
\[
    \ModelSC \iso \ModelSR{C}{\ModelB}
    \emath{, a finitely generated submodel of \ModelSB.}
\]
\par
Monk's Theorem applies directly to $\ModelSR{A}{ \ModelB}$,
\ModelSB, and $\ModelSR{C}{\ModelB}$, so we may conclude
5.1.3:
\Fact{5.1.3}
    There is an isomorphic embedding,
    $g \colon \ModelSC \to \ModelSA$,
    where
    \[
    g(\SIGMA{x}{\mB}) = \SIGMA{x}{\mA} \etext{ for all }
    x \in \DomainC \intersect \DomainA \etext{.}
    \]
\endFact
That is to say, the sizes of elements in \ModelC can be embedded
in the sizes of elements in \ModelA.
It remains to be shown that the elements of \ModelC themselves can
be mapped into \ModelA by a function which preserves boolean
relations as well as size relations.
\par
\ModelC is a finite, and hence atomic, boolean algebra,
though its atoms need not be atoms of \ModelB;
indeed, if \ModelB is infinite, there must be some atoms of \ModelC
which are not atoms of \ModelB since the union of all atoms of \ModelC
is the basis of \ModelB.
\Def{5.1.4}
$d$ is a \introSTerm{molecule} \iff $d \in \DomainC \intersect \DomainA$ and
no proper subset of $d$ is in $\DomainC \intersect \DomainA$.
\endDef
\CIA is a boolean algebra whose basis is the same as the
basis of \ModelC.
So every atom of \ModelC is included in some molecule.
The embedding \f has to map each molecule to itself.
Moreover, \f has to be determined by its values on \ModelC since
\f must preserve unions.
In fact, the atoms of \ModelC can be partitioned among the molecules.
So, if
\[
    d = b_1 \union \dots \union b_n
\]
where \d is a molecule and $b_1, \dots , b_n$ are the atoms of \ModelC
contained in \d, then \d must also be the (disjoint) union of
$f(b_1), \dots , f(b_n)$.
If this condition is satisfied, \f will preserve boolean relations.
\f must also select images with appropriate sizes.
\Proof
Proof of 5.1.2
Given a molecule, \d, let $b_1, \dots, b_n$ be the atoms
of \C contained in \d.
For each \bi, let \ci be some member of $g(\SIGMA{\bi}{\mB})$
These elements will be elements of \dA, since \g yields
sizes in \mA whose members are in \dA.
These elements have the right sizes, as we will show
below, but they are in the wrong places.
We have no guarantee that they are contained in the
molecule \d.
So we still need to show that there are disjoint
elements of \dA, $a_1, \dots, a_n$ whose union is \d
and whose sizes are the same as those of
$b_1, \dots, b_n$, respectively.
Well,
\begin{flalign*}
    \Suppose    d = b_1 \union \dots \union b_n \padleft\\
    \so         \mC \Satisfies \tSUM(b_1,\dots,b_n,d)
                    \\
    \so \mSC \Satisfies \tSUM(\SIGMA{b_1}{\mB},
            \dots,\SIGMA{b_n}{\mB},\SIGMA{d}{\mB})\\
    \so \mSA \Satisfies \tSUM(g(\SIGMA{b_1}{\mB}),
            \dots,g(\SIGMA{b_n}{\mB}), g(\SIGMA{d}{\mB}))\\
    \so \mSA \Satisfies \tSUM(c_1,\dots,c_n,d)
\end{flalign*}
since each $c_i \in g(\SIGMA{\bi}{\mB})$ and
$d \in g(\SIGMA{d}{\mB}) = \SIGMA{d}{\mA}$.
So, the existence of $a_1, \dots, a_n$, as above,
is guaranteed because $\mA \Satisfies \DEFPL$.
\par
Now, let $f(b_i) = a_i$ for each $b_i$ in the molecule \d.
Repeating this procedure for each molecule yields a value of
\f for each atom of \ModelC.
Finally, if $x \in \DomainC$ is non-atomic, then
\[
    x = b_1 \union \dots \union b_k
\]
where each $b_k$ is atomic.
So let
\[
    f(x) = f(b_1) \union \dots \union f(b_k)
\]
\f satisfies the requirements of Theorem 5.1.2:
Boolean relations are preserved by \f because the function
is determined by its values on the atoms of \mC;
\f maps elements of $C \intersect A$ into themselves
because the set of atoms contained in each of these
molecules is mapped into a disjoint collection of
elements of \dA whose union is the same molecule;
so, we need only show that \f preserves size relations.
To do so, we invoke lemma 5.1.6, below:
\begin{align*}
\tag{a}
    \mC \Satisfies x < y &\tiff \mA \Satisfies f(x) < f(y)\\
\tag{b}
    \mC \Satisfies x \eqsize y &\tiff \mA \Satisfies f(x) \eqsize f(y)\\
\tag{c}
    \mC \Satisfies \tSUM(x,y,z) &\tiff \mA \Satisfies \tSUM(f(x),f(y),f(z))\\
\tag{d}
    \mC \Satisfies \tUNIT(x) &\tiff \mA \Satisfies \tUNIT(y)
\end{align*}
\begin{align*}
\intertext{Proof of (a):  (The others are similar).}\\
\mC \Satisfies x < y
&\tiff  \mSB \Satisfies \SIGMA{x}{\mB} < \SIGMA{y}{\mB}\\
&\tiff  \mSA \Satisfies
                g(\SIGMA{x}{\mB}) < g(\SIGMA{y}{\mB})\\
&\tiff  \mSA \Satisfies
                \SIGMA{f(x)}{\mA} < \SIGMA{f(y)}{\mA}
                \by{5.1.6}\\
&\tiff  \mA \Satisfies f(x) < f(y)
\end{align*}
\endProof
So, we may conclude:
\Theorem{5.1.6}
    \CA is model complete.
\endTheorem

\Lemma{5.1.5}
    For all \x in \Domain{C},
    \[
        \SIGMA{f(x)}{\mA} = g(\SIGMA{x}{\mB})
    \]
\endLemma
\Proof
    Suppose
    \[
        x = b_1 \cup \dots \cup b_n
    \]
    where each $b_i$ is an atom of \dC.
    \begin{flalign*}
    \So \mB \Satisfies \tSUM(b_1, \dots, b_n, x)
        \since {\mB \Satisfies \DISJU}\\
    \so \mSB \Satisfies
        \tSUM(\SIGMA{b_1}{\mB},\dots,\SIGMA{b_n}{\mB},
                \SIGMA{x}{\mB})\\
    \so \mSA \Satisfies
        \tSUM(g(\SIGMA{b_1}{\mB}),\dots,g(\SIGMA{b_n}{\mB}),
                g(\SIGMA{x}{\mB}))
            \since{$\mSA \submodel \mSB$}\\
    \but g(\SIGMA{b}{\mB}) = \SIGMA{f(b)}{\mA}
        \by{the choice of $f(b)$.}\\
    \so \mSA \Satisfies
        \tSUM(\SIGMA{f(b_1)}{\mA},\dots, \SIGMA{f(b_n)}{\mA},
                g(\SIGMA{x}{\mB}))\\
    \but \mA \Satisfies
        \tSUM(f(b_1),\dots,f(b_n),
                f(x))
        \since  {\mA\Satisfies\DISJU}\\
   \so \mSA \Satisfies
        \tSUM(\SIGMA{f(b_1)}{\mA},\dots,\SIGMA{f(b_n)}{\mA},
                \SIGMA{f(x)}{\mA})\\
    \so g(\SIGMA{x}{\mB}) = \SIGMA{f(x)}{\mA}
            \since{$\mSA \Satisfies \UNIQUEPL$}
    \end{flalign*}
\endProof

\Section{5.2}{Prime models for CA}
For each total, congruous remainder function, \f, we want to
find a prime model, \mQf, for \CAIf.
All of these prime models can be defined over the class, \Q,
of sets near quasi-congruence classes (see section 3.2).
For different remainder functions, we need to assign
different size relations over \Q.
Section 5.2.1 defines the structures \mQf and verifies that each
satisfies the respective theory, $\CAIf$.
Section 5.2.2 defines, for each model of $\CAIf$, a submodel, or \emph{shell}.
Section 5.2.3 shows that \mQf is isomorphic to the shell of any model
of $\CAIf$.

\Subsection{5.2.1}{The standard protomodels}
The construction here is a more elaborate version of the
construction of \ModelQ in chapter 3 (see 3.2.19).
\ModelQ turns out to be \mQf, where $f(n) = 0$ for
all \n.
As in the case of \ModelQ, the models \mQf and their copies in
arbitrary models of \CAI are  unions of chains.
\par
The sizes assigned to elements of \Q to induce \mQf for a total,
congruous remainder function, \f, are more elaborate than those
used in the definition of \ModelQ (see 3.2.22), but they are employed
in substantially the same way:

\Def{5.2.1}
\Enum
\item
    A \introSTerm{size} is an ordered pair \Pair{\rho}{\delta}, where
    both $\rho$  and $\delta$ are rational.
\item
    If $\theta_1 = \rd{1}$ and $\theta_2 = \rd{2}$ are
    sizes, then
    \Enum
    \item   $\theta_1 < \theta_2$ \iff
                $\rho_1 < \rho_2$
                or $\rho_1 = \rho_2$ and $\delta_1 < \delta_2$.
    \item   $\theta_1 + \theta_2 =
                \Pair{\rho_1 + \rho_2}{\delta_1 + \delta_2}$
    \endEnum
\endEnum
(cf 3.2.22)
\endDef

To assign sizes for \mQf we rely on the representation of sets in \Q
defined in 3.2.24.
\mQf, unlike \ModelQ, assigns different sizes to the \n-congruence
classes for a given \n.
\Def{5.2.2}
    If \f is total and congruous, then
\Enum
\item
    If $x$ is an  \n-congruence class,
    $x = \MS{k}{n}{i}$, then
    \[
        \thetaf(x) = \begin{cases}
                       \Pair{\nth}{\frac{n - f(n)}{n}}, &\text{if $i < f(n)$,} \\
                       \Pair{\nth}{\frac{-f(n)}{n}}, &\text{if $f(n) \leq i$.}

                        \end{cases}
    \]

\item
    If $x \in \QCn$, so that \x is the disjoint union of \n-congruence classes,
    \[
        x_1 \union \dots \union x_k
    \]
    then
    \[
        \thetaf( x) = \thetaf( x_1) + \dots +\thetaf( x_k)
    \]
\item
    If $x$ is finite, then
    \[
            \thetaf(x) = \Pair{0}{\Card{x}}
    \]
\item
    If $x \in Q$, so that \x can be represented as
    \[
        (C(x) \union \Dox) - \Dtx
    \]
    as in 3.2.24, then
    \begin{align*}
        \thetaf( x) &= \thetaf(C(x)) + \thetaf(\Dox) - \thetaf(\Dtx) \\
                       &= \thetaf(C(x)) + \Pair{0}{\Card{\Dox}}
                             - \Pair{0}{\Card{\Dtx}}  \\
                       &= \thetaf(C(x)) + \Pair{0}{\Card{\Dox}-\Card{\Dtx}}
    \end{align*}
\endEnum
\endDef

Intuitively, all \n-congruence classes are assigned sizes
\Pair{1/n}{\delta},
but $\delta$ is no longer 0 in all cases, as in \ModelQ.
Instead, the first $f(n)$ \n-congruence classes are each
one atom larger that the remaining $(n - f(n))$ \n-congruence
classes.
\par
The desired models of \CAI may now be defined:
\Def{5.2.3}
    If \f is total and congruous, then
    the \introSTerm{standard protomodel} for $f$ is \mQf, where
    \begin{align*}
    \Domain{\Q_f} &= \Q\\
    \mQf \Satisfies x < y &\tiff \thetaf(x) < \thetaf(y)\\
    \mQf \Satisfies x \eqsize y &\tiff \thetaf(x) = \thetaf(y)\\
    \mQf \Satisfies \tUNIT(x) &\tiff \thetaf(x) = \Pair{0}{1}\\
    \mQf \Satisfies \tSUM(x,y,z) &\tiff \thetaf(z) = \thetaf(x) + \thetaf(y)
    \end{align*}
    (cf 3.2.27)
\endDef

To verify that the structure \mQf is a model of
\CAIf, for total and congruous \f, we exhibit each such
model as the union of a chain of models.
\Def{5.2.4}
    If \f is total and congruous, then \mQfn is the submodel
    of \mQf whose domain is \Qn.
\endDef

\Fact{5.2.5}
    If \f is total and congruous and $n>0$, then
    $\mQfn \Satisfies \BASIC$.
\endFact
\Proof
The proof can be obtained from the proof of Theorem
3.2.29 by substituting:
\begin{align*}
\mQfn   &\emath{for} \mQn\emath{,}\\
\thetaf(x) &\emath{for} \theta(x)\emath{,}\\
\rho_f(x)   &\emath{for} \rho(x)\emath{, and}\\
\delta_f(x) &\emath{for} \delta(x)
\end{align*}
\endProof

The following notion is helpful in understanding our
constructions.
\Def{5.2.6}
    Suppose $\mA \Satisfies \BASIC$, $x \in \dA$, and
    $0 \leq m < n$.  Then an \introMTerm{$(n.m)$-partition}{nm-partition} of $x$ in $\mA$
    is a sequence $x_1, \dots , x_n$, where

\Enum
\item
       $x =  x_1 \union \dots \union x_n$
\item
    If $0 < i < j \leq n$, then $x_i \intersect x_j = \NullSet$, and
\item   The $x_i$ are approximately the same size:
    \Enum
    \item
        If $0 < i < j \leq m$, then $x_i \eqsize x_j$,
    \item
        If $m < i < j \leq n$, then $x_i \eqsize x_j$,
    \item
        If $0 < i \leq m < j \leq n$, then
        $\tSUM(x_i, a, x_j)$, for any atom $a$,

    \endEnum

\endEnum
In other words, \x is partitioned among \n pairwise
disjoint sets which are roughly the same size:
each of the first \m is one atom larger than each of the
remaining $n-m$.
If $0 < i \leq m$, $x_i$ is called a \introSTerm{charmed n-factor}
of \x; for $i > m$, $x_i$ is a \introSTerm{common n-factor} of \x.

\endDef

The sequence \mQfn does not constitute a chain of models.
For example $\ModelQ_{f,3}$ is not an extension of $\ModelQ_{f,2}$.
But this sequence does harbor a chain of models:
\Fact{5.2.7}
    If $n > m$, then $\mQfsub{n!} \submodel \mQfsub{m!}$.
\endFact
\Proof
    It will be clearer, and easier, to establish this
    by example rather than by formal proof.
    Letting $n = 2$ and $m=3$, we want to show
    that the 2-congruence classes have the same
    size relations in \mQfsub{6} as they do in
    \mQfsub{2}, for any \f.
    The other elements of \mQfsub{6} will then fall
    into place, since size relations are determined
    by the representations of each set, \x, as $C(x)$, \Dox,
    and \Dtx.
    \par
    Suppose that $f(2) = 0$, so that
    \[
        \mQfsub{2} \Satisfies \evens \eqsize \odds
    \]
    Since \f is congruous, $f(6) \in \theSet{0,2,4}$.
    If $f(6) = 0$, then all of the 6-congruence classes
    are common.
    If $f(6) = 2$, then \MSkn{6} and \MSn{6}{1} are
    the only charmed 3-congruence classes.
    If $f(6) = 4$, then all of the 3-congruence classes are charmed except for \MSn{6}{4} and \MSn{6}{5}.
    \par
    In any case, \evens will include the same
    number of charmed 3-congruence classes as
    \odds, so
    \[
        \mQfsub{6} \Satisfies \evens \eqsize \odds
    \]
    \par
    Suppose, however that $f(2) = 1$, so that
    \[
        \mQfsub{2} \Satisfies \evens \emath{is one atom
        larger than} \odds
    \]
    Here, $f(6) \in \theSet{1,3,5}$, since \f is
    congruous.
    In any case, \evens contains exactly one more
    charmed congruence class than \odds, so
    \[
        \mQfsub{6} \Satisfies \evens \emath{is one atom
        larger than} \odds
    \]
    \par
    So it goes in general.
\endProof

Fact 5.2.7 allows us to regard \mQf as the union of a chain of
models:
\Fact{5.2.8}
    If \f is total and congruous, then
    \[
      \mQf = \Unionst{\ModelQ_{f,{n!}}}{n > 0}
    \]
\endFact

\Fact{5.2.9}
    If \f is total and congruous, then
\Enum
\item
    \mQf \Satisfies \BASIC
\item
    \mQf \Satisfies \ADIV{k}, for $k >0$

\item
    \mQf \Satisfies \MODn{f(n)}, for $n >0$

\item
    \mQf \Satisfies \CAIf

\endEnum
\endFact
\Proof
\Enum
\item
    \BASIC is a universal-existential theory; so it
    is preserved under unions of chains. (See Appendix, Fact B.5).
\item
    Each \n-congruence class can be partitioned into
    \k $n k$-congruence classes.
\item
    $\mQfn \Satisfies \MODn{f(n)}$, by definition.
    Since \MODn{f(n)} is an existential sentence, it
    is preserved under extensions.
\item
    Immediate from (a), (b), and (c).
\endEnum
\endProof

\Subsection{5.2.2}{Shells of models}
To embed the model \mQf into an arbitrary model, \ModelA,
of \CAIf, we must find a smallest submodel, \ModelB, of \ModelA
which satisfies \CAI.
Clearly, the basis of \ModelA, call it $x_0$, must be included
in \ModelB, since the symbol $\Identity$ must refer to the same
set in the submodel as it does in the model.
But, if the basis of \ModelA is in \DomainB and \ModelB \Satisfies \CAI,
then \DomainB must contain two disjoint sets of roughly the same
size whose union is the basis of \ModelA, $x_0$.
Pick such a pair, $x_1$ and $x_2$, to include in \DomainB.
Whether these are exactly the same size or differ by an atom
depends on whether \ModelA \Satisfies \MOD{2}{0} or
\ModelA \Satisfies \MOD{2}{1}.
\par
\ModelB must also satisfy \ADIV{3}.
We can aim for this by placing in \DomainB three disjoint sets,
$x_{11}$, $x_{12}$, and $x_{13}$ whose union is $x_1$ and another
three disjoint sets,
$x_{21}$, $x_{22}$, and $x_{23}$ whose union is $x_2$ .
The existence of such sets is assured because
\ModelA \Satisfies \ADIV{3}.
Again, the exact size relations will be determined by
which MOD principles are satisfied in \ModelA.
Insuring that $x_1$ and $x_2$ are each divisible
by 3,  also guarantees that $x_0$ is divisible by 3:
the three unions
\[
    x_{11} \union x_{21},
    x_{12} \union x_{22},\emath{and}
    x_{13} \union x_{23}
\]
will be roughly the same size and will exhaust $x_0$.
\par
We can continue this process indefinitely, dividing each
set introduced at stage \n into $n+1$ roughly equal subsets
at stage $n+1$.
This will produce an infinite tree, bearing sets.
The deeper a node is in this tree, the smaller the set it
bears and the greater the number of successors among which this
set will be partitioned.
\par
This great tree of sets will not form a boolean algebra,
nor will it be closed under finite unions.
A boolean algebra could be obtained by including both the
node sets and their finite unions, but this would still not
be an atomic boolean algebra, which is what we are looking for.
We cannot correct for this problem by including in \DomainB all
atoms of \ModelA:
\ModelA may have uncountably many atoms while \ModelB,
to be a prime model, must be countable.
We leave the solution of this problem to the formal construction.
\par
The formal proof  proceeds as follows:
First, we define a tree, i.e. the set of nodes on which we
shall hang both the components of the successive partitions described
above and the atoms of the submodel
being constructed.
Second, we present the construction which, given a model \ModelA
of \CAIf, assigns a \introSTerm{node set} \NSP and a \introSTerm{node atom},
 \NAP, to each
node, \NP.
Third, we define the \introSTerm{shell} of \ModelA as the submodel of \ModelA
generated by the collection of node sets and node atoms.
In the next section, we show that the shell of \ModelA is isomorphic
to \mQf.
\par
First, the tree:

\Def{5.2.10}
\Enum
\item
    A \introSTerm{node} is a finite sequence \Ntok,
    where $k > 0$ and for all $i \leq k$, $n_i < i$.
    (The variables \NP and \NR
    range over nodes.)
\item
    If \NP = \Ntok, then
    \Enum
    \item
        The \introSTerm{length}, or \introSTerm{depth}, of \NP,
        \LNP is $k$.
    \item
        If $1 \leq i \leq k$, then $\NPC{i} = n_i$.
    \item
        $\NPN{m} = \tuple{ n_1, \dots, n_k, m }$
    \endEnum
\item
    \NP \introSTerm{extends} \NR iff $\LNP > \LNR$ and,
    for $1 \leq i \leq \LNR$,  $\NPC{i} = \NRC{i}$.
\item
    \NP \introMTerm{0-extends}{zero-extends} \NR iff \NP extends \NR
    and, for $\LNR < i \leq \LNP$,
    $\NPC{i} = 0$.
\endEnum
\endDef
Nodes constitute the vertices of an infinite tree
in which \tuplez  is the root and \NP dominates \NR
iff \NR extends \NP.
The number of immediate descendants of a node
grows as the depth of the node increases.
\par
We shall now assign a set to each node by repeatedly partitioning
the basis of \ModelA.
At the same time we shall assign an atom to each node.

\Def{5.2.11}
    Given $\mA = \mA_f \Satisfies \CAIf$, where $f$ is a
    total, congruous remainder function:
\Enum
\item
    Let \NS{\tuplez} be the basis of \ModelA, and let
    \NA{\tuplez} be any atom of \ModelA.
\item
    Suppose that \NSP and \NAP have been chosen.
    Let \m = \LNP and let
    \[
        k = \frac{f((m+1)!) - f(m!)}{m!}
    \]
    Since  \f is congruous,
    \k is an integer (see Definition 3.3.7e).
    Let
    \[
        \NS{\NPN{0}}, \dots, \NS{\NPN{m}}
    \]
    be an $(m+1,k)$ partition of \NSP if \NSP is common or
    an $(m+1,k+1)$ partition of \NSP if \NSP is charmed.
    Fact 5.2.12 guarantees that such partitions exist.
    In either case, choose \NS{\NPN{0}} so that it contains
    \NAP.
    This is always possible because \NSP contains \NAP.
\item
    Let \NA{\NPN{0}} = \NAP.
    If $0 < i \leq m$, let \NA{\NPN{i}} be any atomic subset of
    \NS{\NPN{i}}.
\endEnum
    The sets \NSP will be referred to as \introSTerm{node sets} and
    the atoms \NAP as \introSTerm{node atoms}.
\endDef

\Fact{5.2.12}
    Suppose \ModelA \Satisfies \CAIf, $n > 0$, $m > 0$, and
    \[
        k = \frac{f(nm) - f(n)}{n}
    \]
    Then
\Enum
\item
    \x has an $(m,k)$-partition \iff
        $\ModelA \Satisfies \Modm{k}(x)$
\item
    Every common \n-factor of \DomainA has an $(m,k)$-partition.
\item
    Every charmed \n-factor of \DomainA has an $(m,k+1)$-partition.

\endEnum
\endFact
\Proof
\Enum
\item
    $\Modm{k}(x)$ says that $x$ has an $(m,k)$-partition.
\item
    All common \n-factors are the same size and satisfy the same \Modm{k} predicate.
    So, suppose that each common \n-factor has \k charmed \m-factors
    and $m-k$ common \m-factors.

    By (a), \mA has $f(n)$ charmed \n-factors and $n-f(n)$ common \n-factors.
    Partitioning each of the \n-factors into \m subsets of roughly the same
    size yields an $nm$-partition of \mA; the charmed \m-factors of the \n-factors
    are the charmed $nm$ factors of \mA and the common \m-factors of the
    \n-factors are the common $nm$ factors of \mA.

    Each of the common \n factors has $k$ charmed \m-factors and each of
    the charmed \n-factors has $k+1$ charmed \m-factors.  In all, there are
     \begin{align*}
        &(n - f(n)) k + f(n)(k+1)        \\
       = &nk + f(n)
     \end{align*}
    charmed \m-factors among the \n-factors of \mA.

    But, by (a) again, there are $f(nm)$ charmed $nm$-factors of \mA.
    So
    \begin{align*}
        & f(nm) = nk + f(n) \\
    \AND    k = (f(nm) - f(n))/n
    \end{align*}
\item
    Immediate from (b), since every charmed \n-factor
    is one atom larger than each common \n-factor.
\endEnum
\endProof

The facts listed in 5.2.13 should clarify this constuction.
All of them can be established by induction on the depth of nodes.
\Fact{5.2.13}
\Enum
\item
    $\NSP \subset \NSR$ \iff \NR extends \NP or \NR = \NP.
\item
    If $\NSP = \NSR$, then \NR = \NP.
\item
    There are $n!$ nodes (and, hence, node sets) of depth \n.
\item
    If $i \neq j$, then $\NS{\NPN{i}} \intersect \NS{\NPN{j}}  =\emptyset$.
\item
    Any two node sets of the same depth are disjoint.
\item
    Each node set is the disjoint union of its immediate
    descendants.
\item
    Each node set is the disjoint union of all of its
    descendants at any given depth.
\item
    $\NAP \subset \NSR$ \iff \NP extends \NR.
\item
    $\NAP = \NAR$ \iff one of \NP and \NR 0-extends
    the other.
\item
    Every node set contains infinitely many node atoms.
\item
    For any \n, the node sets of depth \n form an
    $(n!, f(n!))$-partition of the basis of \ModelA.
\endEnum
\endFact
\Proof
We demonstrate (k) by induction on $n$.:
    If $n=1$, then $n!= 1$, $f(n!) = 0$, and $\APSUB{\tuplez}$  is
    the basis of \mA.
    So (\emph{k}) holds because any set is a (1,0)-partition of itself.
    \par
    Assume that (\emph{k}) holds for $n$.  Then (\emph{k}) also holds for $n+1$:
    each node set of depth $n$ has $n+1$ immediate descendants; so, there
    are $(n!)(n+1) = (n+1)!$  node sets of depth $n+1$.
    \par
    Furthermore,
    $f(n!)$ of the $n$-factors are charmed and $n! - f(n!)$ are common.
    By Fact 5.2.12, each charmed $n$-factor has an $(n+1,k+1)$-partition
    and each common $n$-factor has an $(n+1, k)$-partition, where
    \begin{align*}
        k &= \frac{f(n!(n+1)) - f(n!)}{n!} \\
        \\
            &= \frac{f((n+1)!) - f(n!)}{n!}
    \end{align*}
    (substituting $n!$ for $n$ and $n+1$ for $m$).
    \par
    So, there are
    \begin{align*}
        (k+1)f(n!) & \text{ charmed $n+1$ factors from the charmed $n$-factors}\\
        \text{and } k(n!-f(n!)) &\text{ charmed $n+1$ factors from the common $n$-factors}
    \end{align*}
    \par
    In all, then , the number of charmed $n$-factors is:
    \begin{align*}
    (k+1) f(n!) + k(n! - f(n!)) \\
        & = k f(n!) + f(n!) + n!k - f(n!)k \\
        & = f(n!) + n! k   \\
        &= f(n!) + f((n!) (n+1)) - f(n!)  \\
        & = f(n!) (n+1))    \\
        & = f((n+1)!)
    \end{align*}
    So, the $n+1$ factors of the $n$-factors of the basis form an
    $((n+1)!, f((n+1)!))$-partition of the basis.  That is to say,
    (\emph{k}) holds for $n+1$.
\endProof

Given a collection of node sets, \NSP, and node-atoms, \NAP,
from a model, \ModelA, of \CAIf, we can now construct
a submodel, \ModelB, of \ModelA which is isomorphic to \mQf.

\Def{5.2.14}
    Suppose \f is total and congruous, that \ModelA \Satisfies \CAIf
    and \NSP and \NAP are the node sets and node-atoms of \ModelA
    produced by construction 5.1.12.
    Then, the \introSTerm{shell} of \ModelA is
    the submodel of \ModelA generated by \theSet{\NSP,\NAP}.
\endDef
For the remainder of this chapter, we will regard as fixed:
\Enum
\item
    $\f$, a total, congruous remainder function
\item
    $\ModelA$, a model of \CAIf,
\item
    $\theSet{\NSP,\NAP}$, a set of node sets and node-atoms produced
    by the construction above.
\item
    $\ModelAh$, the shell of \ModelA generated from \theSet{\NSP,\NAP}
\endEnum

To show that $\ModelAh \iso \mQf$, we need a sharper characterization
of the elements of \ModelAh.
Recall from 3.2.24 that each member, \x, of \Q has a unique representation as
    \[
        (C(x) \union \Dox) - \Dtx
    \]
where $C(x)$ is a quasi-congruence class,  \Dox and \Dtx
are finite sets and
\begin{align*}
        \Dox \intersect C(x) = \NullSet\\
        \Dox \intersect \Dtx = \NullSet  \\
        \Dtx \subseteq C(x)
\end{align*}
\par
We can obtain a similar representation for elements of \ModelAh:
the node sets play the role of (some of) the congruence classes;
finite unions of node sets correspond to the quasi-congruence
classes; finite sets of node-atoms correspond to the finite
subsets in \Q.

\Def{5.2.15}
\Enum
\item
    \x is a \introSTerm{quasi-nodal} set of \ModelA \iff it is the
    union of finitely many node sets (iff it is the union of
    finitely many node sets at a given depth).
\item
    \x is an \introMTerm{\protect \mA-finite set}{A-finite set} \iff it is a finite
    set of node atoms.
\item
    If $x \in A$ and $y \in A$, then \x is
    \introMTerm{\protect \ModelA-near}{A-near} \y \iff both $x-y$ and $y-x$ are
    \ModelA-finite sets,
\endEnum
(cf. 3.2.14-3.2.18)
\endDef
Still following in the footsteps of chapter 3, we can characterize
\DomainAh as the collection of sets \ModelA-near quasi-nodal sets.
Analogues of 3.2.15 through 3.2.18 obtain for \ModelA-nearness.

\Fact{5.2.16}
\Enum
\item
    $x \in \DomainAh$ \iff \x is \ModelA-near some quasi-nodal set of
    \ModelA.
\item
    \ModelAh is an atomic boolean algebra whose atoms are the
    node-atoms \NAP.
\item
    If $x \in \DomainAh$, then \x has a unique representation as
   \[
        (C(x) \union \Dox) - \Dtx
    \]
    where $C(x)$ is a quasi-nodal set disjoint from \Dox
    and including \Dtx, both of which are \mA-finite sets.
\endEnum
\endFact
\Proof
\Enum
\item
    \IfWay
        \ModelAh is generated from node sets
        and
        node atoms via the boolean operations, each
        of which preserves \mA-nearness to quasi-nodal
        sets.
    \OnlyWay
        \ModelAh must contain finite unions of node sets
        as well as \mA-finite sets;
        so it must also contain sets obtained by
        adding or removing \mA-finite sets from
        quasi-nodal sets.
\item
    The proof parallels that of Theorem 3.2.21 exactly.
\item
    Let $C(x)$ be the quasi-nodal set which is
    \mA-near \x (see 3.2.23);
    Let $\Dox = x - C(x)$;
    and, let $\Dtx = C(x) - x$.
\endEnum
\endProof

\Subsection{5.2.3}{The embeddings}
To embed \mQf into \ModelA, we first describe \Q in terms of node sets
and node atoms.
In effect, we are performing the construction 5.2.12 on \mQf, but with
two differences:
first, we are stipulating which \nm-partitions to use at each level;
second, we are selecting node atoms so that every singleton in \Q is the
node atom for some node.
This latter condition guarantees that the
shell of \mQf will be \mQf itself.

\Def{5.2.17}
    Suppose \NP is a node. Then
\Enum
\item
    If $\LNP = k$, then
    \[
        \QP = \MSn{k!}{m}
    \]
    where
    \[
        m = \sum_{i = 1}^{k}{\NPI{i} (i-1)!}
    \]
\item
    The \introSTerm{depth} of $\QP$ is $\LNP$.
\endEnum
\endDef

\Examples
\begin{align*}
    \Nodal{\BQ}{\tuple{0}}& = \MSkn{}\\
    \Nodal{\BQ}{\tuple{0,0}}& = \MSkn{2}\\
    \Nodal{\BQ}{\tuple{0,1}}& = \MSn{2}{1}\\
    \Nodal{\BQ}{\tuple{0,1,0}}& = \MSn{6}{1}\\
    \Nodal{\BQ}{\tuple{0,0,2}}& = \MSn{6}{4}\\
    \Nodal{\BQ}{\tuple{0,1,2}}& = \MSn{6}{5}\\
\end{align*}
\endExamples

\Def{5.2.18}
\Enum
\item
   \IP = the least $n \in \QP$.
\item
    \qP = \theSet{\IP}
\endEnum
\endDef

\Fact{5.2.19}
\Enum
\item
    If $\NP = \Ntok$, then
    \[
        \IP = \sum_{i = 1}^{k}{(i-1)!n_i}
    \]
\item
    $\IP = \IR$ \iff $\NP = \NR$ or one 0-extends the other.
\item
    For every \n, there's a \NP such that $n = \IP$.
\item
    For every \n, there are infinitely many  \NP such that $n = \IP$.
\item
    For every \n, there's a \NP such that $n = \IR$ \iff
    \NR = \NP or \NR 0-extends \NP.
\endEnum
\endFact

\Fact{5.2.20}
\Enum
\item
    At each depth, \n, the \IP take on all and only values
    less than $n!$.
\item
    If $\IP = k$, then all nodes along the left-most branch
    descending from \NP also have the value \k.
    These are the only nodes below \NP with the value \k.
\item
    Every natural number is the value of all and only
    those nodes along the left-most branch descending from some
    node.
\endEnum
\endFact
Though for a given natural number \n, there will be infinitely many
nodes \NP for which $\IP = n$, we can associate with each
natural number a shortest (i.e. shallowest) node for which
$\IP = n$.

\Def{5.2.21}
    \NODESEQn = the shortest \NP such that $\IP = n$.
\endDef

\Fact{5.2.22}
\Enum
\item
    $\In{\NODESEQn} = n$
\item
    \NP extends $\NODESEQ{\IP}$
\item
    $\NODESEQ{\In{\NODESEQn}} = \NODESEQn$
\item
    $\NODESEQ{\IP} = \NP$ \iff $\NP = <0>$ or
    $\NPC{\LNP} \neq 0$
    (ie a node, \NP, will be the highest node with a certain
    value just in case \NP is not the leftmost immediate descendant
    of its parent.)
\endEnum
Each of the points listed in Fact 5.2.13 hold for the sets \QP
and \qP.
That is to say, the \QP can be regarded as node sets and the
\qP as node atoms for any model \mQf.
Notice, especially, that 5.2.13k holds.
\endFact

We may, finally, define the embedding of \mQf into \ModelA:
\Fact{5.2.23}
    If $x \in Q$, there is a unique $y \in \DomainAh$ such that
    \[
        \forall \NP(\qP \subseteq x \bicond
                \NAP \subseteq y)
    \]
\endFact
\Proof
There is at most one such $y$, by Fact 5.2.16.
To show that there is such a $y$, suppose first that $x \in \QC$.

Then there is some $n$ such that
\[
    x = x_1 \cup \dots \cup x_k
\]
where each $x_i$ is an ($n!$)-congruence class  and hence a node set in
\mQf.  So
\[
    x =   \QPSUB{1} \cup \dots \cup \QPSUB{k}
\]
Now, let
\[
    y =  \APSUB{1} \cup \dots \cup \APSUB{k}
\]
So $y \in \DomainAh$ and:
\begin{align*}
    \qP \subset x & \tiff \qP \subset \QPSUB{i} \etext{, for some $i$}\\
                & \tiff  \NP \etext{ extends } \NP_{i}
                            \by{5.2..13h}\\
                & \tiff  \NAP \subset \APSUB{i}  \\
                & \tiff \NAP \subset y.
\end{align*}
If $x$ is not a quasi-congruence class, then
\[
    x = (\xP \cup \Dox) - \Dtx
\]
where \xP is a quasi-congruence class.
Let \yP be the element of \DomainAh corresponding to \xP, as described above, and
let
\[
    y = (\yP \cup \theSetst{\NAP}{\qP \subseteq \Dox} ) 
        - \theSetst{\NAP}{\qP \subseteq \Dtx}
\]
\endProof

\Def{5.2.24}  The \introSTerm{nodal embedding} of \Q into
    \DomainA is defined by letting  $g(x \in \Q)$ be the $y \in \DomainAh$ such that
    \[
        \forall \NP( \NAP \subseteq y \bicond \qP \subseteq x)
    \]
\endDef

\Fact{5.2.25}
\Enum
\item
    If $x \in Q$, then
    \[
        g(x) = ( g(C(x)) \union g(\Dox)) - g(\Dtx)
    \]
\item
    \g is one-one.
\item
    \g maps \Q onto \DomainAh.
\endEnum
\endFact
\Proof
\Enum
\item
    Immediate from the proof of 5.2.23.
\item
    Suppose $g(x) = y = g(\xP)$. Then
    \[
        \forall \NP (\qP \subseteq x \bicond
            \NAP \subseteq y \bicond \qP \subseteq \xP)
    \]
    But every integer is \IP for some \NP, so $x = \xP$.
\item
    Obvious.
\endEnum
\endProof

\Theorem{5.2.26}
    The nodal embedding, \g, of \Q into \ModelA
     is an isomorphism of \mQf onto \ModelAh.
\endTheorem
\Proof
\Enum
\item
\begin{flalign*}
\mQf \Satisfies x \subseteq y
    &\tiff \forall \NP (\qP \subset x \then \qP \subset y)\\
    &\tiff \forall \NP (\NAP \subseteq g(x) \then \NAP  \subseteq g(y))
                                                &\text{(by 5.2.24)}\\
    &\tiff \ModelAh \Satisfies g(x) \subseteq g(y)
                                                &\text{(by 5.2.16)}
\end{flalign*}
\item
    As in (a),
    $g$ can be shown to preserve \NullSymbol, \Identity, unions,
    intersections, relative complements, and proper subsets.
\item
    Let $\nc(z,n)$ be the number of charmed node sets of level $n$
    contained in $C(z)$.
    And let $n$, below, be the least $k$ such that $x$ and $y$ are both
    unions of node sets of depth $k$.
    \begin{flalign*}
  \mQf \Satisfies x \eqsize y
            &\tiff  \theta_f(x) = \theta_f(y)        \\
            &\tiff  \nc(x,n) - (\Card{\Dox} - \Card{\Dtx})         \\
                 & \quad\quad     =   \nc(y,n) - (\Card{\Doy} - \Card{\Dty})      \\
            &\tiff   \nc(g(x),n) - (\Card{g(\Dox)} - \Card{g(\Dtx)})   \\
                  & \quad\quad    =   \nc(g(y),n) - (\Card{g(\Doy)} - \Card{g(\Dty)}) \\
                & \quad\quad\quad \text{since \QP is charmed iff \NSP is charmed}\\
                 & \quad\quad\quad \text{ and $g$ preserves boolean relations}      \\
            &\tiff   \mA \Satisfies g(x) \eqsize g(y)       \\
            &\tiff  \ModelAh \Satisfies g(x) \eqsize g(y)
                &\text{since } \ModelAh \submodel \mA
    \end{flalign*}
\item
    \begin{flalign*}
  \mQf \Satisfies x < y
            &\tiff   \mQf \Satisfies x \eqsize \xP \subset y  \\
                    & \quad\quad\quad \text{ for some }  \xP \in \Q
                &\text{since } \mQf \Satisfies \REPLT       \\
            &\tiff \mA \Satisfies  g(x) \eqsize  g(\xP) \subset g(y)
                 &\text{by (b) and (c)}      \\
            &\tiff  \mA \Satisfies g(x) < g(y)
                &\text{since } \mA \Satisfies \SUBSET, \INDISCEQ       \\
            &\tiff  \ModelAh \Satisfies g(x) < g(y)
                &\text{since } \ModelAh \submodel \mA
    \end{flalign*}
\item
    \begin{flalign*}
    \mQf \Satisfies \tSUM(x,y,z)
            &\tiff \mQf \Satisfies (\xP \cup \yP = z    \\
              & \quad\quad \Aand   x \eqsize \xP
                \Aand   y \eqsize \yP                   \\
              & \quad\quad    \Aand \xP \cap \yP = \NullSymbol)
                    &\text{since } \mQf \Satisfies \DEFPL \\
            &\tiff \ModelAh \Satisfies (g(\xP) \cup g(\yP) =g( z)  \\
              & \quad\quad    \Aand   g(x) \eqsize g(\xP )   \\
              & \quad\quad   \Aand   g(y) \eqsize g(\yP )                 \\
              & \quad\quad    \Aand g(\xP) \cap g(\yP) = \NullSymbol)
                    &\text{by (b) and (c)}  \\
            &\tiff \mA \Satisfies (g(\xP) \cup g(\yP) =g( z)  \\
               & \quad\quad   \Aand   g(x) \eqsize g(\xP )     \\
               & \quad\quad     \Aand   g(y) \eqsize g(\yP )                   \\
              & \quad\quad    \Aand g(\xP) \cap g(\yP) = \NullSymbol)          \\
            &\tiff \mA \Satisfies \tSUM(g(x), g(y), g(z))
                &\text{since } \mA \Satisfies \DEFPL     \\
            &\tiff \ModelAh \Satisfies \tSUM(g(x), g(y), g(z))
                &\text{since } \ModelAh \submodel \mA
\end{flalign*}
\endEnum
\endProof

So, each model of \CAIf has a submodel isomorphic to \mQf, and
we may conclude:

\Corollary{5.2.27}
    If \f is total and congruous, then \CAIf has a prime model.
\endCorollary

\Section{5.3}{Summary}
We can now draw our final conclusions about \CA, \CS, and their
completions.
\Theorem{5.3.1}
    If \f is a total, congruous remainder function, then \CAIf
    is consistent and complete.
\endTheorem
\Proof
    \CAIf is consistent because \f is congruous, by 3.3.12c.
    Since \CAIf is model complete, by 5.1.6, and
    has a prime model, by 5.2.27, the prime model
    test (Appendix, Fact C.3) applies.
    So \CAIf is complete.
\endProof

\Corollary{5.3.2}
    It \T is a completion of \CAI, then $T \equiv \CAIf$, for some
    congruous \f.
\endCorollary
\Proof
    For each $n >0$, $T \Proves \MODn{i}$ for exactly
    one $i$, $0 \leq i < n$:
    \par
    Since \T is complete, $T \Proves \MODn{i}$ or
    $T \Proves \neg \MODn{i}$ for each such $i$.
    But if $T \Proves (\neg \MODn{0} \Aand \dots
    \Aand \neg \MODn{n-1})$, then
    \T is inconsistent, since $T \Proves \CAI$ and
    $\CAI \Proves \DIVn$.
    Hence $T \Proves \MODn{i}$ for at least one
    one $i$, $0 \leq i < n$.
    \par
    But suppose  $T\Proves \MODn{i}$ and $T\Proves \MODn{j}$ where $0 \leq i \neq j < n$.
    Again, \T would be inconsistent, for
    $\CA \Proves \MODn{i} \then \neg \MODn{j}$
    (see Lemma 3.3.11).
    \par
    So, let $f(n) = m$ iff $T \Proves \MODn{m}$.
    Then $T \Proves \CAIf$ and, since
    \CAIf is complete, $\CAIf \Proves T$.
\endProof

\Corollary{5.3.3}
\Enum
\item
    Every completion of \CA is consistent with \CS.
\item
    $\CA \equiv \CS$
\endEnum
\endCorollary
\Proof
\Enum
\item
    Follows from 5.3.2 and 3.3.13.
\item
    Follows from (a) and 3.3.2b, given that
    $\CS \Proves \CA$.
\endEnum
\endProof

\Theorem{5.3.4}
\Enum
\item
    For $n > 0$, $\CS ; \EXACTLYn$ is decidable.
\item
    \CS is decidable.
\endEnum
\endTheorem
\Proof
\Enum
\item
    $\CS\adjoin\EXn \Proves \phi$ \iff
    $\mAn \Satisfies \phi$.
    But \mAn is a finite model.
\item
    To determine whether $\CS \Proves \phi$, alternate
    between generating theorems of \CA and
    testing whether $\mAn \Satisfies \neg\phi$.
\endEnum
\endProof

\Corollary{5.3.5}
\Enum
\item
    \CSI has $2^\omega$ completions.
\item
    For total \f, \CSIf is decidable \iff \f is decidable.
\endEnum
\endCorollary
\Proof
\Enum
\item
    There are $2^\omega$ remainder functions whose
    domain is the set of prime numbers.
    Each such function is congruous, so each
    corresponds to a consistent extension of \CSI.
    By Lindenbaum's lemma, each of these extensions
    has a consistent and complete extension.
\item
    \IfWay
        If \f is decidable, then \CAIf is recursively
        enumerable.
        But \CAIf is complete, so it is decidable.
    \OnlyWay
        To calculate $f(n)$, determine which \MODn{m}  sentence
       is in \CSIf.
\endEnum
\endProof

\Theorem{5.3.6}
    There is no sentence, \gphi, such that $T = \CA;\phi$ is consistent
    and \T has only infinite models.
\endTheorem
\Proof
    If \gphi is true only in infinite models of \CA,
    then $\neg \gphi$ is true in all finite models
    of \CA, so $\neg\gphi \in \CS$.
    But $CA \equiv \CS$, so $\CA\adjoin\gphi$ is
    inconsistent.
\endProof

\Theorem{5.3.7}
    \CS is not finitely axiomatizable.
\endTheorem
\Proof
    Suppose $\CS \Proves \gphi$.
    So $\CA \Proves \gphi$ and, by compactness,
    $(\BASIC \union T)\Proves \gphi$ for some
    finite set of \ADIVn principles,
    $T = \theSetst{\ADIVn}{n \in J}$.
    Let $K = \theSetst{n}{\mbox{every prime factor of
    \n is a member of \J}}$.
    \par
    Let \mA be a model with domain
    $\Unionst{\Qk}{k \in K}$,
    in which size relations are determined in accordance
    with the size function, \gtheta, defined in 3.2.26.
    We claim the following without proof:
    \begin{align*}
    \mA &\Satisfies \BASIC \tag{1}\\
    \mA &\Satisfies \ADIV{j} \emath{for all} j \in J
        \tag{2}\\
    \mA &\nSatisfies \ADIV{k} \tag{3}
    \end{align*}
    By (1) and (2), $\mA \Satisfies (\BASIC \union T)$,
    so $\mA \Satisfies \gphi$.
    But by (3), $\mA \nSatisfies \CS$.
    Hence $\gphi \nProves \CS$.
\endProof

\Chapter{6}{Sets of Natural Numbers}

\CS has standard finite models since it consists of sentences true in all such models.
It has infinite models \mQn and \mQf.
In this chapter we will show that \CS has infinite standard models over \PowerN.
\par
An ordering of \PowerN that satisfies \CS will not necessarily appear
reasonable.
For example, some such orderings say that there are fewer even numbers
than prime numbers (see below 6.2.13).
To rule out such anomalies, we introduce a principle, \OP, in section 1.
\OP mentions the natural ordering of $\fN$ and applies only to subsets
of \fN.
Section 2 establishes that \OP can be satisfied jointly with any
consistent extension of \CS in a model whose domain is \PN.
So \CS \adjoin \OP does not fix the size relations over \PN.
\Section{6.1}{The outpacing principle}
Throughout this chapter, $x$ and $y$ will range over \PN.
\Def{6.1.1}
$\x$  \introSTerm{outpaces} $\y$ just in case the restriction of \x to
  any sufficiently
large initial segment of $\fN$
 is larger than the corresponding restriction of \y, that is, iff:
\[
\exists n \forall m(m > n \then \Card{\Upto{x}{m}} >\Card{\Upto{y}{m}})
\]
\endDef

Notice that the size comparison between the two restricted sets will
always agree with the comparison of their normal cardinalities since all
initial segments of $\fN$ are finite.
\par
We employ this notion to state a \emph{sufficient} condition for
one set of natural numbers to be larger than another:
\namedSentence{\OP}
If \x outpaces \y, then \x > \y.
\endNamedSentence
The general motivation behind this principle should be familiar.
We extrapolate from well understood finite cases to puzzling infinite cases.
But we should also emphasize, again, that this extrapolation cannot
be done in any straightforward, mechanical way without risking
contradiction.
We cannot, for example, strengthen the conditional to a biconditional,
thus:
\numbSentence{1}
\x > \y \iff \x outpaces \y.
\endNumbSentence
This revised principle conflicts with \CS, for \emph{outpacing} is
not a quasi-linear ordering.
For example, neither \evens nor \odds outpaces the
other since each initial segment \theSet{0, \dots,2n+1} of $\fN$
contains \n evens and \n odds.
But the two are discernible under \emph{outpacing}, since \evens outpaces
\evensb while \odds does not.
\par
There is another point that underlines the need for care in
extrapolating from finite cases to infinite cases: we cannot
just use (2):
\numbSentence{2}
If, given any finite subset \z of $\fN$, \x restricted to \z is
larger than \y restricted to \z, then $x > y$.
\endNumbSentence
Though (2) is true, its antecedent is only satisfied when $\y \subset \x$.
\par
So, there are many statements that assert of infinite cases what is
true of finite cases.
Some of these conflict with one another.
Others are too weak to be helpful.
It is doubtful whether there is any mechanical way to decide which
of these statements are true.
The best we can do is propose plausible theories, determine whether
they are consistent, and see how far they go.
\setcounter{theorem}{-1}
\Def{6.1.0}
\Model{A} is an \introSTerm{outpacing model} \iff
\begin{align*}
\Domain{A} & = \PN \\
\Model{A}  & \Satisfies \BASIC, and \\
\Model{A}  & \Satisfies \OP
\end{align*}
\endDef
There is a slight difficulty in saying that an interpretation of \LCS
\squote{satisfies \OP}.
Since \OP involves the < relation over $\fN$, it cannot
be expressed in \LCS.
We shall finesse this problem by regarding \OP as the (very large)
set of sentences
\[ \theSetst{{\bf b} < {\bf a}}{a \etext{ outpaces } b} \]
\Fact{6.1.1a}
Every outpacing model satisfies the following
\Enum
\item
    $\evens > \MSkn{3}$
\item
    $\MSkn{3} > \MSkn{4} > \MSkn{5} > \dots $
\item
    If $k > 0$, then $\MSkn{k} > \squares$
\endEnum
\endFact
\Proof
\Enum
\item
    \evens has at least $(k-1)/2$ members
    less than or equal to \k for any given \k.
    \MSkn{3} has at most $k/3+1$ such members.
    If $k > 4$, then $(k-1)/2 > k/3 +1$
\item Similar to (1).
\item
    If $m = k^2$, both \MSkn{k} and \squares
    have exactly \k members less than \m.
    For $m > k(k+1)$, \Nto{m} will have more members in
    \MSkn{k} than in \squares.
\endEnum

\endProof

\Fact{6.1.2}
Every outpacing model satisfies the following
\Enum
\item $\evens \geq \odds$
\item $\odds \geq \evensb$
\item If $\evens > \odds$, then $\odds \eqsize \evensb$.
\endEnum
\endFact
\Proof
\Enum
\item By \TRICH, it is sufficient to show that $\neg\evens < \odds$.
    If it were, then, by \REPLT, there is a \y such that $y \eqsize \evens$
    and $y \subset \odds$.
    But any proper subset of \odds is outpaced by, and hence smaller than
    \evens.
    Let \k be the least odd number not in \y.
    Then \evens leads \y at $k+1$ and \y never catches up.
    So there is no \y such that $\evens \eqsize y \subset \odds$.
\item Similar to that of (1).
\item Note that $\evens = \evensb \union \theSet{0}$ and that \BASIC
        \Satisfies (*).
    \[
        (y \subset x \Aand z \leqsize y \Aand \tATOM(z')
            \Aand z \subset x \Aand x = z \union z') \then y \eqsize z  
    \]
    Let $\x = \evens$, $\y = \odds$, $\z = \evensb$ and apply (*).
\endEnum
\endProof
The two alternatives left open in 6.1.2 correspond to the possibilities
that $\fN$ may be odd or even:
if $\evens \eqsize \odds$, then $\fN$ is even,
if $\odds \eqsize \evensb$, then $\fN$ is odd.
In section 6.2, we show that both of these possibilities can be realized
in standard models over $\PN$.
Here, we generalize 6.1.2 to similar cases, including other congruence
classes.
\Def{6.1.3}
    $\Pair{x}{y}$ is an \introSTerm{alternating pair}
    \iff \x and \y are infinite and for all $\i > 0$,
    $x_i < y_i < x_{i+1}$.
\endDef
\Fact{6.1.4}
    If $\Pair{x}{y}$ is an alternating pair, then in any outpacing model:
    \[
        x \eqsize y \Or x > y \eqsize (x - {x_1})
    \]
\endFact
\Proof
The argument for 6.1.2 applies here since the only facts about \evens
and \odds used hold by virtue of these sets forming an alternating pair.
\endProof
\Theorem{6.1.5}
For a given $k > 0$, let $A_i = \MSn{k}{i}$ for each $i < k$.
Then:
\Enum
\item If $0 \leq i < j < k$, then
    \Pair{A_i}{A_j} is an alternating pair.
\item
    There is a \p, $0 < p \leq k$ such that
\Enum
\item
    If $i < j < p$, then $A_i \eqsize A_j$.
\item
    If $p \neq k$, then $A_0 \eqsize A_p \union \theSet{0}$.
\item
    If $p \leq i < k$, then $A_i \eqsize A_p$.
\endEnum
    (See example below).
\endEnum
\endTheorem

\Proof
\Enum
\item
    $A_i(n) = k n +i$,$A_j(n) = k n +j$,
    $A_i(n+1) = k n +i +k$ and
    $i < j < k + i$.
\item
    If $A_0 > A_i$ for some \i, let \p be the
    least such \i;
    otherwise, let $p = k$.
    \Enum
    \item
        For $0 \le i < p$, \Pair{A_0}{A_i} is an
        alternating pair.
        So either $A_0 > A_i$ or $A_0 \eqsize A_i$ by 6.1.4.
        But $A_0 \leqsize A_i$ by the selection of \p,
        so $A_0 \eqsize A_i$.
        (i) follows by \TRANSEQ.
    \item
        Immediate from 6.1.4 since \Pair{A_0}{A_p}
        is an alternating pair and $A_0 > A_p$.
    \item
        $A_0 > A_p \geq A_i$ if $i \geq p$.
        So $A_0 > A_i$.
        Hence $A_i \eqsize A_0 - \theSet{0} \eqsize A_p$.
        So $A_i \eqsize A_p$.
    \endEnum
\endEnum
\endProof
\Examp
Let $k = 4$, so $A_i = \MSn{4}{i}$ for $i = 0, 1,2,3$.
Then one of the following situations obtains:
\begin{gather*}
A_0 \eqsize A_1  \eqsize A_2 \eqsize A_3 > \MSn{4}{4}\\
A_0 > A_1  \eqsize A_2 \eqsize A_3 \eqsize \MSn{4}{4}\\
A_0 \eqsize A_1  > A_2 \eqsize A_3 \eqsize \MSn{4}{4}\\
A_0 \eqsize A_1  \eqsize A_2 > A_3 \eqsize \MSn{4}{4}
\end{gather*}
\endExamp
\Section{6.2}{Models of CS and Outpacing}
In this section, we construct models of \CS over \PN that satisfy
\OP.
\par
Outpacing models will be constructed out of finite models of \CS
by a technique which is very much like the ultraproduct construction
common in model theory, though the application here demands some
important differences.
\Def{6.2.1}
\Enum
\item
    \introLanguage{\LN}{subsets of \protect \fN} is the first order language which results from adding to
    \LCS, as individual constants, a name for each subset of \fN.
\item
    \ModelAn is the standard finite interpretation of \LN
    over \PNn in which:
    \[
        \ModelAn(\Constant{a}) = a \intersect \fNn = \Upto{a}{n}
        \mbox{\ for each\ } a \subseteq \fN.
    \]
\endEnum
\endDef

\Def{6.2.2}
(cf. \cite[Def. 18.15, p318]{Monk}
If \X is a set and $\FF \subseteq \PowerSet{X}$, then
\Enum
\item
    \FF has the \introSTerm{finite intersection property} iff the
    intersection of any finite subset of \FF is non-empty.
\item
    \FF is a \introSTerm{filter} over \X\iff
    \Enum
    \item
        $\FF \neq \NullSet$
    \item
        If $a \in \FF$ and $a \subseteq b$, then $b \in \FF$, and
    \item
        If $a \in \FF$ and $b \in \FF$, then
        $a \intersect b \in \FF$
    \endEnum
\item
    \FF is an \introSTerm{ultrafilter} over \X\iff
    \Enum
    \item
        \FF is a filter over \X
    \item
        $X \in \FF$, and
    \item
        if $Y \subseteq X$, then either $Y \in \FF$ or
        $(X - Y) \in \FF$.
    \endEnum
\item
    An ultrafilter, \FF, over \X is \introSTerm{principal} iff
    there is some $x \in \FF$ such that
    $\FF = \theSetst{a \subseteq X}{x \in a}$
\endEnum\endDef

\Fact{6.2.3}
\Enum
\item
    A non-principal ultrafilter contains no finite sets.
    \cite[Ch.6, lemma 1.3, p.108]{Bell}
\item
    A non-principal ultrafilter over \X contains all cofinite
    subsets of \X.
\item
    If $\FF \subseteq \PowerSet{X}$ and \FF has the finite intersection
    property, then there is an ultrafilter over \X which includes
    \FF.
    \cite[Prop. 18.18, p.319]{Monk}
\item
    If $Y \subseteq X$ and \Y is infinite, then there is a
    non-principal ultrafilter over \X which contains \Y.
\endEnum
\endFact

\Def{6.2.4}
    If \FF is an ultrafilter over \fN, then $\ModelAF$ is the
    interpretation of \LCS in which
    \Enum
    \item
        $\DomainAF = \PN$,
    \item
        $\ModelAF \Satisfies \Constant{a} < \Constant{b}$
       \iff
        $\theSetst{k}{\Upto{a}{k} < \Upto{b}{k}} \in \FF$, and
        similarly for other predicates.
    \item
        Boolean symbols receive the usual interpretation.
    \endEnum
\endDef

Our main result is that if \FF is non-principal, then \ModelAF is
an outpacing model.
\Theorem{6.2.5}
    If \FF is a non-principal ultrafilter over $\fN$, then
\Enum
\item
    If \gtau is a term of \LN, then
    \[
        \ModelAk(\tau) = \ModelAF(\tau) \intersect A_k = \Upto{\ModelAF(\tau)}{k}
    \]
\item
    If \gphi is a quantifier free formula of \LN, then
    \[
        \ModelAF \Satisfies \gphi \emath{iff}
            \theSetst{k}{\ModelAk \Satisfies \gphi} \in \FF
    \]
\item
    If \gphi is a universal  formula of \LN and
       $ \ModelAk \Satisfies \gphi $ for every $k$, then
       $\ModelAF \Satisfies \gphi $.
\item
    \ModelAF \Satisfies \REPLT.
\item
    $\ModelAF \Satisfies \BASIC $.

\item
    $\ModelAF \Satisfies \ADIVn$, for every \n.

\endEnum
\endTheorem
\Proof
\Enum
\item      
    By induction on the structure of \gtau:
    \Enum
    \item
        If \gtau is a constant, $\tau = a$ for some $a \subset \fN$.
        So  $A_k(\tau)  = \Upto{a}{n}$ by 6.2.1b.
    \item
        If $\tau = \schema{\tau_1 \union \tau_2}$,
        \begin{align*}
              \ModelAk(\tau) &= \ModelAk(\tau_1)  \union  \ModelAk(\tau_2)\\
                         &=  (\ModelAF(\tau_1) \intersect A_k ) \union (\ModelAF(\tau_2) \intersect A_k)\\
                         &=  (\ModelAF(\tau_1) \union \ModelAF(\tau_2)) \intersect A_k    \\
                         &=  \ModelAF(\tau_1 \union \tau_2)  \intersect A_k
        \end{align*}

    \endEnum
    The proofs for intersections and relative complements
    are similar.
\item      
    By induction on the structure of \gphi:
    \Enum
 \item
    If $\gphi = \schema{a \subset b}$,
    then $\ModelAF \Satisfies \gphi$ iff $a \subset b$

    If $a \subset b$, then there is  a $k \in b$ but not $\in a$.
    So if $n > k$, $\Upto{a}{n} \subset \Upto{b}{n}$.
    Hence, \theSetst{n}{\ModelAF \Satisfies \gphi} is cofinite and,
    by 4.2.3b, in \FF.

    Conversely, if $\theSetst{n}{\Upto{a}{n} \subset \Upto{b}{n}} \in \FF$,
    then it is infinite.
    So, there cannot be a $k$ in $a$ but not in $b$; otherwise \Upto{a}{n}
    would not be included in \Upto{b}{n} for any $n$ greater than $k$.
    So $a \subseteq b$.  But, clearly $a \neq b$, so $a \subset b$.
 \item
    If $\gphi = \schema{a = b}$, then
    \begin{align*}
            \ModelAF \Satisfies \gphi & \tiff a = b\\
                &\tiff \Upto{a}{n} = \Upto{b}{n} \etext{ for all $n$.}\\
                &\tiff \theSetst{n}{\mAn \Satisfies \noschema{a = b}} = \fN \\
                &\tiff \theSetst{n}{\mAn \Satisfies \noschema{a = b}} \in \FF
    \end{align*}
    since if $\mAn \nSatisfies \noschema{a = b}$ and $k > n$, then
          $\mAk \nSatisfies \noschema{a=b}$.

    \item
    $\ModelAF \Satisfies \noschema{a < b}$ iff
    $\theSetst{k}{A_k \Satisfies (a < b)} \in \FF$. Immediate from 6.2.4b.
    \item
    If \gphi is non-atomic, then, since  \FF is an ultrafilter:
    \begin{align*}
    \ModelAF \Satisfies \noschema{\gphi_1 \Aand \gphi_2}
        & \tiff \ModelAF \Satisfies \gphi_1 \emath{and} \ModelAF \Satisfies \gphi_2 \\
        & \tiff \theSetst{k}{\mAk \Satisfies \gphi_1} \in \FF \emath{and}
     \theSetst{k}{\mAk \Satisfies \gphi_2} \in \FF\\
        & \tiff \theSetst{k}{\mAk \Satisfies \noschema{\gphi_1 \Aand \gphi_2}} \in \FF \\
\\
    \ModelAF \Satisfies \noschema{\neg \gphi}
        & \tiff   \ModelAF \nSatisfies \noschema{\gphi} \\
        & \tiff \theSetst{k}{\mAk \Satisfies \gphi} \notin \FF \\
        & \tiff \theSetst{k}{\mAk \Satisfies \noschema{\neg \gphi}} \in \FF
    \end{align*}

    \endEnum
\item     
    \begin{flalign*}
    \Suppose \mAk  \Satisfies \noschema{\forall x \gphi(\x)} \tforall \k   \\
    \Then   \mAk  \Satisfies \noschema{\gphi(a)}  \tforall a \tforall k \\
    \so     \mAk  \Satisfies \noschema{\gphi(a)}  \tforall k \tforall a \\
    \so     \ModelAF  \Satisfies   \noschema{\gphi(a)}\tforall a  \by{(b)} \\
    \so     \ModelAF \Satisfies \noschema{\forall x \gphi(\x)}
    \end{flalign*}
\item     
    Suppose $\ModelAF \Satisfies (a < b)$.
    Construct $\aP$, a subset of $b$, for which $\ModelAF \Satisfies
    (a \eqsize \aP)$ as follows:
    \begin{align*}
    \emath{Let} K &= \theSetst{k}{\Upto{a}{k} < \Upto{b}{k}} \emath{so, $K \in \FF$.}\\
            &= \theSet{k_1, \dots , k_i, \dots} \emath{where the $k_i$'s
        are in strictly increasing order.}\\
        \emath{Let} a_0 &= \NullSet      \\
        a_{i+1} &= a_i \cup \theSet{\emath{the $n$ greatest members of}
                \Upto{b}{k_{i+1}}  \emath{not in} a_i} \emath{,}\\
                & \quad \quad\emath{where} n = \Card{\Upto{a}{k_{i+1}}} - \Card{a_i} \\
        a^\prime & = \Union a_i
    \end{align*}
   Then $a^\prime \subset b$ since each $a_i$ draws its new members from $b$.
   \par
    Claim: If $k \in K$, then $\Upto{a^\prime}{k} \eqsize \Upto{a}{k}$.
    \par
    Hence: $\ModelAF \Satisfies a^\prime \eqsize a$ since they are the same size over
    some set which contains $K$ and is, thus, if $\FF$.
\item     
    Immediate from (c) and (d): \BASIC is equivalent to a set of universal
    sentences, together with \ATOM and \REPLT. ($\ModelAF \Satisfies \ATOM$ because
    it contains all singletons of natural numbers.)
\item     
    For $0 < i \leq n$, $\x$  an infinite subset of \fN, let
    $x_i = \theSetst{x_{kn + i - 1}}{k \in \fN}$.

    The \n sets, $x_i$, partition \x.  Furthermore, these form an
    alternating \n-tuple, in the manner of the congruence classes modulo \n (see
    Theorem 6.1.5).
    As in 6.1.5, these sets are approximately equal in size and
    $\ModelAF \Satisfies \ADIVn(x)$.
\endEnum
\endProof

\Corollary{6.2.6}
    If \FF is a non-principal ultrafilter over \fN, then
    \Enum
    \item
        $\ModelAF \Satisfies \CAI$, and
    \item
        $\ModelAF \Satisfies \CSI$
    \endEnum
\endCorollary
\Proof
\Enum
\item
    Immediate from 6.2.5e and 6.2.5f.
\item
    Immediate from (a) and 5.3.3b.
\endEnum
\endProof
The proof of 6.2.6 is modeled on the usual ultraproduct construction,
but is not quite the same.
In the usual construction(see, for example, \cite[pp. 87-92]{Bell}), a model is
built by first taking the product of all factors (in this case, the
$\ModelA_k$), which results in a domain whose elements are functions
from the index set (\fN, here) to elements of the factors.
These functions are then gathered into equivalence classes (by
virtue of agreeing almost everywhere, i.e. on some member of the
filter) and the reduced ultraproduct is defined by interpreting
the language over these equivalence classes.
The model so constructed, which we'll call \PrAKF, has the
handy property that it satisfies any formula which is satisfied
by almost all factors, and certainly any formula which is satisfied
in all the factors.
This is handy because, given that each of the $\ModelA_k$ satisfies
\CS, we can immediately conclude that \PrAKF satisfies \CS.
\par
Unfortunately, \PrAKF is not the model we wanted:
its elements are not subsets of \fN, but equivalence classes
of functions from $\fN$ to finite subsets of \fN.
There is, indeed, a natural mapping from subsets of $\fN$ to such
elements, and this mapping would have allowed us to identify
a model over \PowerN as a submodel of \PrAKF; but only a submodel.
So, had we constructed the reduced ultraproduct, we would have then
been able to infer that the part of that model which held
our interest satisfied all universal formulas of \CS; we still
would have had to resort to special means to show that the
non-universal formulas were likewise satisfied.
\par
Fortunately, these special means were available;
the only non-universal axioms of \CAI could be verified in
the constructed model more or less directly, and the
completeness proof of the last chapter allowed us to
infer that all formulas true in all of the factors are true
in the model \ModelAF, after all.

\par
For a given  \FF, there will be many cases where
$\ModelAF \Satisfies \Constant{a} < \Constant{b}$ even though
\vb does not outpace \va.
This will happen whenever
$\theSetst{k} {\Upto{a}{k} <\Upto{b}{k}} \in \FF$ but is not cofinite.
Consider a familiar example:
Let $a = \odds$ and $b = \evens$.
Then $\Upto{a}{k} < \Upto{b}{k}$ iff $k \in \evens$, for
if we count the even numbers and the odd numbers up to some
even number, there will always be one more even and if we
count up to some odd number, there will always be the same
number of evens and odds.
\evens is neither finite nor cofinite, so it may or may not
be in \FF,
If $\evens \in \FF$, $\ModelAF \Satisfies \odds < \evens$.
Otherwise, $\odds \in \FF$,
so $\ModelAF \Satisfies \odds \eqsize \evens$.
\par
The construction of a model \ModelAF from any non-principal
ultrafilter, \FF, suggests that there are many outpacing
models unless different ultrafilters can yield the same model.
We will first show that this qualification is not needed.

\Theorem{6.2.7}
    If $\FF_1$ and  $\FF_2$ are distinct non-principal ultrafilters over \fN, then
    $\Model{A}_{\FF_1} \neq \Model{A}_{\FF_2}$.
\endTheorem
\Proof
(See below.)
\endProof

To show this, we will show that the presence of a set in an
ultrafilter makes a direct, personalized contribution to
the model \ModelAF.
Putting this in another way, there is a set of decisions that
must be made in constructing an outpacing model;
each decision may go either way, though the decisions are not
independent of each other.
Furthermore, each decision is made for a model \ModelAF
by the presence or absence of a particular set in \FF.

\Def{6.2.8}
    If $x \subseteq \fN$, then $\xp = \theSetst{i+1}{i \in x}$.
\endDef

A pair of sets, \xxp, can sometimes be an alternating pair,
but this is not always the case.
\Fact{6.2.9}
    \xxp is an alternating pair iff \x is infinite
    and there is no $n \in x$ such that $(n+1) \in x$.
    That is, no two consecutive numbers are in \x.
\endFact

Nevertheless, pairs \xxp are like alternating pairs in the
following way:
\Lemma{6.2.10}
If \x is infinite, $x1 = (x - x^+)$, $x2 = (x^+ - x)$,
\FF is a non-principal ultrafilter, and $\Model{A} = \ModelAF$, then
\Enum
\item $\Pair{x1}{x2}$ is an alternating pair.
\item $\ModA \Satisfies x \eqsize x^+$ or
    $\ModA \Satisfies x > x^+ \eqsize x - \{x_1\}$
\item $\ModA \Satisfies x > x^+ \tiff \x \in \FF$.
\endEnum
\endLemma
\Proof
\Enum
\item
    Let a \introSTerm{run} of \x be a maximal consecutive
    subset of \x. (So \evens has only 1-membered runs,
    while $\fN - \MSn{10}{1}$ has only 9-membered runs.)
    So, $x_1(n)$ is the first element in the \n-th run
    of \x and $x_2(n)$ is the first element after the
    \n-th run of \x.
\item
    We know from (a) that $x_2 \eqsize x_1$ or
    $x_2 \eqsize x_1 - x_1(1)$.
    But the disjoint union of $x\intersect\xp$ with
    $x_1$ or $x_2$, respectively, yields \x or \xp.
    So (b) follows from \DISJU.
\item
    We need only prove (*):
    \begin{gather*}
    \tag{*}
    \Card{\Upto{x}{n}} > \Card{\Upto{\xp}{n}}
    \tiff n \in x
    \end{gather*}
    Informally:
    \Upto{x}{n} first becomes greater than
    \Upto{\xp}{n} for $n = x(1)$, since $x(1) \notin \xp$
    because $\neg(x_1 - x_1(1) \subset x)$.
    Throughout the first run of \x, \x maintains
    its lead, losing it only at the least \n for which
    $n \notin x$ (since $(n-1) \in x$, so $n\in \xp$).
    This pattern repeats during successive runs of \x.
\endEnum
\endProof

We can now prove our main result, Theorem 6.2.7:
\Proof
    Without loss of generality, we can suppose there
    is a set, \x, such that $x \in \FFo$ and
    $x \notin \FFt$.
    By 6.2.10c, $\ModelAFo \Satisfies x > \xp$
    and $\ModelAFt \Satisfies x \eqsize \xp$
\endProof
Theorem 6.2.7 allows us to improve upon some previous
results.
For example, we can show that either of the alternatives
in 6.2.10b can obtained for any alternating pair.

\Theorem{6.2.11}
\Enum
\item
    If neither \x nor \y outpaces the other, then there is
    an ultrafilter \FF, such that $\ModelAF \Satisfies (x \eqsize y)$.
\item
    If \Pair{x}{y}  is an alternating pair, then there is
    an ultrafilter \FF, such that $\ModelAF \Satisfies (x > y)$.
\endEnum
\endTheorem
\Proof
\Enum
\item
    Let $J = \theSetst{k}
        {\Card{\Upto{x}{k}} =  \Card{\Upto{y}{k}}}$.
    Since neither \x nor \y outpaces the other,
    \J is infinite.
    By 6.2.3d, let \FF be a non-principal ultrafilter
    which contains \J.
    Then $\ModelAF \Satisfies x \eqsize y$.
\item
    If \Pair{x}{y} is an alternating pair, so is
    \Pair{y}{x-x(1)}.
    So, by (a) there is an \FF such that $\ModelAF
    \Satisfies y \eqsize (x - x(1)$.
    But then $\ModelAF \Satisfies x > y$.
\endEnum
\endProof

\Theorem{6.2.12}
Every infinite completion of \CS has an outpacing model.
\endTheorem
\Proof
    Recall that every infinite completion of \CS
    is equivalent to \CSIf for some total and
    congruous remainder function, \f.
    (See 3.6.2.)
   \par
    Given such an \f, let
    $G_k = \MSn{k}{f(k) + 1}$ for each $k > 0$.
    Then (*) holds for each \k:
    \begin{gather*}
    \tag{*}
    {G_k = \theSetst{n}{\ModelAn \Satisfies \MOD{k}{f(k)}}} 
    \end{gather*}
    Let $G = \theSetst{G_k}{k > 0}$.
    \par
    The intersection of any finite subset, $H$, of $G$
    is infinite.
    If $H$ is a finite subset of $G$, then
    $H = \theSetst{G_k}{k \in J}$, where \J is some
    finite subset of \fNp.
    So $H = \theSetst{n+1}{k \in J \then n \Cmod{f(k)}{k}}$.
    But \f is congruous, so the restriction of \f to
    the finite domain  \J has infinitely many solutions (see 3.3.8a).

\par
Since the intersection of any finite subset of \G
is infinite, there is a non-principal ultrafilter \FF
such that $\G \subseteq \FF$.
By (*), $\ModelAF\Satisfies\MODk{f(k)}$, so
$\ModelAF\Satisfies\CSIf$.
\endProof

On the basis of 6.2.11 we noted that there are even and odd
outpacing models;
we can now extend that observation to moduli other than 2.
More specifically, all of the possibilities listed
in 6.1.5 for the relative sizes of the \k-congruence classes
are obtainable in outpacing models.
\par
\par
This section has explored the existence and variety of outpacing
models.
Three comments are in order before we turn to the common structure
of outpacing models.
\par
First, even it \T is an infinite completion of \CS, there is no
unique outpacing model which satisfies \T.
This would be true only if fixing the congruence classes
determined whether $x \eqsize \xp$ or $x > \xp$
for every $x \subseteq N$.
That all such choices are not determined by a remainder theory
can be seen intuitively, perhaps, by considering
$x = \squares$:
any finite set of congruences has infinitely many solutions that
are squares and infinitely
many that are not; so whether $x \in \FF$ is an independent
choice.
Also, there are only $2^\omega$ remainder functions
while there are $2^{2^\omega}$ non-principal ultrafilters
over \fN\cite[Ch. 6, Theorem 1.5]{Bell},
 each yielding a different outpacing model.

\par
Second, it is not clear whether every outpacing model can be
obtained by the construction of 6.2.4.
Lemma 6.2.10c may suggest that any outpacing model, \ModelA, is
\ModelAF for
$\FF_{\ModelA} = \theSetst{x}{\ModelA \Satisfies x > \xp}$,
but it should not.
To establish that \ModelA = \ModelAF, both (1) and (2) are
necessary.
\numbSentence{1}
If \ModelA is an outpacing model, then $\FF_{\ModelA}$ is
a non-principal ultrafilter.
\endNumbSentence
\numbSentence{2}
If $\FF_{\ModelA} = \FF_{\ModelB}$, then $\ModelA = \ModelB$
\endNumbSentence
\par
I have not been able to prove (1) or (2). If (1) is false, then clearly
$\ModelA \neq \ModelAF$.
But even if (1) is true, two outpacing models may agree about all pairs
\xxp but disagree elsewhere.
At most one of them is obtainable by our construction.
So (1) and (2) are open problems.
\par
Finally, though it may be extraneous to show that \OP is independent,
we will do so.
\Theorem{6.2.13}
There are standard models of \CS over \PN which do not satisfy
\OP.
\endTheorem
\Proof
Suppose that \LOA is a linear ordering under which
$\fN$ forms an \gomega-sequence.
Then we could define \x \LOA-outpaces \y as follows:
\[
\exists n \forall m[n \LOA m \then
    \Card{\theSetst{k}{k \in x \Aand k \LOA m}}
    >
    \Card{\theSetst{k}{k \in x \Aand k \LOA n}}
    ]
\]
and define the principle:

\namedSentence{\OPLOA}
If $x \LOA\emath{-outpaces} y$, then $x>y$.
\endNamedSentence
Modifying 6.2.4, we could produce standard models of
\CS over \PN which satisfy \OPLOA and these will not, in
general, satisfy \OP.
\par
Suppose, for example, that \LOA is the ordering:
\[
    p_1, q_1,\dots,p_k,q_k,\dots
\]
where \p is the set of prime numbers and \q is its complement.
In any model of \OPLOA, \p and \q will be nearly the
same size, as the evens and the odds are in normal
outpacing models.
But the evens are much smaller than \q, so $p > \evens$
and \OP is false in \OPLOA models.
\endProof
\Section{6.3}{Size and density}
When number theorists talk about the sizes of sets of natural
numbers, they do not content themselves with speaking of the
(Cantorian) cardinalities of these sets.
Since they often want to compare infinite subsets of \fN, they
need a more discriminating notion.
\par
One notion they use is \introSTerm{asymptotic density}.
The asymptotic density of a set, \x, of natural numbers is the
limit, if there is one, of \Card{\Upto{x}{n}} as \n grows.
For example, the asymptotic density of \evens is $1/2$.
From now on we shall use the term \qword{density} for asymptotic density.
\par
In this section, we compare the ordering of \PN given by
density to the orderings given by \CS and \OP.
\Def{6.3.1}
\Enum
\item
    $\fract{x}{i} = \frac{\Card{\Upto{x}{i}}}{i}$,
    the fraction of numbers
    less than or equal to \i that are members of \x.
\item
    If $x \subseteq y \neq \emptyset$, then $\rho(x,y)$,
    the \introSTerm{density} of \x \emph{in} \y, is the limit,
    if it exists of
    \[
       \lim_{i = y_1}^{\omega} \fract{x}{i}/\fract{y}{i}
    \]
    That is,
    \[
        \rho(x,y) = r \tiff
        \forall( \delta > 0) \exists n (\forall i > n)
        (r -\delta < \frac{\fract{x}{i}}{\fract{y}{i}} < r + \delta)))
    \]
\item
    The \introSTerm{density} of \x, $\rho(x)$, is the density of
    of \x in \fN, if \x has a density in \fN.
\item
    If $x \subseteq y \neq \emptyset$, then
    \x \introSTerm{converges} in \y, $\cvg(x,y)$, iff \x
    has a density in \y;
    otherwise, \x \introSTerm{diverges} in \y.
\item
    \x \introSTerm{converges}, $\cvg(x)$, iff \x converges in \fN.
    \x \introSTerm{diverges} iff \x diverges in \fN.
\endEnum
\endDef
\Fact{6.3.2}
\Enum
\item
    If \x is finite, $\rho(x) = 0$.
\item
    If \x is cofinite, $\rho(x) = 1$.
\item
\begin{align*}
\rho(\evens)    &=  1/2\\
\rho(\MSkn{4}, \evens)   &=  1/2\\
\rho(\MSkn{4})   &=  1/4
\end{align*}
\item
    If $\cvg(x,y)$ and $\cvg(y,z)$, then $\cvg(x,z)$ and
    $\rho(x,z) = \rho(x,y) \rho(y,z)$.

\endEnum
\endFact

\Fact{6.3.3}
\Enum
\item
    There are divergent sets.
\item
    If $0 \leq r \leq 1$, there is a set with density \r.
\item
    If $0 \leq r \leq 1$ and \y is infinite,
    there is a set with density \r in \y.

\endEnum
\endFact
\Proof
\Enum
\item
    Let $x = \theSetst{i}{\exists n (10^{2n} \le i < 10^{2n+1}}$.
    So \x contains all numbers between 0 and 9, between 100 and 999,
    between 10,000 and 99,999, and so
    forth.
    If $n > 1$, then $\fract{x}{10^{2n}} \leq .1$ and
    $\fract{x}{10^{2n+1}} \geq .9$.
    So $\fract{x}{k}$ cannot have a limit.
\item
    Suppose \r is given. Construct the set \x as
    follows:
    \[
        x_0 = \NullSet
    \]
    \[
        x_{i+1} =  \begin{cases}
                x_i, &\text{if $\fract{x}{i} \geq r$;}\\
                x_i;i+1, &\text{otherwise.}
                \end{cases}
    \]
    \[
        x   = \Union {x_i}
    \]
\item
    Modify the construction for (b) in the obvious ways.
\endEnum
\endProof

\Theorem{6.3.4}
    Suppose that both \x and \y converge in \z.
    Then, if $\rho(x,z) < \rho(y,z)$, \y outpaces \x.
\endTheorem
\Proof
\begin{flalign*}
    \Let         b = (\rho(y,z) - \rho(x,z))/3\\
    \Let         \n_1 = \emath{the least \n such that for \i > \n,}
                   -b < \fract{x}{i}/\fract{z}{i} - \rho(x,z) < b \\
    \andLet        \n_2 = \emath{the least \n such that for \i > \n,}
                   -b < \fract{y}{i}/\fract{z}{i} - \rho(y,z) < b \\
    \So         \emath{for any} \i > \n_1 + \n_2,\\
    \wehave      \fract{x}{i}/\fract{z}{i} < \rho(x,z) +  b \\
    \AND        \fract{y}{i}/\fract{z}{i} > \rho(y,z) -  b \\
    \But        \rho(x,z) + b < \rho(y,z) - b,\\
    \so         \fract{x}{i}/\fract{z}{i} < \fract{y}{i}/\fract{z}{i}\\
    \so         \fract{x}{i} < \fract{y}{i}\\
    \so         \Card{\Upto{x}{i}} < \Card{\Upto{y}{i}}
\end{flalign*}
\endProof

We can use the relation between density and outpacing to draw
conclusions about densities that have nothing to do with outpacing,
as in Theorem 6.3.5.

\Theorem{6.3.5}
    If $\rho(x,z_1) < \rho(y,z_1)$, and both \x and \y converge
    in $z_2$, then
    $\rho(x,z_2) \leq \rho(y,z_2)$
\endTheorem
\Proof
Since $\rho(x,z_1) < \rho(y,z_1)$, \y outpaces \x.
But if $\rho(x,z_2) > \rho(y,z_2)$, then
\x outpaces \y.
So, $\rho(x,z_2) \leq \rho(y,z_2)$.
\endProof

We cannot strengthen the consequent of 6.3.5 to say
that $\rho(x,z_2) < \rho(y,z_2)$:
Let $z_1$ be the set of primes, let \x contain every third member
of $z_1$, and let \y be $z_1 - x$.
Then $\rho(x,z_1) = 1/3$ and $\rho(y,z_1) = 2/3$, but
$\rho(x,\fN) = \rho(y,\fN) = 0$.
\par
Theorem 6.3.4 implies that in any outpacing model, sets with
distinct densities will have distinct sizes.
Even if two sets have the same density, they will differ in
size if they have the same density in some common set.
So, from Facts 6.3.b and (c), we can begin to appreciate
how precise an ordering outpacing models provide:
\Fact{6.3.6}
    If \ModelA is an outpacing model, then
    \Enum
    \item
        there are uncountably many sizes of sets in \ModelA, and
    \item
        if $0 \leq r \leq 1$, then even among sets with density \r,
        there are uncountably many sizes in \ModelA.
    \endEnum
\endFact
\Proof
\Enum
\item
    Immediate from Fact 6.3.3b and Theorem 6.3.4.
\item
    Let \x be an infinite set with density \r, let \y be an infinite
    subset of \x, where $\rho(y,x) = 0$, and let $z = x - y$.
    There are uncountably many subsets of \y with distinct densities in \y,
    though $\rho(y_1, x) = 0$, for any $y_1 \subset y$.

    Suppose now that $y_1 \subset y$, $y_2 \subset y$, and $\rho(y_1) < \rho(y_2)$.
    Then $\mA \Satisfies y_1 < y_2$, by 6.3.4.
    so $\mA \Satisfies z \union y_1 < z \union y_2$, by \DISJU.

    But $\rho(z \union y_1) = \rho(z \union y_2) = \rho(x)$, since \z
    was obtained by removing from \x a set with density 0 relative to \x.
\endEnum
\endProof

\Subsection{6.3.1}{The Convexity Problem}
It's tempting to infer from these results that the extremely
fine ordering of sets by size (or, rather, any such ordering
which is realized in an outpacing model) is both a refinement
and a completion of the ordering suggested by density:
a refinement because it preserves all differences in size which
are captured by the notion of density, a completion
because all sets are located in a single, linear ordering of
sizes.
\par
But the situation is really not so clear.
It is evident that the size ordering over \PowerN in any outpacing
is a refinement of the ordering by cardinality:
if \x has a smaller cardinal number than \y, then \x is smaller than
\y, though two sets with the same cardinal number may have
different sizes.
We can regard the cardinality of a set in \PowerN
as determined by its size, though
different sizes may yield the same cardinal number.
We shall express this fact by saying that, at least when
we focus on \PowerN, cardinality is a function of size.
(It is not at all clear that this is true in any power set.)
\par
Now, we want to know whether the density of a set is a function of
its size.
It turns out that a negative answer is compatible with our theory
(\CS and \OP), while an affirmative answer may or may not
be consistent.
First, we will give a more precise formulation of this problem;
second, we will show that a negative answer is consistent;
finally, we'll discuss the consistency of an affirmative answer.
In passing, we'll explain why this is called \squote{the convexity problem}.
\par
Consider (1) and (1a):
\numbSentence{1}
If $x \eqsize y$ and $\rho(x) = r$, then $\rho(y) = r$.
\endNumbSentence
\numbSentence{1a}
If $x \eqsize y$ and $\rho(x) = r$ and \cvg(y), then $\rho(y) = r$.
\endNumbSentence
(1a) is an immediate consequence of Theorem 6.3.4.
For if \y converges, there is some $r_2 = \rho(y)$;
if $r < r_2$, then $x < y$ and
if $r > r_2$, then $x > y$.
But $x \eqsize y$, so $r_2 = r$.
\par
So, the questionable part of (1) can be expressed as (2):
\numbSentence{2}
If $x \eqsize y$ and \x converges, then \y converges.
\endNumbSentence
Theorem 6.3.4 insures that (1) just in case (2).
Recalling that sets may have the same density even though they differ
in size, we may consider two additional formulations:

\numbSentence{3}
If $\rho(x) = \rho(y)$ and $ x  < z < y$, then
$\rho(z) = \rho(x)$
\endNumbSentence
\numbSentence{4}
If $\rho(x) = \rho(y)$ and $ x  < z < y$, then
\z converges.
\endNumbSentence
(3) and (4) are equivalent for the same reasons that (1) and (2)
are equivalent.
\Fact{6.3.7}
    An outpacing model satisfies (1) iff it satisfies (3).
\endFact
\Proof
\BothWays
\IfWayB
    Suppose that
    $x < y < z$ and $\rho(x) = \rho(z)$.
    \par
    By \REPLT, there are two sets, \xP and \yP, for which
    \[
        \xP \subset \yP \subset z, x\eqsize\xP,
        y \eqsize \yP.
    \]
    \begin{flalign*}
    \Since  x \eqsize \xP\\
    \AND    \rho(x) = \rho(z)\\
    \Then   \rho(\xP) = \rho(z) \by{(1)}\\
    \So     \rho(z-\xP) = 0\\
    \AND    \rho(z-\yP) = 0\\
    \But    \yP = z - (z-\yP)\\
    \so     \rho(\yP) = \rho(z)\\
    \so     \rho(y) = \rho(z)   \by{(1)}
    \end{flalign*}
\OnlyWayB
    Suppose that
    $ x\eqsize y$ and $\rho(x) = r$.
    \par
    If $x = \NullSet$ or $ x=\fN$, then $y=x$,
    so $\rho(y) = \rho(x)$
    \par
    To apply (3), we need to find two sets, $x_1$ and
    $x_2$ such that
    \[
        \rho(x_1) = \rho(x_2) = r = \rho(x)
    \]
    and
    \[
        x_1 < y < x_2
    \]
    Assuming that $x \neq \NullSet$ and $x \neq \fN$,
    let
    \[
        x_1 = x - {a}\emath{, for some} a \in x
    \]
    and
    \[
        x_2 = x \union {b}\emath{, for some} b \notin x.
    \]
    Then $x_1 < x < x_2$, by \SUBSET, and
    $x_1 < y < x_2$, by \INDISCEQ.
    \par
    Since adding or removing a single element has no
    effect on the density of a set,
    \[
        \rho(x_1) = \rho(x) = \rho(x_2)
    \]
    So, by (3), $\rho(y) = r$.
\enBoth
\endProof
The question at hand is called the \qword{convexity problem} because of
the formulation in (3).
In geometry, a figure is convex if, given any two points in the figure,
any point between them (ie on the line segment from one to the other)
is also in the figure.
Applying this notion in the obvious way to the size ordering,
(3) says that the class of sets having a given density is convex.
By theorems 6.3.4 and 6.3.7, (1), (2) and (4) say the same thing.
\par
We regret that the only thing we know about (3) is that it may
be false:
\Theorem{6.3.8}
The negation of (1) is satisfied in some outpacing model.
\endTheorem
\Proof
    Let \x be a set with density \r and let \y be a divergent
    set which neither outpaces nor is outpaced by \x.
    Let \K be
    \[
        \theSetst{n}{\Upto{x}{n} = \Upto{y}{n}}
    \]
    \K is infinite, so \K is a member of some non-principal ultrafilter,
    \FF.
    $\ModelAF \Satisfies x \eqsize y$, so (1) is false in \ModelAF.
\endProof
\par
Open Problem:
Is (1) consistent with \CS and \OP?

\chapter*{Appendix}
\markboth{APPENDIX}{RightOddAppendix}
\addcontentsline{toc}{chapter}{Appendix}
\renewcommand{\thesection}{\Alph{section}}
\setcounter{section}{0}
\Section{A.1}{Notation}
\Subsection{A.1.1}{Predicate Logic}
The formal theories discussed in this thesis are theories with standard
formalization, in the sense of \cite[p.5]{Tarski}.
That is, they are formalized in first order predicate logic with identity
and function symbols.  The following notation is used for the predicate
calculus:
\begin{align*}
\Aand &\quad \etext{and} \\
\Or &\quad \etext{or}      \\
\neg    &\quad  \etext{not} \\
 \then  &\quad \etext{if ... then}  \\
\bicond         &\quad   \etext{if and only if}\\
=       &\quad  \etext {identical} \\
 \neq      &\quad \etext{not identical}  \\
\exists x &\quad \etext{there is an $x$}  \\
\forall x &\quad \etext{for all $x$}
\end{align*}
Conjunctions and disjunctions of sets of sentences are represented by:
\[
\Andst{\phi}{\etext{... $\phi$ ...}}
\]
and
\[
\Orst{\phi}{\etext{... $\phi$ ...}}
\]
\par
A first order language is determined by its non-logical constant symbols,
in the usual way.
These may be predicates, individual constants, or function symbols.
In most cases, the ranks of the symbols will accord with their familiar
uses and we do not stipulate them.
\par
A \introSTerm{schematic function} is a function whose range is a set
of formulae.
\DIVn, for example, maps natural numbers into sentences (see 3.2.8).
The arguments of such functions are indicated by subscripts.

 \Subsection{A.1.2}{Set Theory}
 For first order languages with boolean operations and predicates,
 we use the following notation:
 \begin{align*}
 \Identity &\quad \etext{the universe}      \\
  \NullSymbol &\quad \etext{the empty set}      \\
   x \cup y &\quad \etext{the union of $x$ and $y$}      \\
    x \cap y &\quad \etext{the intersection of $x$ and $y$}      \\
     x - y &\quad \etext{the relative complement of $y$ in $x$}      \\
      x \subseteq y &\quad \etext{$x$ is contained in $y$}      \\
       x \subset y &\quad \etext{$x$ is a proper subset of $y$}
 \end{align*}
 These symbols are also used for set theoretic relations, outside of
 first order languages.
 
 In addition, we use the following:
 \begin{align*}
        x \in y &\quad \etext{$x$ is a member of $y$}      \\
         \PowerSet{x} &\quad \etext{the power set of $x$}      \\
         \theSetst{x}{\phi(x)} &\quad \etext{the set of $x$'s which are $\phi$} \\
  x \adjoin a  &\quad \etext{the union of $x$ and ${a}$} \\
  \Upto{x}{n} &\quad \etext{the members of $x$ less than or equal to $n$}\\
  \fN &\quad \etext{the set of natural numbers: 0, 1, ...}\\
  \fNp &\quad \etext{$\fN - {0}$} \\
  \fZ &\quad \etext{the set of integers} \\
  \omega &\quad \etext{the smallest infinite cardinal}\\
  2^\omega &\quad \etext{2 to the $\omega$}
\end{align*}
 
 \Subsection{A.1.3}{Arithmetic}
For arithmetic, we use the following:
\begin{align*}
i + j   &\quad \etext{ $i$ plus $j$} \\
i - j  &\quad \etext{ $i$ minus $j$} \\
i j   &\quad \etext{  $i$ times $j$} \\
i^j &\quad \etext{   $i$ to the  $j$-th power} \\
i!  &\quad \etext{$i$ factorial} \\
i \divides j &\quad \etext{$i$ divides $j$}\\
\etext{gcd}(i,j) &\quad \etext{greatest common divisor of $i$ and $j$}\\
\n \Cmod{m}{j}  &\quad \etext{$n$ is congruent to $m$ modulo $j$}
\end{align*}
We use the following notation for sets of natural numbers:
\[
    \theSeq{...\Schema{n}...} = \theSetst{k}{\exists n(k = ... n ...)}
\]
For example:
\begin{align*}
\MSn{k}{j} &= \etext{the set of numbers congruent to $j$ modulo $k$}   \\
\squares &= \etext{the set of squares of natural numbers}
\end{align*}
\Section{A.2}{Model Theory}
This section lists the model-theoretic notions and
results assumed in the text and presents our notation
for these notions.
\par
An \introSTerm{interpretation}, \mA, of a first order language,
\L, consists of a domain, \dA, and a function which
assigns to each individual constant of \L a member of
\dA, to each \n-place predicate of \L, a set of
\n-tuples over \dA, and to each \n-argument function
symbol of \L, an \n-ary function defined over and yielding values
in \dA.
\par
We assume the notion of \introSTerm{satisfaction} as
defined in \cite[section 1.3]{Chang} and use
\[
    \mA \Satisfies \gphi \emath{to mean that \mA
    satisfies \gphi}.
\]
\Def{A.2.1}Familiar notions about models.
\Enum
\item
    \mA is a \introSTerm{submodel} of \mB
    ($\mA \submodel \mB$).
\item
   \mB is an \introSTerm{extension} of \mA
    ($\mA \submodel \mB$).
\item
    The \introSTerm{submodel} of \mB \introSTerm{generated}
    by \X, where \X is a subset of \dB.
\item
   \mA is  \introSTerm{isomorphic} to \mB
    ($\mA \iso \mB$).
\item
    \f is an \introSTerm{isomorphic embedding} of \mA into
    \mB.
\item
   \mB is an \introSTerm{elementary extension} of \mA.
\item
    \theSet{\mA_i} is a \introSTerm{chain} of models.
\item
    $\mA = \Union{\theSet{\mA_i}}$; the
    \introSTerm{union of a chain} of models.
\endEnum
\endDef
\Def{A.2.2}Notions about theories:
\Enum
\item
    A \introSTerm{theory} is a set of first-order sentences.
\item
    The \introSTerm{language of a theory}, \theLang{T}.
\item
    \T \introSTerm{proves} \gphi ($\T \Proves  \gphi$);

    \T \introSTerm{proves} \Tt ($\T \Proves  \Tt$).
\item
    \T is \introSTerm{complete}.
\item
    \T is \introSTerm{consistent}.
\item
    \T is \introSTerm{categorical}.
\item
    \T is \introSTerm{equivalent} to \Tt ($\T \equiv  \Tt$).
\endEnum
\endDef

\Fact{A.2.3}The following are well known facts:
\Enum
\item
    If $\T \Proves \gphi$, then some finite subset
    of \T  proves \gphi(\introSTerm{compactness}).
\item
    If \T is complete and $T\adjoin\gphi$ is
    consistent, then $T\Proves\gphi$.
\item
    If \T is categorical, then \T is complete.

\endEnum
\endFact

\Def{A.2.4}More familiar notions:
\Enum
\item
    \T is \introSTerm{existential} iff it is equivalent
    to a set of prenex sentences,
    none of which have universal quantifiers.
\item
    \T is \introSTerm{universal} iff it is equivalent
    to a set of prenex sentences,
    none of which have existential quantifiers.
\item
    \T is \introSTerm{universal-existential} iff it is equivalent
    to a set of prenex sentences,
    each of which has all of its universal quantifiers
    preceding any of its existential quantifiers.
\item
    \gphi is a \introSTerm{primitive formula} iff
    \gphi is an existential formula in prenex form
    whose matrix is a conjunction of atomic formulas
    and negations of atomic formulas.

\endEnum
\endDef

\Fact{A.2.5}Fairly familiar facts:
\Enum
\item
    If \T is existential, $\mA\Satisfies T$, and
    $\mA\submodel\mB$, then $\mB\Satisfies T$.
\item
    If \T is universal and $\mA\Satisfies T$,
    then any submodel of \mA also satisfies \T.
\item
    If \T is universal-existential and $\mA_i \Satisfies T$, for all $i > 0$, then
    $\Union{\theSet{\mA_i} }\Satisfies \T$.
\endEnum
\endFact

\Section{A.3}{Model Completeness}
This section presents the definition of \emph{model
completeness} and some basic results needed in Chapters 4 and 5.
\par
There are several ways to show that a theory is
complete.  We use only one, which is based on
A. Robinson's notion of \emph{model
completeness}.
\Def{A.3.1}
\T is \introSTerm{model complete} iff \T is consistent
and for any two models, \mA and \mB of \T,
$\mA \submodel \mB$ iff \mA is an elementary submodel
of \mB.
\cite[p. 355]{Monk}.
\endDef
A model complete theory is not necessarily complete,
unless the theory has a prime model.
\Def{A.3.2}
\mA is a \introSTerm{prime model} of \T if $\mA \Satisfies T$
and \mA can be embedded in any model of \T.
\cite[p.359]{Monk}
\endDef
\Fact{A.3.3}
If \T is model complete and \T has a prime model,
then \T is complete.
\endFact

To show that a theory is model complete we rely,
directly or indirectly, on a theorem of Monk's
(see 4.2.1a), which is based on some
equivalent formulations of model completeness:
\Def{A.3.4}
\Enum
\item
    The \mA-\introSTerm{expansion} of \L is the result
    of adding to \L a constant for each element of
    \dA.
\item
    The \introSTerm{diagram} of \mA is the set of
    atomic sentences and negations of atomic sentences
    of the \mA-expansion of the language of \mA which
    are true in \mA.
\endEnum
\endDef

\Fact{A.3.5}The following are equivalent \cite[p. 356]{Monk}:
\Enum
\item
    \T is model complete.
\item
    For every model, \mA, of \T and every \mA-expansion
    \Lpr of \L, the $T \union \Lpr\emath{-diagram of}L$
    is complete.
\item
    If \mA and \mB are models of \T, $\mA \submodel\mB$,
    \gphi is a universal formula, $x \in \mA$,
    and $\mA \Satisfies \gphi(x)$, then
    $\mB \Satisfies \gphi(x)$.
\item
    If \mA and \mB are models of \T, $\mA \submodel\mB$,
    \gphi is a primitive formula, $x \in \mA$,
    and $\mB \Satisfies \gphi(x)$, then
    $\mA \Satisfies \gphi(x)$.
\endEnum
\endFact

\backmatter


\renewcommand{\thesection}{}
\renewcommand{\thesubsection}{}

  \newcommand\glohead[1]{\indexspace \textbf{#1}}
\newcommand\refto[1]{\pageref{L#1}}
\newcommand\gloitem[2]{\item #1, \refto{#2}}
 \renewcommand\indexname{Index of Symbols}
 \begin{theindex}
 \addcontentsline{toc}{chapter}{Index of Symbols}

\glohead {Logic}
\gloitem {$\phi \Aand \psi$}{SA.1.1}
\gloitem {$\phi \Or \psi$}{SA.1.1}
\gloitem {$\phi \then \psi$}{SA.1.1}
\gloitem {$\neg \phi$}{SA.1.1}
\gloitem {$\phi \bicond \psi$}{SA.1.1}
\gloitem {$\forall x \phi$}{SA.1.1}
\gloitem {$\exists x \phi$}{SA.1.1}

\glohead{Set Operations}
\gloitem {$x \cup y$}{SA.1.2}
\gloitem {$x \cap y$}{SA.1.2}
\gloitem {$x - y$}  {SA.1.2}
\gloitem {$x\adjoin y$}{SA.1.2}
\gloitem {$x \Cup y$}{SA.1.2}
\gloitem {\PowerSet{x}}{SA.1.2}
\gloitem {\Card{x}}{SA.1.2}

\gloitem {\Upto{x}{k}}{SA.1.2}
\gloitem {$\fract{x}{k}$}{6.3.1}

\glohead{Arithmetic}
\gloitem {$k \divides n$}{SA.1.3}
\gloitem {$k!$}{SA.1.3}

\gloitem {\MSn{i}{j}}{SA.1.3}

\glohead {Size ordering}
\gloitem {$x \eqsize y$}{S1.1}
\gloitem {$x < y$}{S1.1}
\gloitem {$x \leq y$}{S1.1}
\gloitem {$x > y$}{S1.1}
\gloitem {$x \geq y$}{S1.1}

\glohead{Model Theory}
\gloitem{$T \Proves \phi$}{A.2.2}
\gloitem{$T \nProves \phi$}{A.2.2}
\gloitem{$T_1 \equiv T_2$}{A.2.2}
\gloitem {$\mA \Satisfies \phi$}{A.2.1}
\gloitem {$\mA \nSatisfies \phi$ }{A.2.1}

\gloitem {$\mA \submodel \mB$}{A.2.1}
 
 \glohead{Models and Domains}
\gloitem {\mA}{A.2.1}
\gloitem {\Domain{A}}{A.2.1}
\gloitem {\Basis{A}}{3.1.2}

\gloitem {$\mAn$}{6.2.1}
\gloitem{\ModelAF }{6.2.4}
 
\gloitem{\ModelAh }{5.2.14}
\gloitem {\Q }{3.2.19}
\gloitem {\mQf }{5.2.3}
\gloitem {\mQfn }{5.2.4}
\gloitem {\Qn }{3.2.19}
\gloitem {\mQn }{3.2.27}
\gloitem {\QC }{3.2.12}
\gloitem {\QCn }{3.2.12}
 \gloitem {\Naturals }{3.2.8}
\gloitem {\ModelSA }{4.1.1}
\gloitem {\ZgmA }{4.3.2}
\gloitem {\mZgastar}{4.3.4}
\glohead{Sizing}
\gloitem{$\Comp{x}$}{4.1.1}
\gloitem {$C(x)$}{3.2.24}
\gloitem {\Dox}{3.2.24}
\gloitem {\Dtx}{3.2.24}
\indexspace

\gloitem {$\rho(x)$}{3.2.26}
\gloitem {$\delta(x)$}{3.2.26}
\gloitem {$\theta(x)$}{3.2.26}
\gloitem {$\thetaf(x)$}{5.2.2}
\gloitem {$\alpha(x)$}{3.2.25}
\gloitem {$\beta(x)$}{3.2.25}

\gloitem {$\rho(x,y)$}{6.3.1}
\gloitem {$\SIGMA{x}{\mA}$}{4.1.1}
\gloitem {$\cvg(x,y)$}{6.3.1}
\gloitem {$\cvg(x)$}{6.3.1}

\indexspace
\glohead{Node Construction}
\gloitem {\NP}{5.2.10}
\gloitem {\LNP}{5.2.10}
\gloitem {$\NPN{m}$}{5.2.10}

\gloitem {\NSP}{5.2.11}
\gloitem {\NAP}{5.2.11}

\gloitem {\ModelAh}{5.2.14}
 
\gloitem {\IP}{5.2.18}
\gloitem {\qP}{5.2.18}

\gloitem {\NODESEQn}{5.2.21}
\indexspace
 \end{theindex}

 \renewcommand\indexname{Index of Theories}
 \begin{theindex}
 \addcontentsline{toc}{chapter}{Index of Theories}

\glohead{Languages}
\gloitem {\LC }{3.1.1}
\gloitem {\LCS }{3.1.1}
\gloitem {\LS }{3.1.1}
\gloitem {\LN }{6.2.1}

\indexspace
\glohead{Predicates}
\gloitem {$\tABST(x)$}{S1.1}
\gloitem {$\tATOM(x)$}{3.1.1}
\gloitem{ $\Divn(x)$ }{3.2.33}
\gloitem {$\tINDISC(x,y)$}{2.1.1}
 \gloitem {$\Modn{m}(x) $}{3.2.33}
 \gloitem {$\tSm(x,y)$}{4.3.3}
 \gloitem {$\tSUM(x,y,z)$}{3.1.1}
 \gloitem {$\tUNIT(x)$}{3.1.1}
 \gloitem {$\Timesn(x,y)$}{3.2.33}

\indexspace
\glohead{Theories}
\gloitem{\ADIVn }{3.2.33}
\gloitem {\ASSOC}{4.1.3}
\gloitem {\ASSOCT}{4.3.3}
\gloitem {\ASYMLT}{2.1.1}
\gloitem {\ASYMF}{2.1.1}
\gloitem {\ASYMGT}{2.1.1}
\gloitem {\ATLEASTn }{3.2.2}
\gloitem {\BA }{3.2.1}
\gloitem {\BAn }{3.2.2}
\gloitem {\BAI }{3.2.2}
\gloitem {\BASIC }{3.2.7}

\gloitem {\BDIVJ }{3.2.11}

\gloitem {\BDIVn }{3.2.11}

\gloitem {\CA}{3.2.37}
\gloitem {\CAI}{3.3.6}
\gloitem {\CAIf}{3.3.7}
\gloitem {\CORE }{2.1.2}
\gloitem {\CS }{3.1.4}
\gloitem {\CSI }{3.3.6}

\gloitem {\CANTORLT}{S2.1}
\gloitem {\DEFPL}{2.2.1}
 \gloitem {\DEFPL}{3.2.6}
\gloitem {\DEFEQ}{2.1.2}
\gloitem {\DEFGT}{2.1.1}
\gloitem {\DISJU}{3.2.6}
\gloitem {\DISJPL}{2.2.1}

\gloitem{\FUNCPL}{2.2.1}
\gloitem {\DIVJ }{3.2.11}

 \gloitem{ \DIVn}{3.2.8}

\gloitem {\EVEN }{S2.1}
\gloitem {\EXn }{3.2.2}
\gloitem {\EXCORE }{2.2.4}

\gloitem {\INDISCEQ}{2.1.1}
\gloitem {\IRREFLT}{2.1.2}

\gloitem {\INF }{3.2.2}

\gloitem {\MAX}  {4.1.3}
\gloitem {\MIN}  {4.1.3}
\gloitem {\MODn{m}}{3.2.8}
\gloitem {\MONOT}{2.2.1}
\gloitem {\ODD   }{S2.1}

\gloitem {\OP }{6.1.1}
\gloitem {\PCA }{4.1.6}
\gloitem {\PCS }{4.1.2}
\gloitem {\PSIZE }{4.1.3}
 \gloitem {\QUASI }{2.1.1}
 
 \gloitem {\REFEQ}{2.1.1}
 
\gloitem {\REPLT}{S2.1}
\gloitem {\RTf }{3.3.7}
\gloitem {\SIZE }{3.2.6}

\gloitem {\SYMEQ}{2.1.1}
 
 \gloitem {\Tf }{3.3.7}
\gloitem {\TRANSEQ}{2.1.1}
\gloitem {\TRANSLT}{2.1.1}
\gloitem {\TRANSGT}{2.1.1}
\gloitem {\TRICH}{2.1.2}
\gloitem {\UNIQUEMN}  {4.3.3}
\gloitem {\UNIQUEPL}  {4.3.3}
\gloitem {\UNIQEQ}  {4.1.3}

\gloitem {\Zg }{4.2.10}
\gloitem {\Zgm }{4.2.10}

\indexspace

\end{theindex}
 

\renewcommand\indexname{Subject Index}
 \addcontentsline{toc}{chapter}{Subject Index}
\printindex

\end{document}